\documentclass[]{interact}

\usepackage{epstopdf}
\usepackage[caption=false]{subfig}

\usepackage[natbibapa,nodoi]{apacite}
\renewcommand\bibliographytypesize{\fontsize{10}{12}\selectfont}

\theoremstyle{plain}
\newtheorem{theorem}{Theorem}[section]
\newtheorem{lemma}[theorem]{Lemma}
\newtheorem{corollary}[theorem]{Corollary}

\theoremstyle{definition}
\newtheorem{definition}[theorem]{Definition}

\theoremstyle{remark}

\usepackage{graphicx} 
\usepackage[space]{grffile}
\usepackage{amsmath, amsthm, amssymb, gensymb}
\usepackage{tcolorbox,colortbl}
\usepackage{enumitem}
\usepackage{tikz}
\usetikzlibrary{shapes.geometric}
\usetikzlibrary{decorations.pathreplacing}
\usepackage{hyperref}
\usepackage{url}
\usepackage{multirow}

\definecolor{glideblue}{HTML}{05ADE2}
\definecolor{refpink}{HTML}{DC2780}
\definecolor{rotyel}{HTML}{FCB017}

\makeatletter
\def\squarecorner#1{
    \pgf@x=\the\wd\pgfnodeparttextbox%
    \pgfmathsetlength\pgf@xc{\pgfkeysvalueof{/pgf/inner xsep}}%
    \advance\pgf@x by 2\pgf@xc%
    \pgfmathsetlength\pgf@xb{\pgfkeysvalueof{/pgf/minimum width}}%
    \ifdim\pgf@x<\pgf@xb%
        \pgf@x=\pgf@xb%
    \fi%
    \pgf@y=\ht\pgfnodeparttextbox%
    \advance\pgf@y by\dp\pgfnodeparttextbox%
    \pgfmathsetlength\pgf@yc{\pgfkeysvalueof{/pgf/inner ysep}}%
    \advance\pgf@y by 2\pgf@yc%
    \pgfmathsetlength\pgf@yb{\pgfkeysvalueof{/pgf/minimum height}}%
    \ifdim\pgf@y<\pgf@yb%
        \pgf@y=\pgf@yb%
    \fi%
    \ifdim\pgf@x<\pgf@y%
        \pgf@x=\pgf@y%
    \else
        \pgf@y=\pgf@x%
    \fi
    \pgf@x=#1.5\pgf@x%
    \advance\pgf@x by.5\wd\pgfnodeparttextbox%
    \pgfmathsetlength\pgf@xa{\pgfkeysvalueof{/pgf/outer xsep}}%
    \advance\pgf@x by#1\pgf@xa%
    \pgf@y=#1.5\pgf@y%
    \advance\pgf@y by-.5\dp\pgfnodeparttextbox%
    \advance\pgf@y by.5\ht\pgfnodeparttextbox%
    \pgfmathsetlength\pgf@ya{\pgfkeysvalueof{/pgf/outer ysep}}%
    \advance\pgf@y by#1\pgf@ya%
}
\makeatother

\pgfdeclareshape{square}{
    \savedanchor\northeast{\squarecorner{}}
    \savedanchor\southwest{\squarecorner{-}}

    \foreach \x in {east,west} \foreach \y in {north,mid,base,south} {
        \inheritanchor[from=rectangle]{\y\space\x}
    }
    \foreach \x in {east,west,north,mid,base,south,center,text} {
        \inheritanchor[from=rectangle]{\x}
    }
    \inheritanchorborder[from=rectangle]
    \inheritbackgroundpath[from=rectangle]
}

\tikzset{%
    2fold/.style={thick,diamond,fill=black,draw,minimum width = 1mm,minimum height = 4mm,inner sep=0.5mm},
    2foldcolorswap/.style={thick,diamond,black,fill=red,draw,minimum width = 1mm,minimum height = 4mm,inner sep=0.5mm},
    4fold/.style={thick,square,fill=black,draw,minimum width = 2.5mm,inner sep=0.5mm},
    4foldcolorswap/.style={thick,square,black,fill=red,draw,minimum width = 2.5mm,inner sep=0.5mm},
    reddot/.style={circle,red,fill=red,draw,inner sep=0.5mm},
}

\title{A Whole New Side: The Duality Symmetries of Hitomezashi}
\author{Megan A. Martinez}

\begin{document}

\title{A Whole New Side: The Duality Symmetries of Generalized Hitomezashi Patterns}

\author{
\name{Megan A.~Martinez\textsuperscript{a}\thanks{CONTACT Megan A.~Martinez. Email: mmartinez@ithaca.edu}}
\affil{\textsuperscript{a}Department of Mathematics, Ithaca College, Ithaca, NY, USA}
}

\maketitle

\begin{abstract}
Hitomezashi is a particular style of sashiko stitching characterized by crossing lines of stitches. Certain hitomezashi designs can be encoded with two binary strings and have been the subject of a wide-array of mathematical research. Dubbed generalized hitomezashi patterns (GHPs) by Sen \& Martinez, these designs are reversible and can be `self-dual.' We investigate the possible `dual-' or `flip-symmetries' in GHPs; the inclusion of these symmetries is equivalent to considering two-colour symmetry groups. Following the work of Sen \& Martinez, who found the wallpaper symmetry groups that are compatible with GHPs, this paper investigates the possible two-colour symmetries compatible with GHPs. We prove that there are four possible rosette two-colour symmetry types, completely describe the properties on binary strings that create these symmetries, and prove that exactly 17 of the 46 two-colour wallpaper groups are compatible with GHPs.
\end{abstract}


\section{Introduction}

Sashiko embroidery originated in Japan during the Edo era (1615-1868). It began as a way to mend worn fabric by sandwiching together multiple layers of fabric and stitching across the entire work. As the craft developed, it became a way to make beautiful designs on clothing without running afoul of the sumptuary laws of the time. In the modern day, sashiko has found its way into the crafting culture---it is fairly easy to find authentic sashiko embroidery materials and helpful videos through the magic of the internet. For those interested in reading more about the history, evolution, and techniques of sashiko, we point the reader to the wonderful book of Susan Briscoe \citeyearpar{Briscoe}.

There are a few different styles of sashiko stitching. One popular form of stitching is the \emph{moyozashi} stitching where lines of running stitch wander all over the fabric to make lovely geometric designs. Another type (and the focus of this paper) is \emph{hitomezashi} stitching; this kind of stitching is characterized by running stitches where threads cross each other and stitching is generally done along straight lines \citep{Briscoe}.

\begin{figure}[h!tbp]
\centering
\includegraphics[width=1.75in]{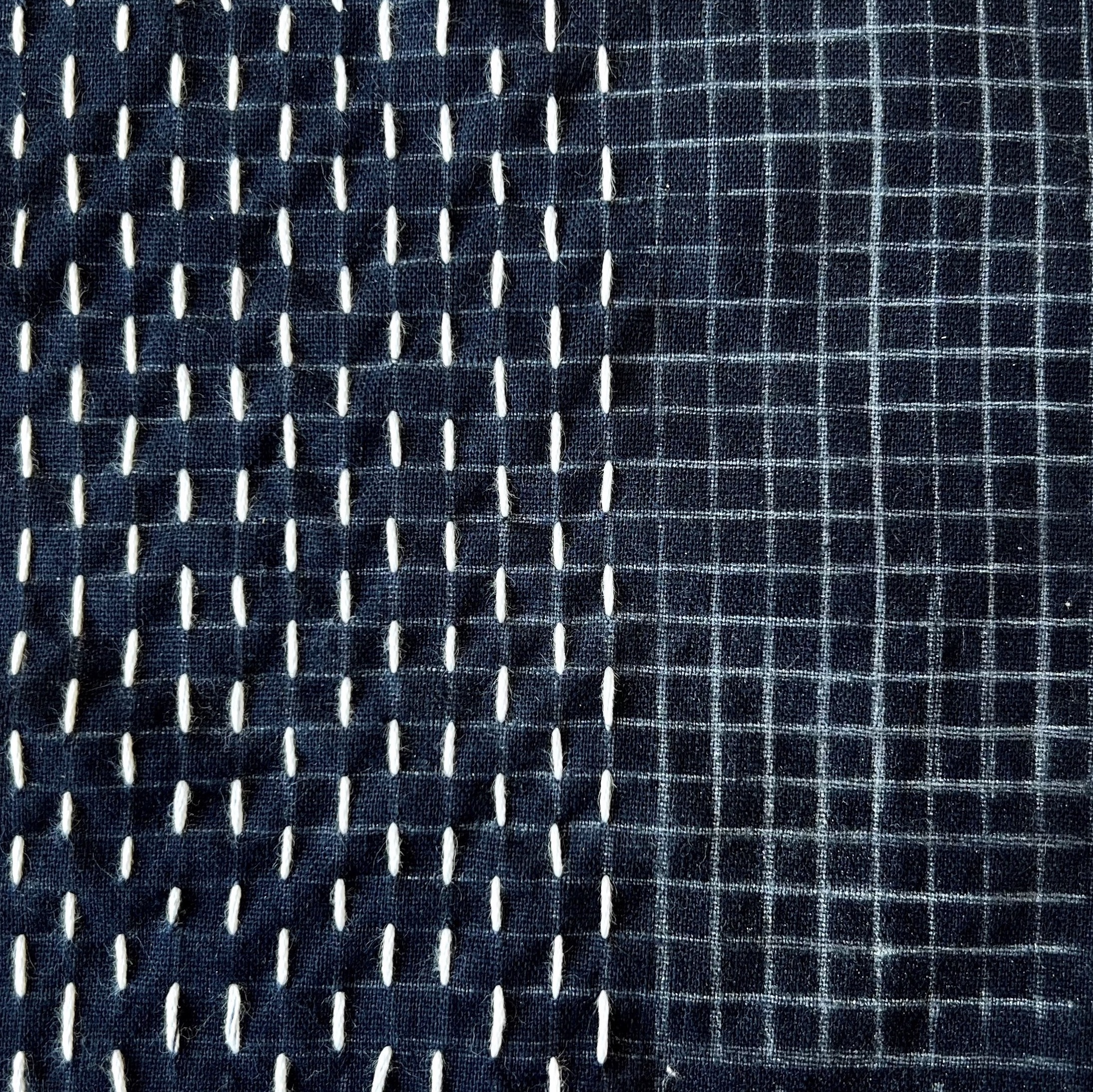} \quad
\includegraphics[width=1.75in]{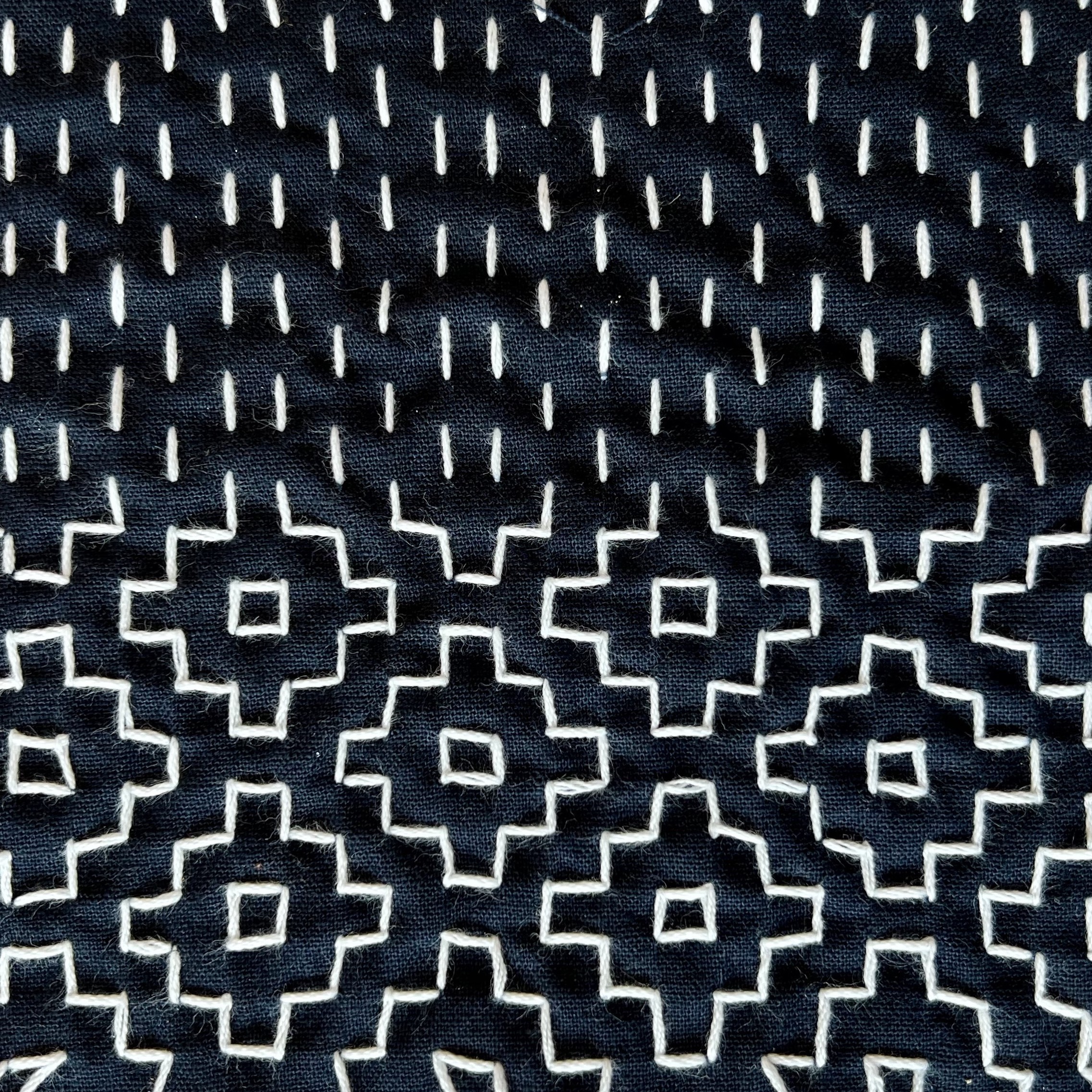} \quad
\includegraphics[width=1.75in]{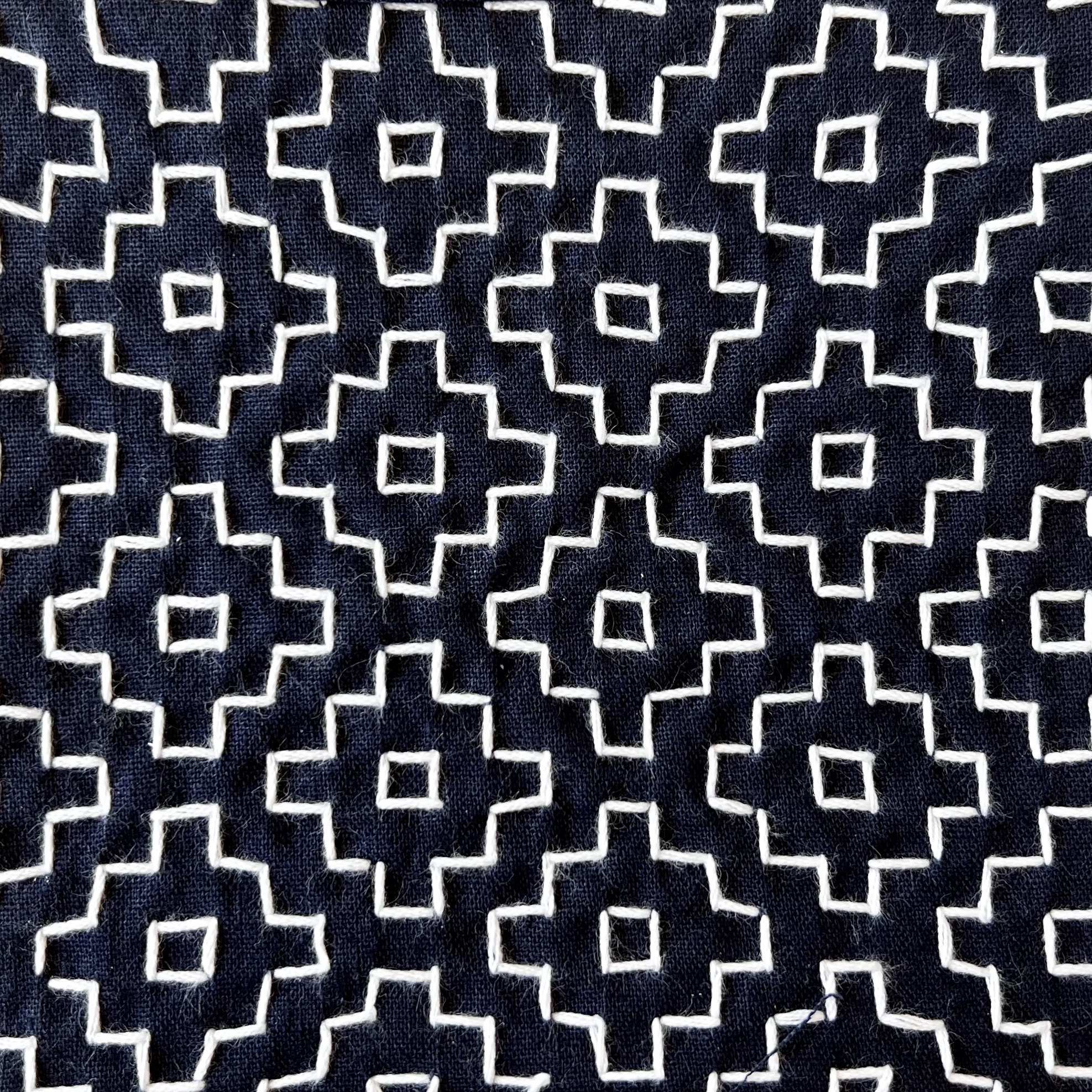}
\caption{The stages of stitching the kakinohanazashi, or persimmon flower, stitch is shown here. In this example, the vertical lines of stitches are stitched first, then the horizontal lines of stitches \citep{Sen_Martinez}.} \label{Process}
\end{figure}

A particular style of hitomezashi, which can be encoded with two binary strings, has come to the attention of mathematicians because of its relatively simple construction that yields intricate results. For these hitomezashi patterns, stitching is done on a square grid in rows and columns (one stitch is the side length of a square), and along any line of stitching we use a running stitch where the thread alternates above and below the fabric. Whether you start your stitching above or below the fabric determines the entire line of stitches. The steps of an in-progress hitomezashi are given in Figure~\ref{Process}.

\begin{figure}[h!tbp]
\centering
\subfloat[]{%
\resizebox*{6.4cm}{!}{\includegraphics{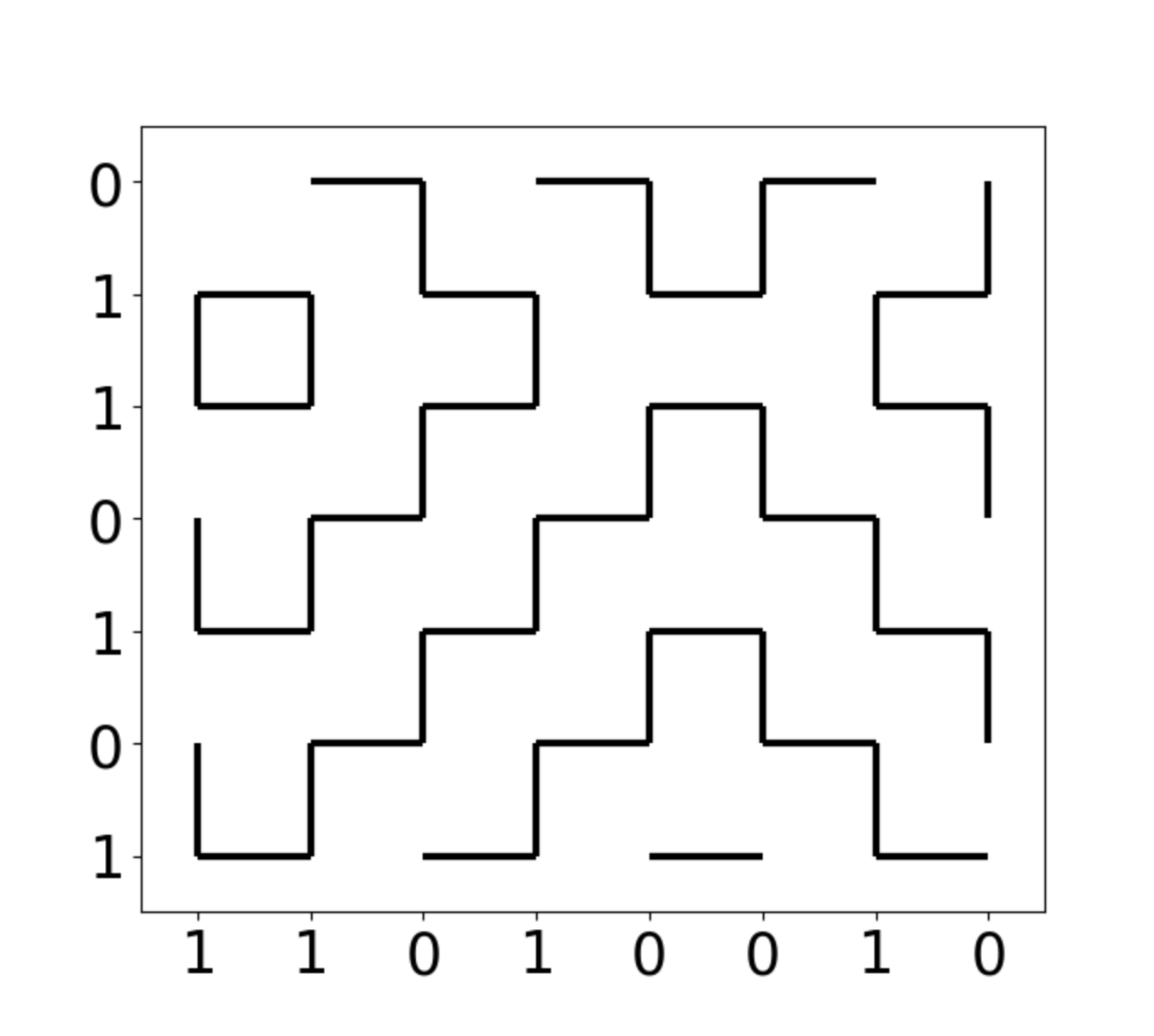}}} \hspace{5pt}
\subfloat[]{%
\resizebox*{5.75cm}{!}{\includegraphics{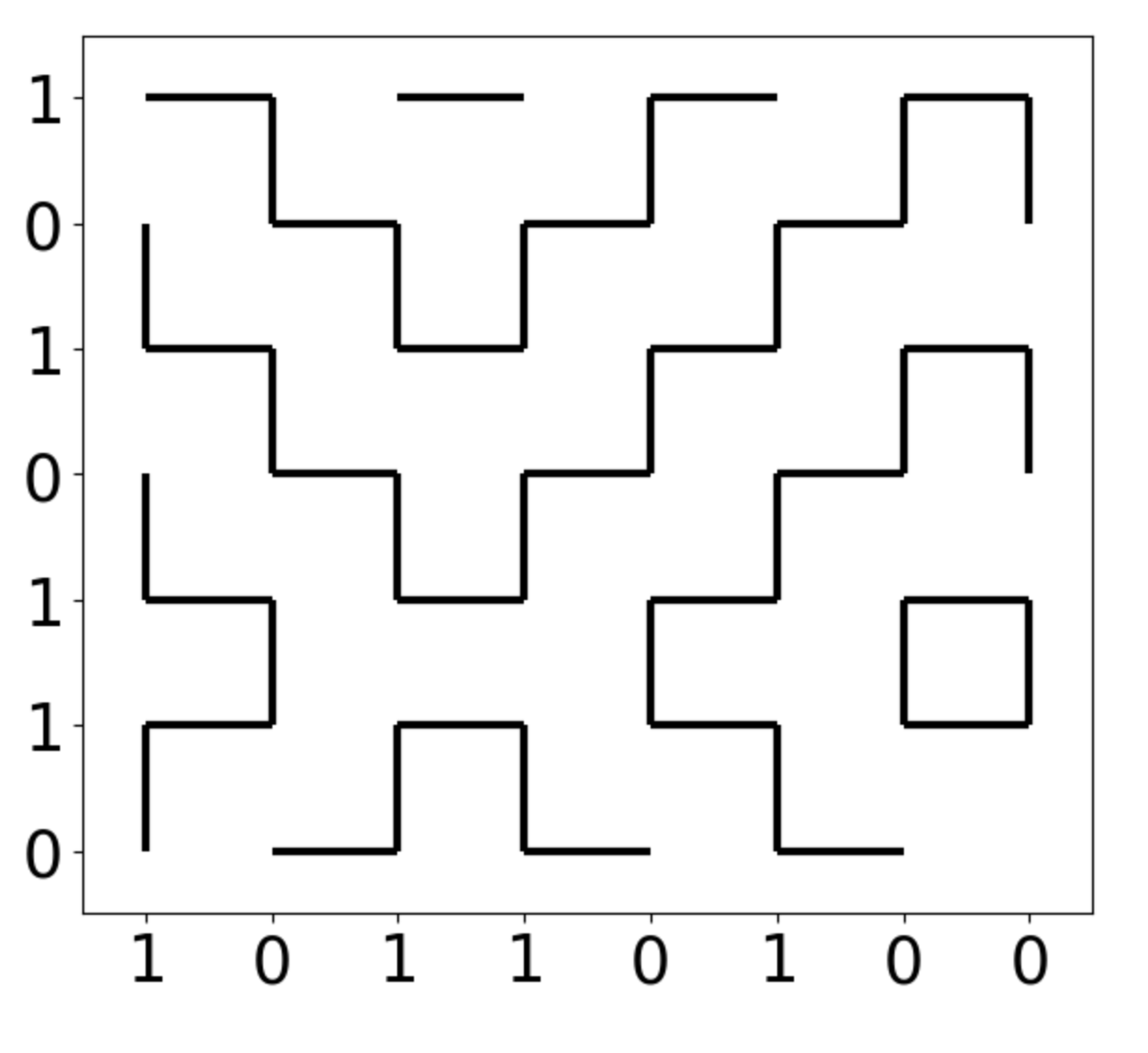}}}
\caption{Figure (a) shows a GHP generated by $(11010010,1010110)$; Figure (b) shows the same GHP rotated by $180^{\circ}$ which is generated by $((11010010)^{RC},(1010110)^R)=(10110100,0110101)$.} \label{first_example}
\end{figure}

These kinds of hitomezashi patterns can be encoded using two binary strings, where a `0' indicates that stitching starts below the fabric and `1' indicates the stitching starts above the fabric. One binary string is associated with the bottom of the design and determines the vertical rows of stitches; the other is associated with the left of the design and determines the horizontal rows of stitches. Two examples of this correspondence are provided in Figure~\ref{first_example}.

Following the convention of Sen \& Martinez \citeyearpar{Sen_Martinez} we refer to any hitomezashi design generated by binary strings as a \emph{generalized hitomezashi patterns} (GHPs for short). This recognizes the fact that we are expanding our array of designs beyond the traditional hitomezashi patterns.

Generalized hitomezashi patterns were brought to the attention of the mathematical community by \cite{Hayes_Seaton} and gained broader popularity through a \emph{Numberphile} YouTube video inspired by Seaton's work \citep{numberphile}. Since these introductions, mathematicians have investigated a great deal about the structure of hitomezashi. Hayes and Seaton have pushed on the math and art connection by recognizing self-dual Fibonacci Snowflakes in the designs \citeyearpar{seaton_hayes2}. Seaton has also explored results on D-forms and crafting biscornu using hitomezashi \citeyearpar{Seaton}. 

The regions present in GHPs have also been the subject of interest. Defant and Kravitz found properties of the loops and regions formed in hitomezashi \citeyearpar{Defant_Kravitz} (with a shortened proof presented by \cite{Ren_Zhang2}). In this work, they also recognized that hitomezashi patterns make an appearance in work studying `corner percolation configurations,' which are equivalent to randomly-generated hitomezashi patterns. Research has been done with these configurations examining the expected values of statistics related to cycles and paths \citep{Pete, Marchand}. Further work investigating the loops and regions in hitomezashi has been applied to different stitch variants \citep{Defant_Kravitz_Tenner} and carried out on different surfaces \citep{Ren_Zhang, Xie}.

\begin{figure}[h!tbp]
\centering
\subfloat[Complete symmetry sampler with fronts pictured]{%
\resizebox*{10cm}{!}{
\includegraphics{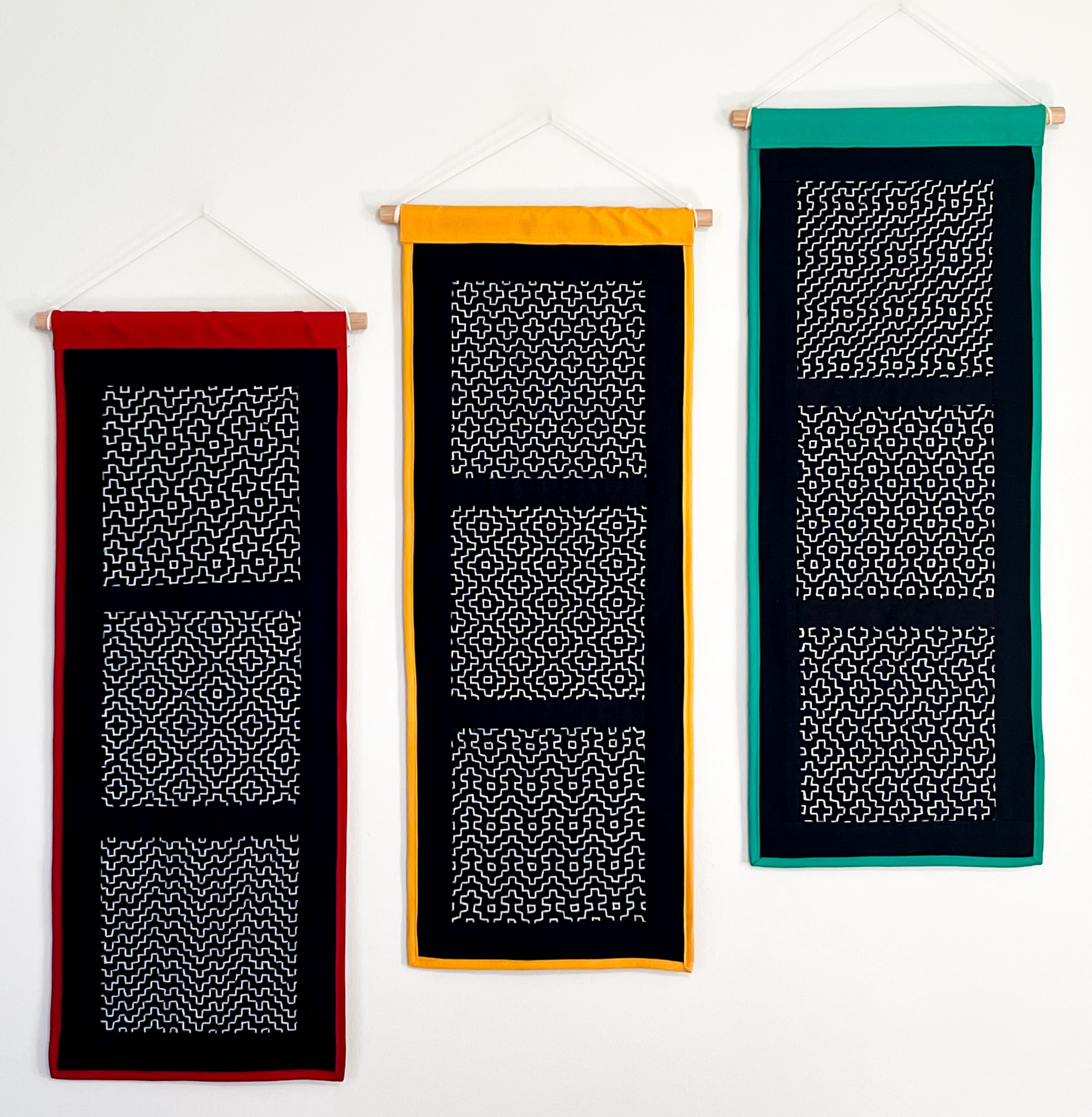}
}}

\subfloat[The front of the green panel]{%
\resizebox*{4cm}{!}{
\includegraphics{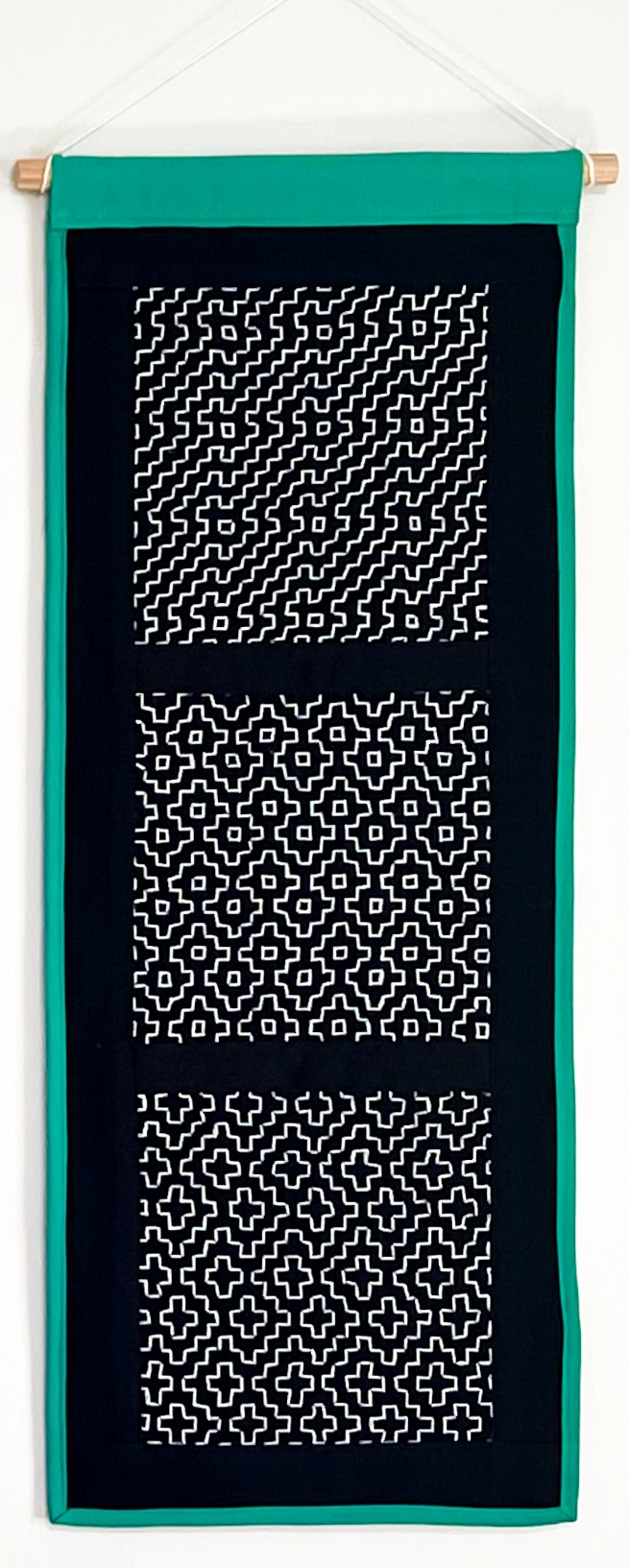}
}} \hspace{1.5cm}
\subfloat[The back of the green panel]{%
\resizebox*{4cm}{!}{
\includegraphics{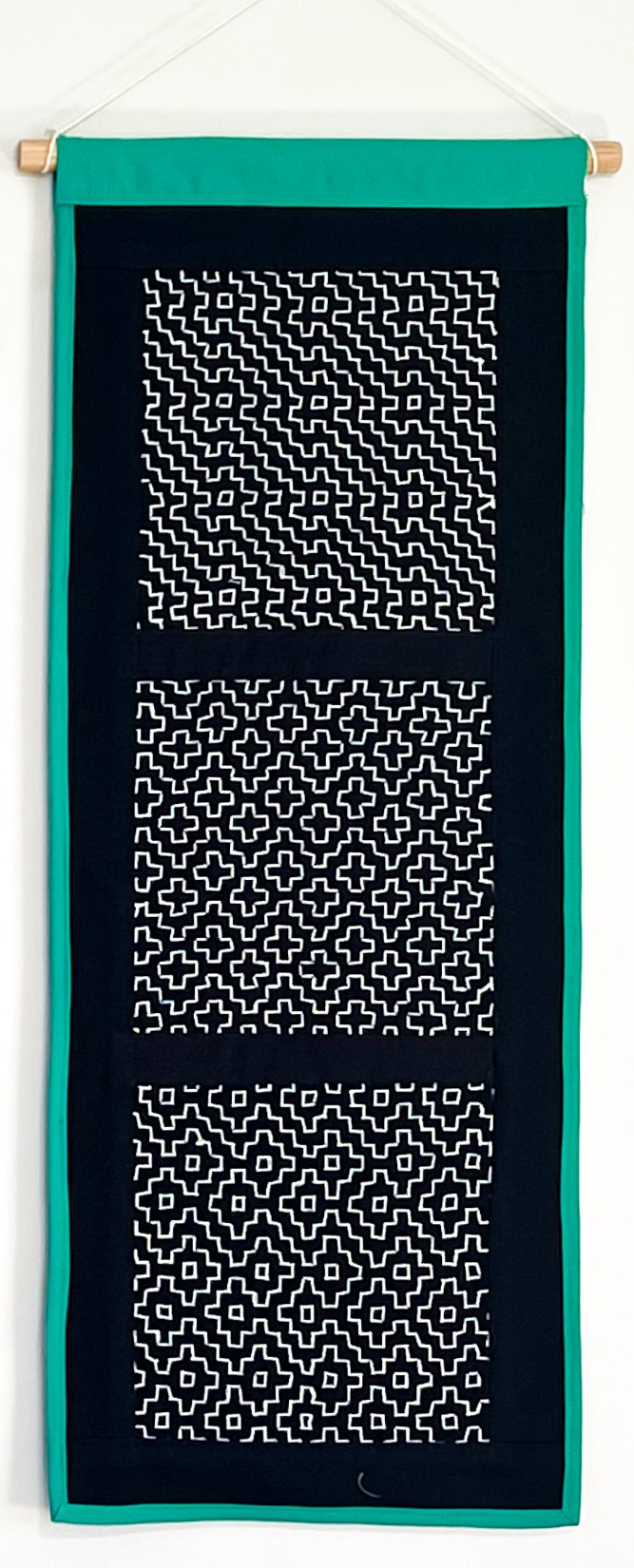}
}}

\caption{Images of the \emph{Hitomezashi Wallpaper Triptych} exhibited at the Bridges 2024 Exhibition of Mathematical Art, Craft, and Design}
\label{fig:wallpaper_triptych}
\end{figure}

As GHPs are a `constrained fiber art form' as described by \cite{Goldstine}, determining the kinds of symmetries that can occur in GHPs is a mathematically-rich endeavour as well. Sen \& Martinez proved that 9 of the 17 wallpaper groups are compatible with GHPs \citep{Sen_Martinez}. Examples of the compatible wallpaper symmetry groups were stitched into a `symmetry sampler' by this author and displayed in the Bridges 2024 Exhibition of Mathematical Art, Craft, and Design \citep{Wallpaper_Triptych}. The front of each of the three panels are presented in Figure~\ref{fig:wallpaper_triptych}(a). The symmetries of GHPs adapted to an isometric grid have also been investigated by Seaton; in this work, particularly beautiful results arise when the stitching is done in a `dilute' way (i.e.\ every other line is stitched) \citep{Seaton_Isometric}.

Around the same time that Sen \& Martinez published their work on symmetry, Seaton released a book on hitomezashi covering a wide array of topics; included are the notion of duality (Chapter 4) and a brief overview of symmetries (Chapter 6) \citep{Seaton_Book}. The dual nature of GHPs is particularly interesting, as it is easy to forget, as we admire hitomezashi designs on the pages of this paper (or a computer screen), that the back of our design is just as beautiful as the front. This was demostrated by Seaton and Hayes as they designed `self-dual' pell persimmons \citep{seaton_hayes2}. It can also be seen in the \emph{Hitomezashi Wallpaper Triptych}, which was designed to be reversible. Figures~\ref{fig:wallpaper_triptych}(b),(c) show the front and back of the green panel. Notice that the top wallpaper in this panel is one that looks similar on the front and back; they appear to be reflections of each other. Meanwhile, the bottom two wallpapers are distinct designs. The top wallpaper is an example of a `self-dual' GHP.

The work of this paper is a natural sequel to the symmetry work of Sen \& Martinez. We will combine the ideas of duality and symmetry to systematically characterize the possible `dual' or `flip'-symmetry groups that can occur in GHPs. We will relate dual symmetries to two-colour symmetry groups; this connection will allow us to use already explored two-colour rosette and wallpaper symmetry groups to describe our results. 

Seaton recently investigated the dual symmetries compatible with GHPs that are periodic in one direction and form frieze patterns \citep{Seaton_Friezes}. Instead of doing this work by investigating two-colour frieze groups, she considered the 'two-sided' frieze symmetry groups introduced by \cite{Cromwell} that are compatible with hitomezashi designs. This is a natural way to analyze the friezes, as the two-sided frieze groups were devised with the intention of analyzing a frieze that can physically be turned over. A certain subset of two-sided friezes are analogous to two-colour friezes and it transpires that each of the possible two-sided frieze symmetries compatible with hitomezashi can actually be described as a two-colour frieze group. In this way, the results of this paper, in addition to Seaton's results, form a complete picture of the dual symmetries compatible with GHPs.

The main results of this paper describe the two-colour rosette groups and two-colour wallpaper groups that are compatible with the symmetries of a GHP. In addition, the properties of binary strings that produce particular rosette symmetries will be laid out in detail. In Section \ref{past_results}, we include the essential results of Sen \& Martinez that are needed for our work. Section \ref{flip_symm_rosettes} describes the connection between `flip-symmetries' and two-colour symmetries in addition to working out the compatible two-colour rosette symmetries. Section \ref{wallpaper} then investigates which two-colour wallpaper groups are compatible with GHPs, proving that 17 of the 46 groups can be realised. Examples of each of the 17 two-colour wallpaper symmetries are also included in Section \ref{wallpaper}.

\section{Rosette Symmetries} \label{past_results}

We will heavily utilize results of Sen and Martinez \citeyearpar{Sen_Martinez}, so will use this section to introduce the important definitions and results that form the building blocks for the rest of the paper. In Section~\ref{subsection:rosette}, we will collate the Sen \& Martinez results to completely describe the rosette symmetries that are compatible with GHPs.

\subsection{Notation}

Since GHPs are generated by two binary strings, we must introduce notation related to these strings that will be used throughout the paper. We consider two binary strings $x=x_1\ldots x_n$ and $y=y_1\ldots y_m$, where $x$ is placed along the bottom and defines the vertical stitches and $y$ is placed along the left and defines the horizontal stitches; both words start in the lower left-hand corner. This yields a design that measures $n-1$ by $m-1$. For example, the design generated by $(x,y)=(11010010,1010110)$ is pictured in Figure~\ref{first_example}(a).

We will be able to understand the symmetries of a design solely from the generating binary strings. It is useful to describe some properties we can look for on binary strings to understand the resulting symmetries.

\begin{definition}
Given a binary string $a=a_1\ldots a_n \in \{0,1\}^n$ we define the \emph{reverse} as $a^R=a_n a_{n-1} \ldots a_1$ and the \emph{complement} as $a^C = (1-a_1)\ldots (1-a_n)$. Additionally, the \emph{reverse-complement} is $a^{RC} = (1-a_n)(1-a_{n-1})\ldots (1-a_1)$. 
\end{definition}

\noindent For example, when $a=1011001$, we have $a^R=1001101$, $a^C=0100110$, and $a^{RC}=0110010$. 

\subsection{The symmetry results of Sen \& Martinez}

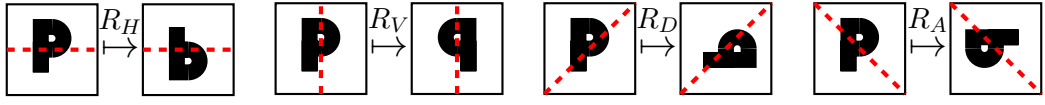
\begin{figure}[h!tbp] 
\centering
\begin{tikzpicture}[scale=0.3]
\draw[thick] (0,0) rectangle (4,4);
\node at (5,2) {\Large$\mapsto$};
\node at (5,3.25) {$R_H$};
\draw[line width=6pt, black] (2,2)--(1.5,2)--(1.5,1)--(1.5,3)--(2,3);
\draw[line width=6pt, black] (2,2) arc (-90:90:0.5);
\draw[ultra thick,dashed,red] (0,2)--(4,2);
\draw[thick] (6,0) rectangle (10,4);
\draw[line width=6pt, black] (8,1)--(7.5,1)--(7.5,2)--(8,2)--(7.5,2)--(7.5,3);
\draw[line width=6pt, black] (8,1) arc (-90:90:0.5);
\draw[ultra thick,dashed,red] (6,2)--(10,2);
\end{tikzpicture}
\quad
\begin{tikzpicture}[scale=0.3]
\draw[thick] (0,0) rectangle (4,4);
\node at (5,2) {\Large$\mapsto$};
\node at (5,3.25) {$R_V$};
\draw[line width=6pt, black] (2,2)--(1.5,2)--(1.5,1)--(1.5,3)--(2,3);
\draw[line width=6pt, black] (2,2) arc (-90:90:0.5);
\draw[ultra thick,dashed,red] (2,0)--(2,4);
\draw[thick] (6,0) rectangle (10,4);
\draw[line width=6pt, black] (8,2)--(8.5,2)--(8.5,1)--(8.5,3)--(8,3);
\draw[line width=6pt, black] (8,3) arc (90:270:0.5);
\draw[ultra thick,dashed,red] (8,0)--(8,4);
\end{tikzpicture}
\quad
\begin{tikzpicture}[scale=0.3]
\draw[thick] (0,0) rectangle (4,4);
\node at (5,2) {\Large$\mapsto$};
\node at (5,3.25) {$R_D$};
\draw[line width=6pt, black] (2,2)--(1.5,2)--(1.5,1)--(1.5,3)--(2,3);
\draw[line width=6pt, black] (2,2) arc (-90:90:0.5);
\draw[ultra thick,dashed,red] (0,0)--(4,4);
\draw[thick] (6,0) rectangle (10,4);
\draw[line width=6pt, black] (8,2)--(8,1.5)--(7,1.5)--(9,1.5)--(9,2);
\draw[line width=6pt, black] (9,2) arc (0:180:0.5);
\draw[ultra thick,dashed,red] (6,0)--(10,4);
\end{tikzpicture}
\quad
\begin{tikzpicture}[scale=0.3]
\draw[thick] (0,0) rectangle (4,4);
\node at (5,2) {\Large$\mapsto$};
\node at (5,3.25) {$R_A$};
\draw[line width=6pt, black] (2,2)--(1.5,2)--(1.5,1)--(1.5,3)--(2,3);
\draw[line width=6pt, black] (2,2) arc (-90:90:0.5);
\draw[ultra thick,dashed,red] (0,4)--(4,0);
\draw[thick] (6,0) rectangle (10,4);
\draw[line width=6pt, black] (8,2)--(8,2.5)--(9,2.5)--(7,2.5)--(7,2);
\draw[line width=6pt, black] (7,2) arc (180:360:0.5);
\draw[ultra thick,dashed,red] (6,4)--(10,0);
\end{tikzpicture}
\caption{The four different reflections in $D_4$.} \label{reflections}
\end{figure}

When considering the attainable symmetries of a GHP, we only need to consider possible symmetries of a square; that is, the group $D_4$. These are the counterclockwise rotations by $0^{\circ}$, $90^{\circ}$, $180^{\circ}$, and $270^{\circ}$ notated as $e$, $r_{90}$, $r_{180}$, and $r_{270}$, respectively; and reflections $R_H$, $R_V$, $R_D$, and $R_A$ as pictured in Figure~\ref{reflections}. Certain symmetries, such as $R_D$, $R_A$, $r_{90}$, and $r_{270}$ require a square design where $|x|=n=m=|y|$; otherwise, we do not enforce this restriction.

Sen \& Martinez presented results that determined the symmetries of a GHP from only the binary strings. This requires determining how each isometry in $D_4$ will alter the binary strings of a design. These are described in the following two lemmas.

\begin{lemma} \label{lem_reflections}
\citep{Sen_Martinez} The reflections $R_H$, $R_V$, $R_D$, and $R_A$ have the following effects on the binary strings that generate a GHP:
\begin{center}
$
R_H(x, y) =
    \begin{cases}
        (x, y^R) & \text{$m$ even}\\
        (x^C, y^R) & \text{$m$ odd}
    \end{cases}
$
\hspace{1cm}
$
R_V(x, y) =
    \begin{cases}
        (x^R, y) & \text{$n$ even}\\
        (x^R, y^C) & \text{$n$ odd}
    \end{cases}
$

$
    R_D(x, y) = (y, x)
$
\hspace{1cm}
$
R_A(x, y) =
    \begin{cases}
        (y^R, x^R) & \text{$n$ even}\\
        (y^{RC},x^{RC}) & \text{$n$ odd}.
    \end{cases}
$
\end{center}
\end{lemma}

\begin{lemma} \label{lem_rotations}
The counterclockwise rotations $r_{90}$, $r_{180}$, and $r_{270}$ have the following effects on the binary strings that generate a GHP:

\begin{minipage}[h]{0.45\textwidth}

$
    r_{90}(x, y) =
    \begin{cases}
        (y^R, x) & \text{$m$ even}\\
        (y^R, x^C) & \text{$m$ odd}
    \end{cases}
$

\vspace{4mm}
$
    r_{270}(x, y) =
    \begin{cases}
        (y, x^R) & \text{$n$ even}\\
        (y^C, x^R) & \text{$n$ odd}
    \end{cases}
$
\end{minipage}
\begin{minipage}[h]{0.5\textwidth}
$
    r_{180}(x, y) =
    \begin{cases}
        (x^R, y^R) & \text{$n, m$ even}\\
        (x^{RC},y^R) & \text{$n$ even, $m$ odd}\\
        (x^R,y^{RC}) & \text{$n$ odd, $m$ even}\\
        (x^{RC},y^{RC}) & \text{$n,m$ odd}
    \end{cases}
$
\end{minipage}

\end{lemma}

Figure~\ref{first_example}(b) shows the effect of applying $r_{180}$ when $n$ is even and $m$ is odd. 

To determine when one of the isometries above is a symmetry, we need properties that must hold for $x$ and $y$ so that, for some isometry $I$, $I(x,y)=(x,y)$. Note that it is entirely possible for $x=x^R$; this occurs when $x$ is a palindrome. However, there are some conditions that can never occur.

\begin{lemma} \label{basic_facts}
\citep{Sen_Martinez} Let $a \in \{0,1\}^n$. Then for all $n$, $a \neq a^C$. Additionally, if $n$ is odd, then $a \neq a^{RC}$.
\end{lemma}

When you combine the results of Lemmas~\ref{lem_reflections}, \ref{lem_rotations}, and \ref{basic_facts} you can determine the symmetries of a GHP based only on the generating binary strings using the following theorem.

\begin{theorem}\label{ros_symm}
\citep{Sen_Martinez} Consider a GHP generated by binary strings $(x,y)$ where $|x|=n$ and $|y|=m$. Then:

\begin{itemize}
\item $R_H$ (resp.\ $R_V$) is a symmetry if and only if $y=y^R$ (resp.\ $x=x^R$) and $m$ (resp.\ $n$) is even
\item $R_D$ is a symmetry if and only if $n=m$ and $x=y$
\item $R_A$ is a symmetry if and only if $x=y^R$ when $n=m$ is even and $x=y^{RC}$ when $n=m$ is odd
\item $r_{90}$ and $r_{270}$ are symmetries if and only if $n=m$ is even and $x=x^R=y=y^R$
\item $r_{180}$ is a symmetry if and only if one of the following is true: $n$ and $m$ are even, $x=x^R$, and $y=y^R$; $n$ is odd, $m$ is even, $x=x^R$, and $y=y^{RC}$; or $n$ is even, $m$ is odd, $x=x^{RC}$, and $y=y^R$
\end{itemize}
\end{theorem}

From Theorem~\ref{ros_symm}, Sen \& Martinez were able to extract two important facts about the occurrence of rotational symmetries in GHPs that will be essential for our work going forward.

\begin{corollary} \label{cor:4Fold}
\citep{Sen_Martinez} Consider a point $P$ in a GHP (with a potentially unbounded motif).
\begin{enumerate}
\item $P$ is a center of 4-fold rotational symmetry if and only if $P$ lies at the intersection of 4 lines of reflective symmetry.
\item $P$ is a center of 2-fold rotational symmetry that is not on gridlines if and only if $P$ lies at the intersection of at least 2 lines of reflective symmetry.
\end{enumerate}
\end{corollary}

\begin{figure}[h!tbp]
\centering
\subfloat[Two-fold rotational symmetry with lines of reflection]{%
\resizebox*{4cm}{!}{\includegraphics{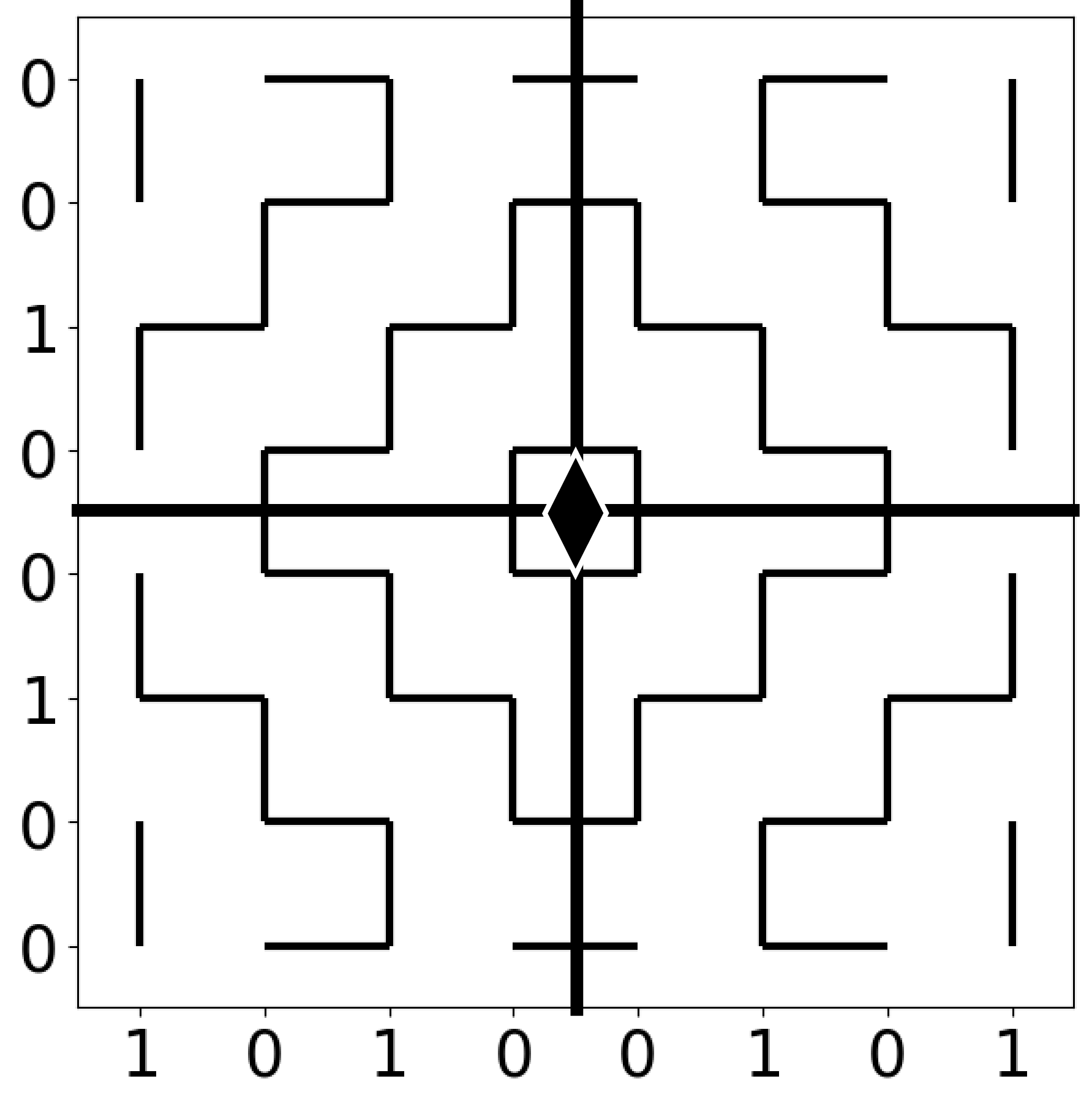}}}\hspace{5pt}
\subfloat[Two-fold rotational symmetry without lines of reflection]{%
\resizebox*{4cm}{!}{\includegraphics{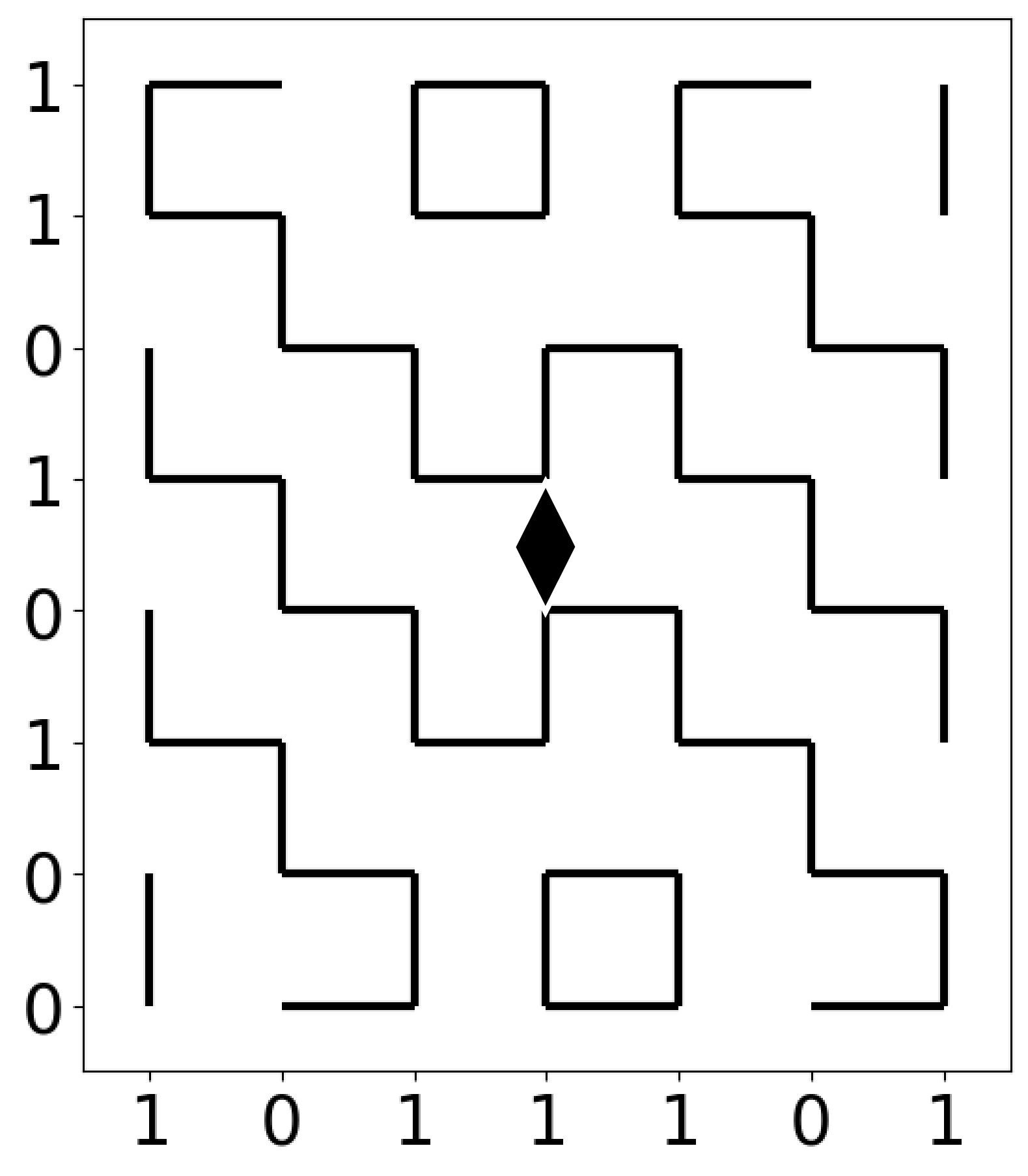}}}\hspace{5pt}
\subfloat[Four-fold rotational symmetry]{%
\resizebox*{4cm}{!}{\includegraphics{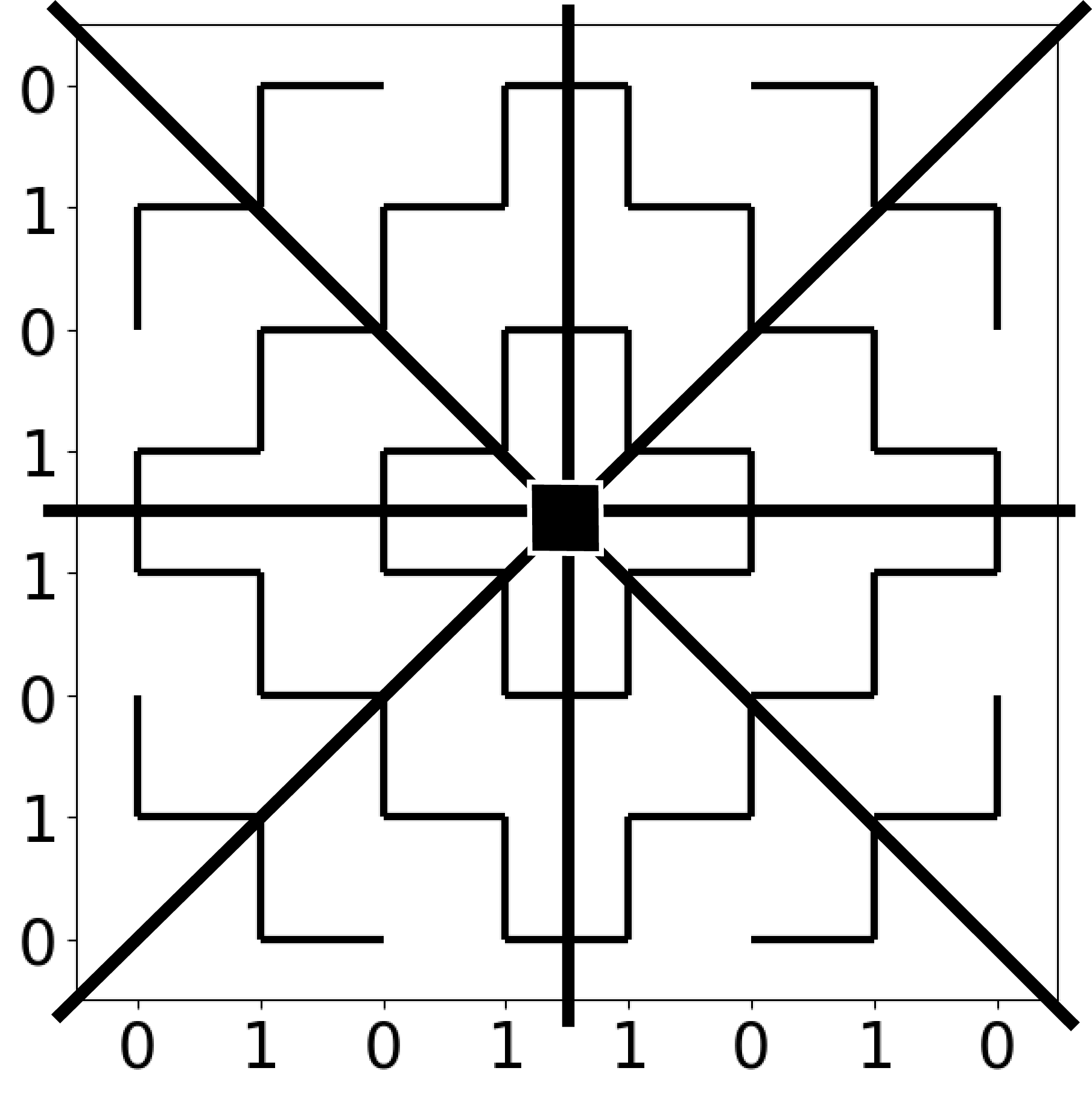}}}
\caption{Examples of the interaction of lines of reflection with rotational symmetries.} \label{fig:Rotation_Examples}
\end{figure}

Examples of how rotational symmetry interacts with reflective symmetry are provided in Figure~\ref{fig:Rotation_Examples}. Notice that the two-fold rotational symmetry in Figure~\ref{fig:Rotation_Examples}(a) has a center off the gridlines---it is this feature that results in lines of reflective symmetry. When the center of two-fold rotational symmetry is on a gridline, the lines of reflection can be avoided, as depicted in Figure~\ref{fig:Rotation_Examples}(b). Four-fold rotational symmetry will always imply lines of reflective symmetry, as depicted in Figure~\ref{fig:Rotation_Examples}(c). In these figures, the lines represent reflections, the diamonds represent centers of two-fold rotational symmetry, and the square represents a center of four-fold rotational symmetry.

\subsection{Rosette Symmetry Groups} \label{subsection:rosette}
The \emph{rosette symmetry groups} are the groups of symmetries that fix a point, and they naturally exclude translations and glide-reflections. There are two kinds of rosette symmetry groups: $cn$ contains the set of rotations that make up $n$-fold rotational symmetry with no reflections, and $dn$ contains the set of rotations that make up $n$-fold rotational symmetry along with $n$ lines of reflection.

From the results of Sen \& Martinez, we can quickly understand which of the rosette symmetry groups are compatible with GHPs. If we allow any parity for $|x|$ and $|y|$, we can only attain the following symmetry types:

\begin{theorem}
There are 5 rosette symmetry groups that are compatible with GHPs. They are $c1$, $d1$, $c2$, $d2$, and $d4$.
\end{theorem}

Notice that it is impossible to have a center of 4-fold rotation with no lines of reflection by Corollary~\ref{cor:4Fold}, so $c4$ is not a compatible symmetry type in a GHP. We naturally exclude any rosette group where $n$ is not equal to 1, 2, or 4, since we are working with symmetries that must realign a square grid.

\section{Incorporating Flip-Symmetries} \label{flip_symm_rosettes}

\begin{figure}[h!tbp]
\centering
\subfloat[Front of design]{%
\resizebox*{4cm}{!}{\includegraphics{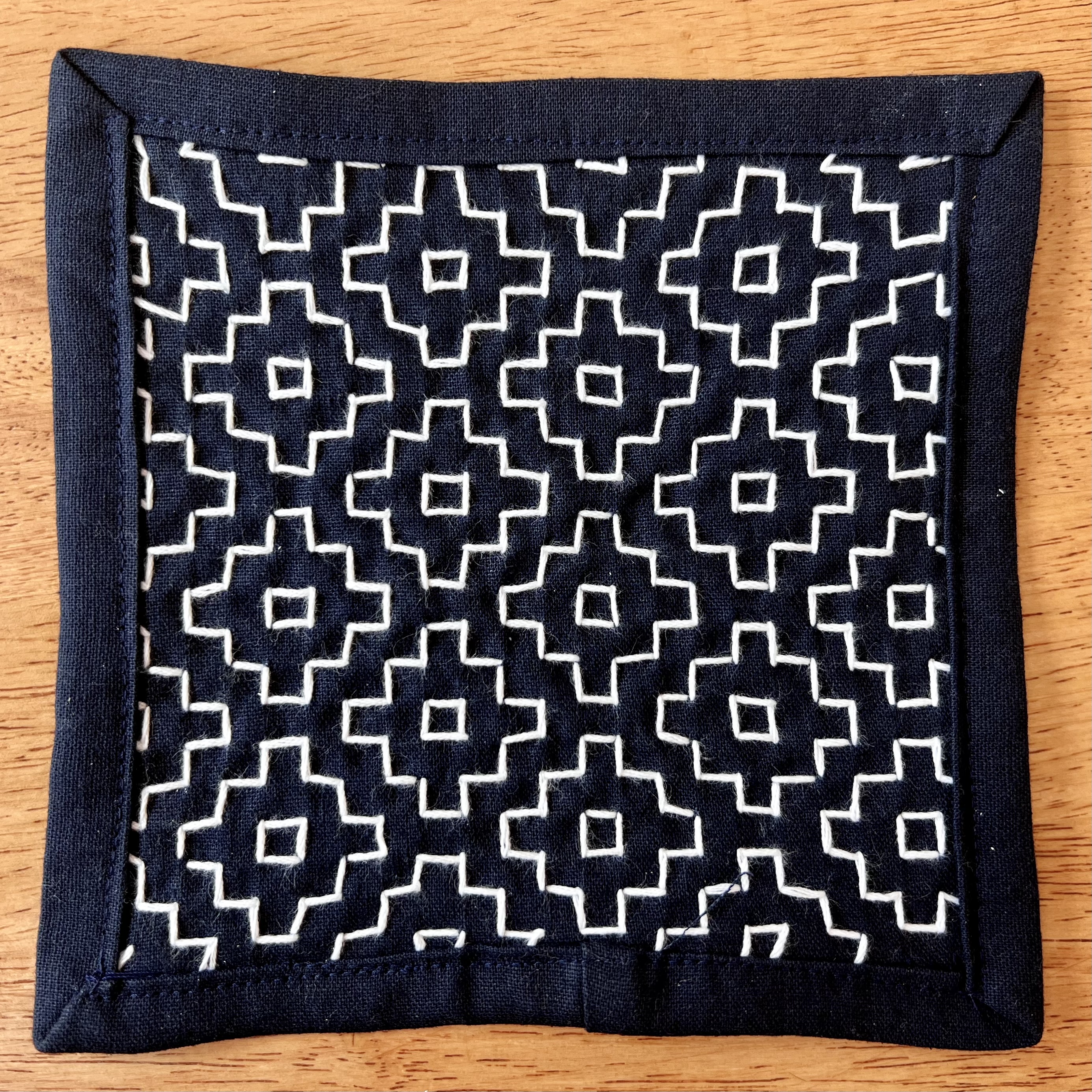}}}\hspace{5pt}
\subfloat[Flipping over]{%
\resizebox*{4cm}{!}{\includegraphics{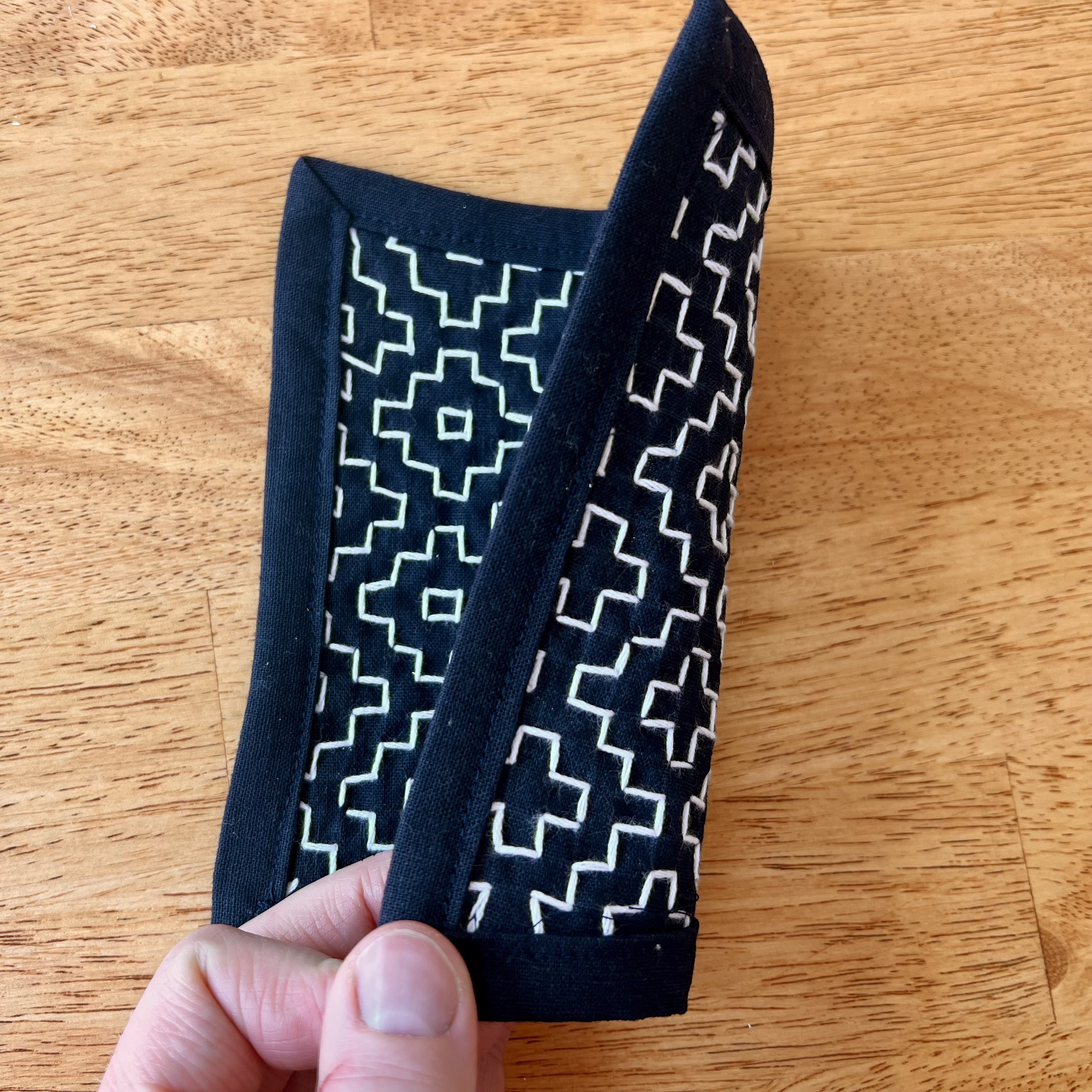}}}\hspace{5pt}
\subfloat[Back of design]{%
\resizebox*{4cm}{!}{\includegraphics{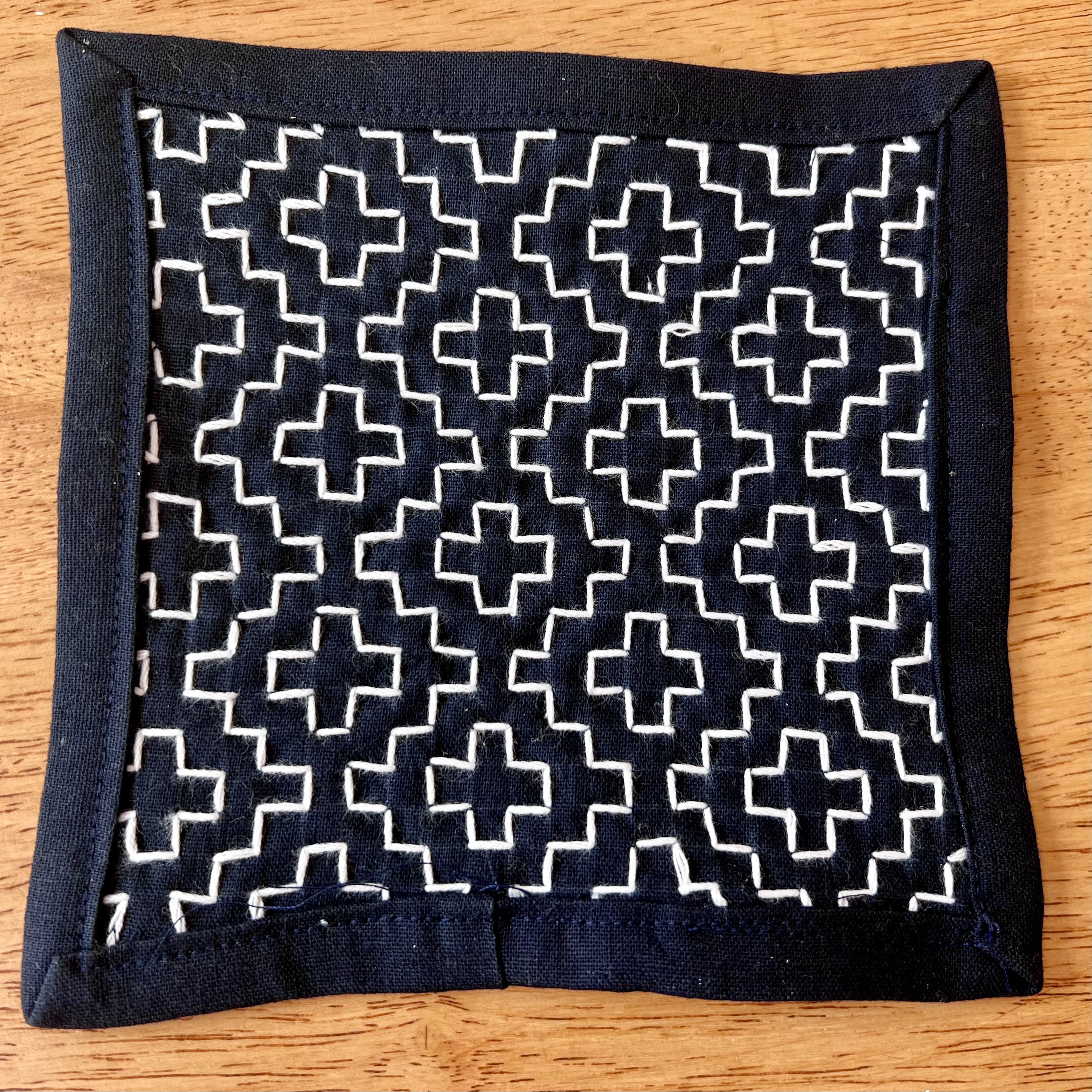}}}
\caption{An example of flipping over a hitomezashi embroidery so that the right edge moves to the left edge} \label{Flip_Example}
\end{figure}

\begin{figure}[h!tbp]
\centering
\subfloat[Front of design]{%
\resizebox*{4cm}{!}{\includegraphics{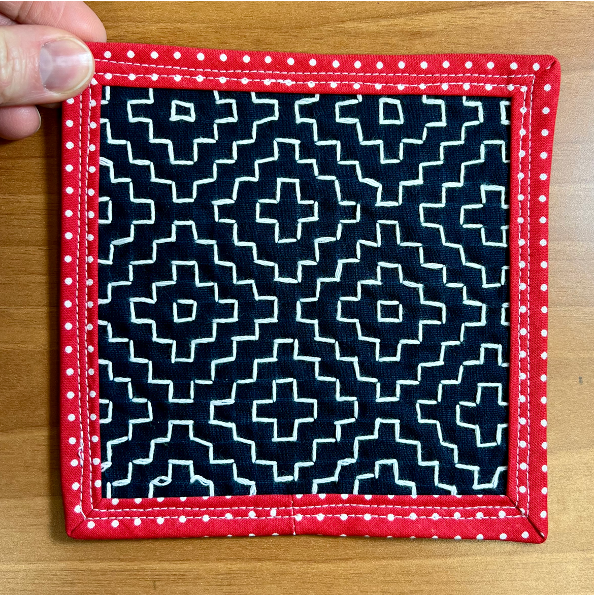}}}\hspace{5pt}
\subfloat[Flipping over]{%
\resizebox*{4cm}{!}{\includegraphics{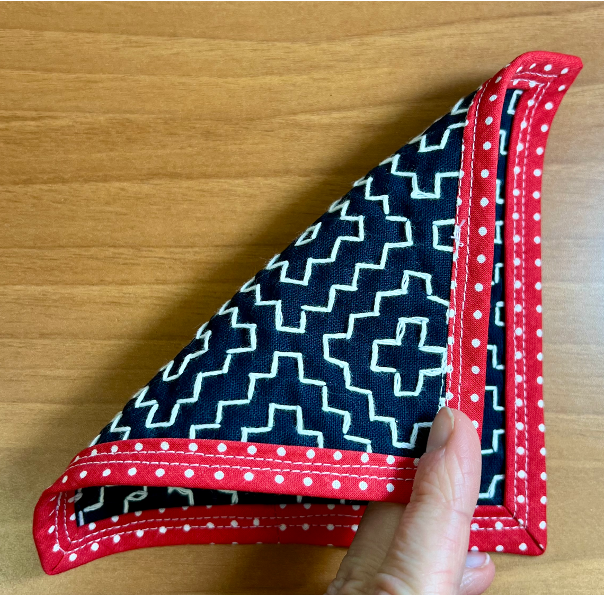}}}\hspace{5pt}
\subfloat[Back of design]{%
\resizebox*{4cm}{!}{\includegraphics{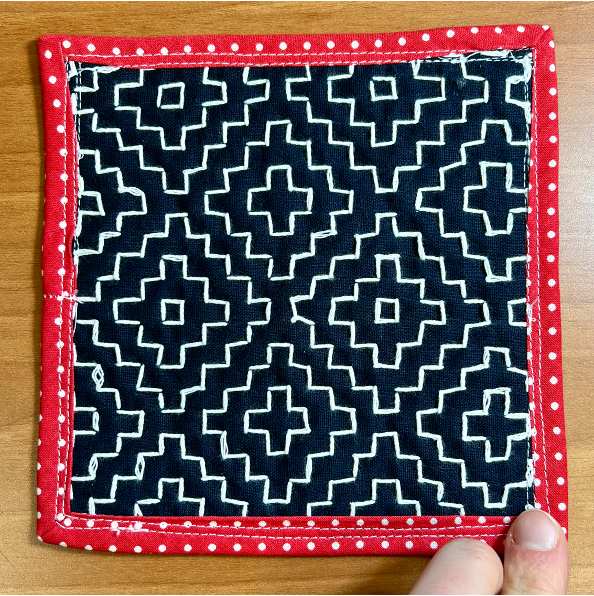}}}
\caption{An example of flipping over a hitomezashi embroidery so that the top left corner moves to the bottom right corner} \label{Flip_Example2}
\end{figure}

In this section we lay out the theory of what we dub `flip-symmetries.' Imagine you would like to turn over your carefully stitched GHP to see the back. How would you do it? You might grasp the right edge and flip the fabric so the right edge is now on the left, as in Figure~\ref{Flip_Example}. Alternatively, you might turn the embroidery over so the top left corner moves to the bottom right corner as in Figure~\ref{Flip_Example2}. These are both examples of a special kind of reflection that converts the front stitches into back stitches---we will call this kind of reflection a \emph{flip-reflection}. The flip-reflection in Figure~\ref{Flip_Example2} is particularly intriguing because the design remains unchanged after the flip-reflection is applied. This flip-reflection is a symmetry of the GHP!

Any isometry can be combined with a `flip' moving the front stitches to the back stitches to produce a \emph{flip-isometry}. In this way we will produce an entire set of isometries corresponding to the different ways to look at the back of our design; there will be eight of these, each corresponding to a symmetry of the square.

Define a function $B$ such that $B(x,y) = (x^C, y^C)$. This mapping lets us look through the cloth so that we can see the stitches on the back, rather than the front. Composed with our four reflections, this describes our flip-reflections (i.e.\ the four different ways to flip over our fabric to see the back). For notation, we will use a bar over the top of an isometry to indicate that this is a flipped version of that isometry. For instance, the flip depicted in Figure \ref{Flip_Example} is $\overline{R_V}=B \circ R_V$ while the flip depicted in Figure \ref{Flip_Example2} is $\overline{R_D}=B \circ R_D$. Note that $B$ and any isometry commute with each other.

With this definition, we can quickly determine how a flip-reflection transforms the binary strings of a GHP in the style of Sen \& Martinez. We simply take the complement of each output from Lemma~\ref{lem_reflections} to obtain the following lemma.

\begin{lemma} \label{lem_flip_reflections}
The reflections $\overline{R_H}$, $\overline{R_V}$, $\overline{R_D}$, and $\overline{R_A}$ have the following effects on the binary strings that generate a GHP:
\begin{center}
$
\overline{R_H}(x, y) =
    \begin{cases}
        (x^C, y^{RC}) & \text{$m$ even}\\
        (x, y^{RC}) & \text{$m$ odd}
    \end{cases}
$
\hspace{1cm}
$
\overline{R_V}(x, y) =
    \begin{cases}
        (x^{RC}, y^C) & \text{$n$ even}\\
        (x^{RC}, y) & \text{$n$ odd}
    \end{cases}
$

$
    \overline{R_D}(x, y) = (y^C, x^C)
$
\hspace{1cm}
$
\overline{R_A}(x, y) =
    \begin{cases}
        (y^{RC}, x^{RC}) & \text{$n$ even}\\
        (y^{R},x^{R}) & \text{$n$ odd}.
    \end{cases}
$
\end{center}
\end{lemma}

From a group-theoretic perspective, the introduction of these flip-reflections, means we should also introduce \emph{flip-rotations}, defined by a composition of the map $B$ with our rotation maps. Complementing the results of Lemma~\ref{lem_rotations} leaves us with the mappings described in the next lemma.

\begin{lemma} \label{lem_flip_rotations}
The counterclockwise rotations $\overline{r_{90}}$, $\overline{r_{180}}$, and $\overline{r_{270}}$ have the following effects on the binary strings that generate a GHP:

\begin{minipage}[h]{0.45\textwidth}

$
    \overline{r_{90}}(x, y) =
    \begin{cases}
        (y^{RC}, x^C) & \text{$m$ even}\\
        (y^{RC}, x) & \text{$m$ odd}
    \end{cases}
$

\vspace{4mm}
$
    \overline{r_{270}}(x, y) =
    \begin{cases}
        (y^C, x^{RC}) & \text{$n$ even}\\
        (y, x^{RC}) & \text{$n$ odd}
    \end{cases}
$
\end{minipage}
\begin{minipage}[h]{0.5\textwidth}
$
    \overline{r_{180}}(x, y) =
    \begin{cases}
        (x^{RC}, y^{RC}) & \text{$n, m$ even}\\
        (x^{R},y^{RC}) & \text{$n$ even, $m$ odd}\\
        (x^{RC},y^{R}) & \text{$n$ odd, $m$ even}\\
        (x^{R},y^{R}) & \text{$n,m$ odd}
    \end{cases}
$
\end{minipage}
\end{lemma}

We now turn our attention to determining when a GHP exhibits a \emph{flip-symmetry}; that is, when a flip-reflection or flip-rotation is a symmetry for the design. First, we note a few basic facts.

\begin{lemma} \label{lem:flip_not}
If $s$ is a symmetry for a GHP generated by $(x,y)$, then $B \circ s = s \circ B$ is not a symmetry of that design.
\end{lemma}

\begin{proof}
This follows immediately from the fact that $a \neq a^c$ for any binary string $a$.
\end{proof}

Note that the above lemma means $B$ alone is never a symmetry for a GHP.

\begin{lemma}\label{not_symm}
The mappings $\overline{R_H}$ and $\overline{R_V}$ will never be symmetries for a GHP.
\end{lemma}

\begin{proof}
In $\overline{R_H}$, we know that $x \neq x^C$ when $m$ is even and $y \neq y^{RC}$ when $m$ is odd by Lemma~\ref{basic_facts}. The argument is similar for $\overline{R_V}$.
\end{proof}

With this in mind, we obtain flip-symmetries in the scenarios described in the following theorem.

\begin{theorem} \label{flip_mappings}
Consider a GHP generated by binary strings $(x,y)$ where $|x|=n$ and $|y|=m$. Then:
\begin{itemize}
\item $\overline{R_D}$ is a symmetry if and only if $n=m$ and $x=y^C$
\item $\overline{R_A}$ is a symmetry if and only if $x=y^{RC}$ when $n=m$ is even and $x=y^R$ when $n=m$ is odd
\item $\overline{r_{90}}$ and $\overline{r_{270}}$ are symmetries if and only if $n=m$ is even and $x=x^R$, $y=y^R$, and $x=y^{RC}$
\item $\overline{r_{180}}$ is a symmetry if and only if one of the following is true:
\begin{itemize}
\item $n$ and $m$ are even, $x=x^{RC}$ and $y=y^{RC}$
\item $n$ and $m$ are odd, $x=x^R$ and $y=y^R$
\end{itemize}
\end{itemize}
\end{theorem}

\begin{proof}
With $\overline{r_{90}}$, we note that $n=m$ odd means $x=y^{RC}=y$ which is impossible. When $n=m$ is even and $\overline{r_{90}}$ is a symmetry, we get $x=y^{RC}$ and $y=x^C$. So $x=(x^C)^{RC}=x^R$ and $y=(y^{RC})^C=y^R$. A similar argument follows for $\overline{r_{270}}$.

The results for $\overline{R_D}, \overline{R_A}, \overline{r_{180}}$ follow immediately from the definitions of these mappings and Lemma~\ref{basic_facts}.
\end{proof}

\subsection{The connection to two-colour symmetry groups}

As we search for the possible symmetry groups containing flip-symmetries, we note that the product of two flipped symmetries gives a standard symmetry, while the product of a standard and flip-symmetry gives a flip-symmetry. Composition between standard and flip-symmetries follow the rules of $D_4$. For instance, $R_A \circ \overline{r_{180}} = \overline{R_D}$ and $\overline{r_{90}} \circ \overline{r_{90}} = r_{180}$.  Hence, the incorporation of the flipped mappings \emph{theoretically} creates a group of transformations that has 16 elements and is isomorphic to $D_4 \times \mathbb{Z}_2$. However, in light of Lemma~\ref{lem:flip_not}, we know that this full group will never be the set of symmetries for a GHP. 

\begin{figure}[h!tbp]
\centering
\subfloat[Design with front stitches in black and back stitches in orange.]{%
\resizebox*{4cm}{!}{\includegraphics{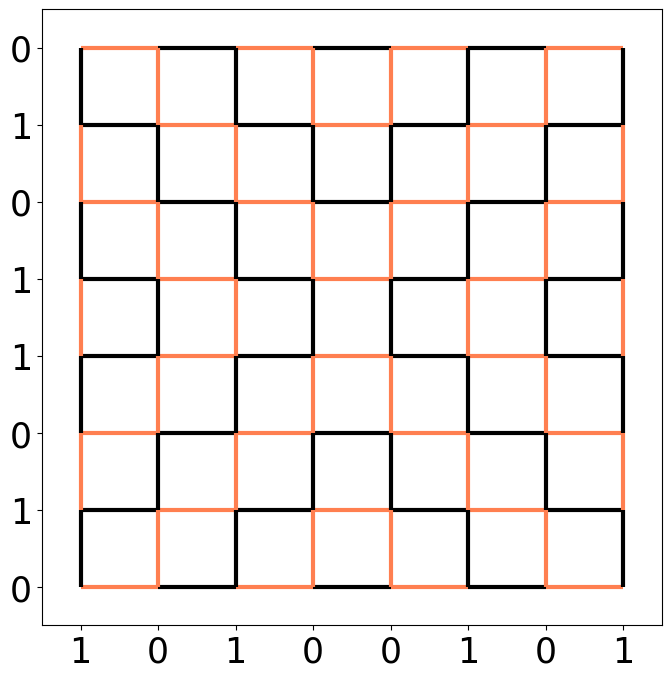}}}\hspace{5pt}
\subfloat[Front stitches only]{%
\resizebox*{4cm}{!}{\includegraphics{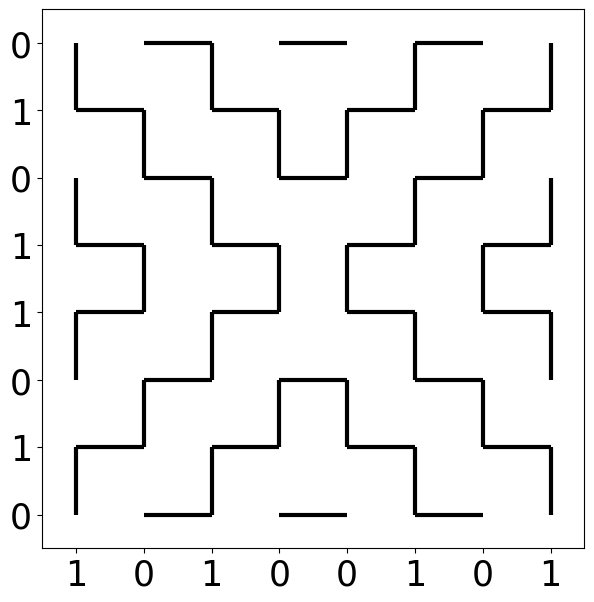}}}\hspace{5pt}
\subfloat[Back stitches only]{%
\resizebox*{4cm}{!}{\includegraphics{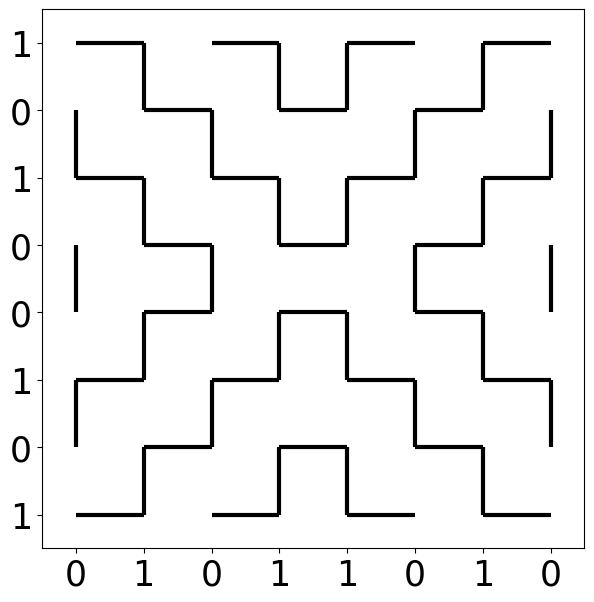}}}

\subfloat[Design with front stitches in black and back stitches in orange with symmetry markings.]{%
\resizebox*{4cm}{!}{\includegraphics{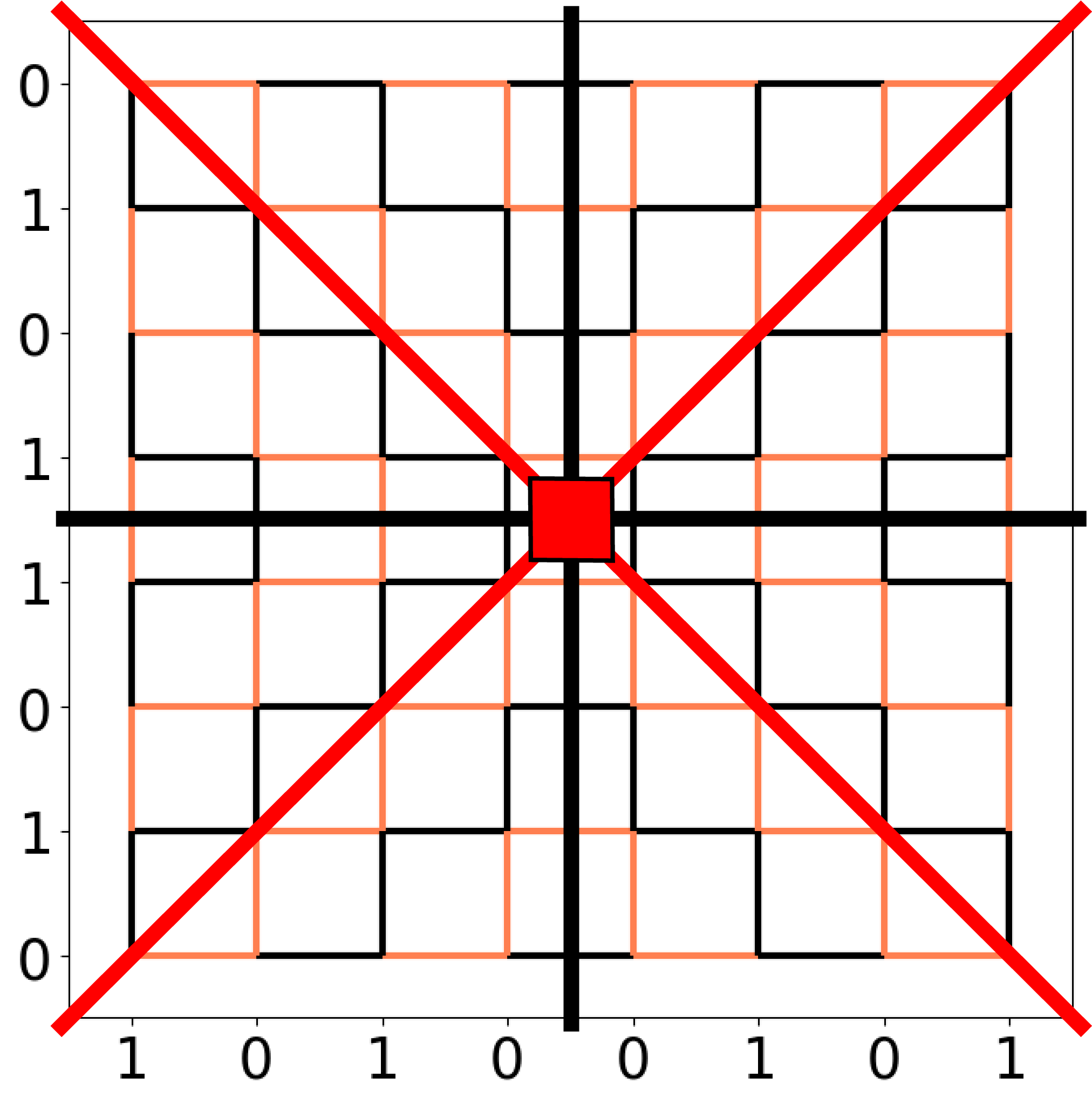}}}\hspace{5pt}
\subfloat[Front stitches only with symmetry markings]{%
\resizebox*{4cm}{!}{\includegraphics{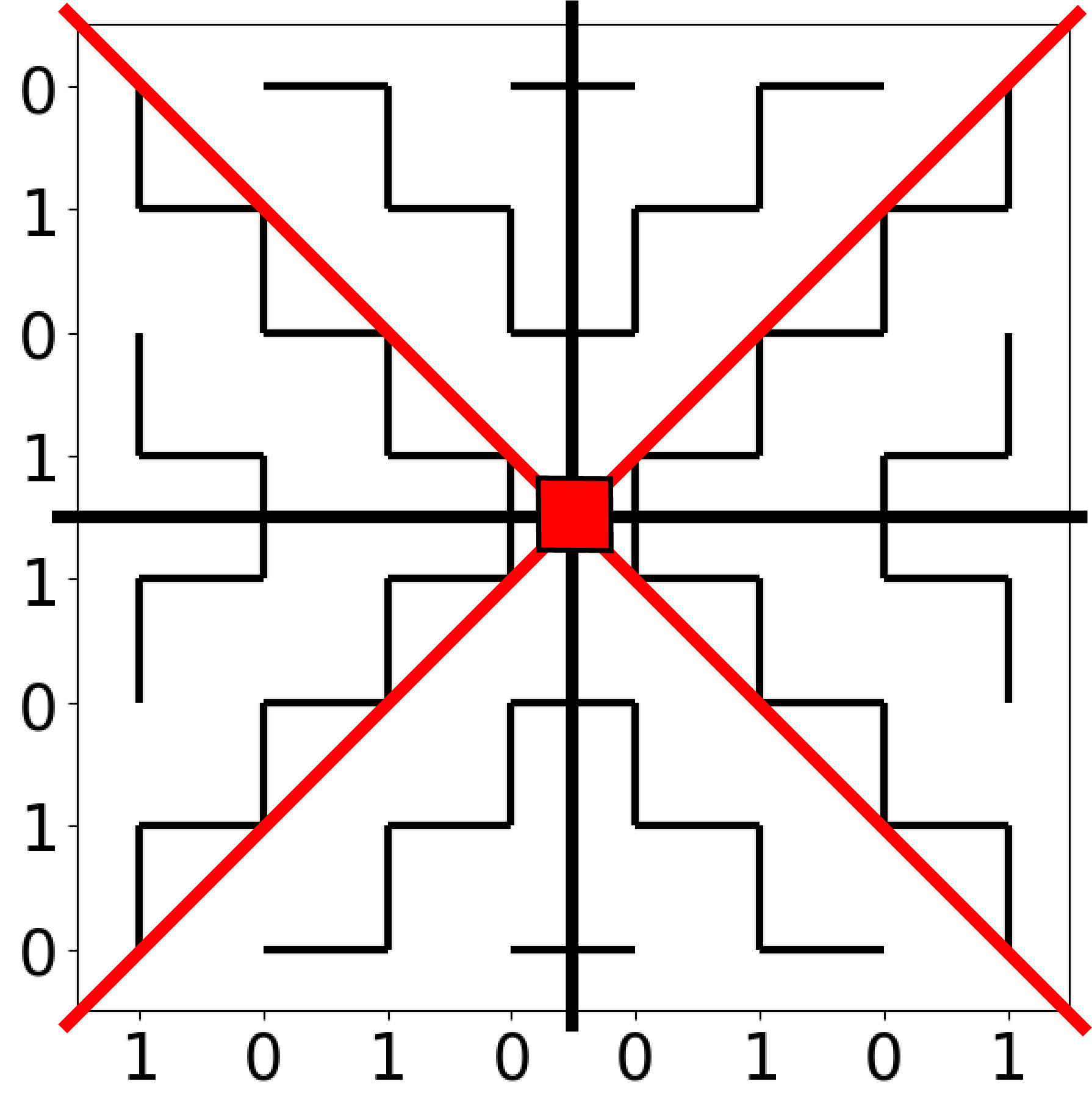}}}\hspace{5pt}
\subfloat[Back stitches only with symmetry markings]{%
\resizebox*{4cm}{!}{\includegraphics{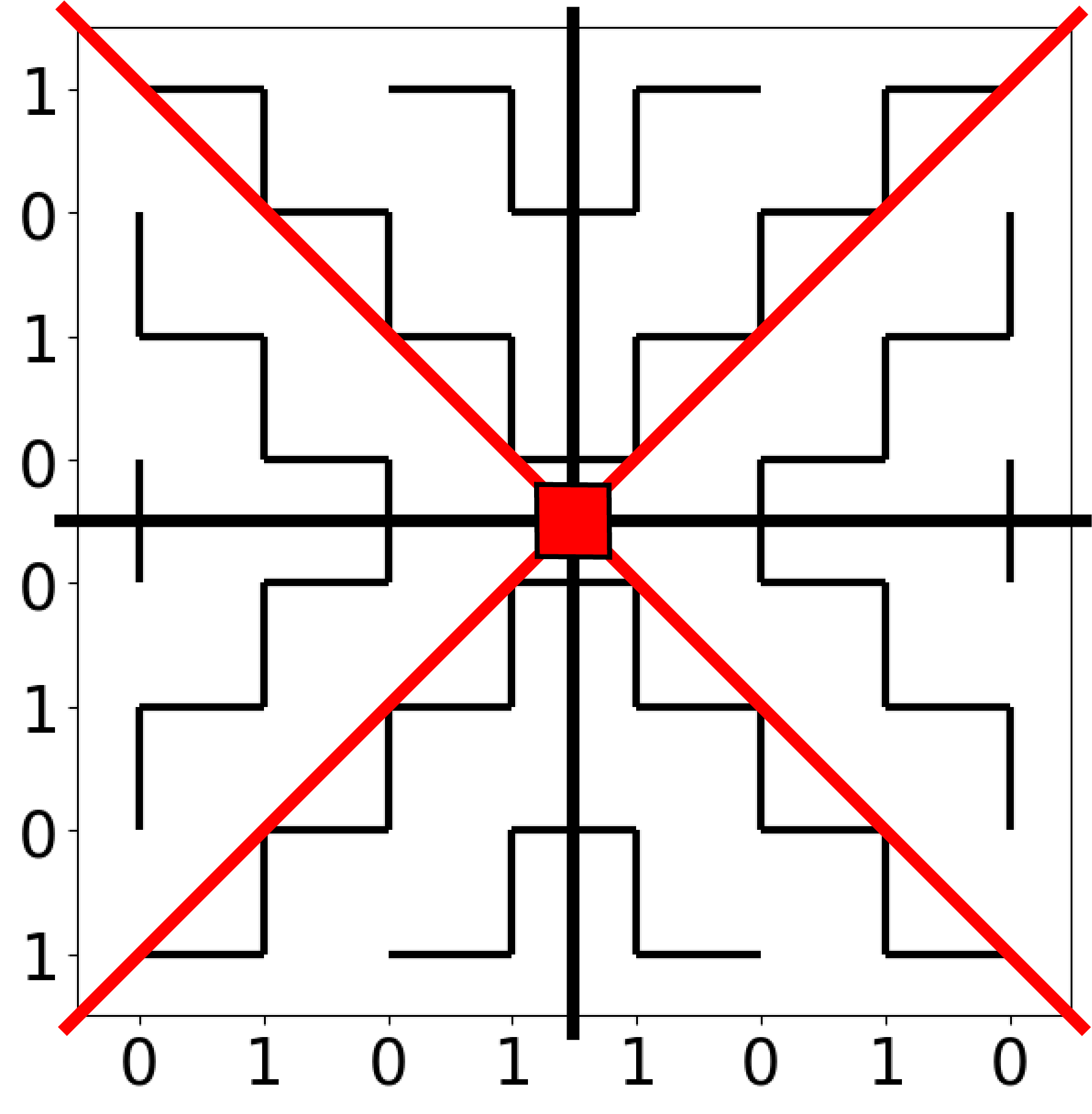}}}
\caption{Examples of how to visualize the front and back stitches of a GHP with symmetries and flip-symmetries.} \label{fig:Two_Color_Ex}
\end{figure}

Describing symmetry groups that incorporate flip-symmetries is equivalent to considering two-colour symmetry groups. Consider Figure~\ref{fig:Two_Color_Ex}(a); this shows the GHP generated by $(10100101,01011010)$ where the stitches on the front are in black. Where there would normally be white space representing the thread passing to the back of the fabric, there are now orange line segments. This is akin to having a translucent fabric where we can see the stitches on the other side. Viewing this two-colour design, a flip-symmetry is equivalent to an isometry in $D_4$ that reverses the black and orange lines.

The two-colour symmetry groups take into account symmetries that are \emph{consistent with colour}; this means that the symmetry either completely preserves or completely reverses the colouring of the pattern. In our language, a symmetry that completely reverses colour is a flip-symmetry. In this way, when we talk about symmetry groups that contain flip-symmetries, we can defer to the theory and notation that has already been developed for two-colour symmetry groups. 

There are two sources that were utilized heavily in this research to understand the theory of two-colour symmetry groups. The first is the beautiful book of Washburn and Crowe, \emph{Symmetries of Culture}, which presents examples and notation for the different two-colour rosette, frieze, and wallpaper groups \citep{Washburn_Crowe}. The second is a book self-published by George Baloglou, \emph{Isometrica}, which is immensely useful for the symmetry diagrams presented for two-colour wallpaper groups \citep{Isometrica}. Any symmetry diagram presented in Section~\ref{wallpaper} was found in this book and verified with the examples provided by Washburn \& Crowe.

From a practical perspective, Figure~\ref{fig:Two_Color_Ex}(a) is very hard to parse visually and can be difficult to use to understand which flip-symmetries are present in a design. An alternative way to present the information in Figure~\ref{fig:Two_Color_Ex}(a) is with two linked images, like those in Figures~\ref{fig:Two_Color_Ex}(b) and \ref{fig:Two_Color_Ex}(c). Figure~\ref{fig:Two_Color_Ex}(b) shows only the front stitches, while Figure~\ref{fig:Two_Color_Ex}(c) shows only the back stitches as they would be seen as if looking through a transparent fabric. Using this method, we can find the colour-preserving symmetries (those symmetries that map front stitches to front stitches) using just one of the figures. To find the colour-reversing or flip-symmetries, we must look for isometries that convert the image of the front stitches into the image of the back stitches. 

The symmetry markings have been provided for both presentations of this design in Figures~\ref{fig:Two_Color_Ex}(d), \ref{fig:Two_Color_Ex}(e), and \ref{fig:Two_Color_Ex}(f). Markings in black denote standard symmetries; markings in red denote flip-symmetries. For centers of rotational symmetry, a black shape means that all the rotational symmetries are not flip-symmetries. When the center of rotation is red, this implies that the smallest rotational symmetry is a flip-symmetry. In the case of a \emph{flip-two-fold rotation} (marked with a red diamond), $\{e,\overline{R_{180}}\}$ would be the corresponding set of symmetries.  In the case of a \emph{flip-four-fold rotation} (marked with a red square), $\{e, \overline{R_{90}}, R_{180}, \overline{R_{270}} \}$ is the corresponding set of symmetries. The legend for our symmetry markings are given in Figure~\ref{fig:legend}; in preparation for analyzing the symmetries of wallpaper patterns, this legend also includes markings for glide-reflection and translation symmetries.

\begin{figure}[h!tbp]
\centering
\begin{minipage}{0.4\textwidth}
\begin{tikzpicture}[scale = 1]
\draw[ultra thick](0,0)--(2,0)node[right]{= Reflection};
\end{tikzpicture}

\begin{tikzpicture}[scale = 1]
\draw[ultra thick, red](0,0)--(2,0)node[right]{\color{black}= Flip-reflection};
\end{tikzpicture}

\begin{tikzpicture}[scale = 1]
\draw[ultra thick, dashed](0,0)--(2,0)node[right]{= Glide-reflection};
\end{tikzpicture}

\begin{tikzpicture}[scale = 1]
\draw[ultra thick, dashed, red](0,0)--(2,0)node[right]{\color{black}= Flip-glide-reflection};
\end{tikzpicture}

\begin{tikzpicture}[scale = 1]
\draw[ultra thick, ->] (0,0)--(2,0)node[right]{= Translation};
\end{tikzpicture}

\begin{tikzpicture}[scale = 1]
\draw[ultra thick, ->,red] (0,0)--(2,0)node[right]{\color{black}= Flip-translation};
\end{tikzpicture}

\end{minipage} \hspace{13mm}
\begin{minipage}{0.4\textwidth}
\begin{tikzpicture}[scale = 1]
\draw(0,0)node[2fold]{};
\draw(0.2,0)node[right]{= Two-fold rotation};
\end{tikzpicture}

\begin{tikzpicture}[scale = 1]
\draw(0,0)node[2foldcolorswap]{};
\draw(0.2,0)node[right]{= Flip-two-fold rotation};
\end{tikzpicture}

\begin{tikzpicture}[scale = 1]
\draw(0,0)node[4fold]{};
\draw(0.18,0)node[right]{= Four-fold rotation};
\end{tikzpicture}

\begin{tikzpicture}[scale = 1]
\draw(0,0)node[4foldcolorswap]{};
\draw(0.18,0)node[right]{= Flip-four-fold rotation};
\end{tikzpicture}
\end{minipage}

\caption{The legend for symmetry notation that will be used throughout this paper.}
\label{fig:legend}
\end{figure}

\subsection{Two-colour rosette symmetries}

With the results of Theorem~\ref{flip_mappings} we can determine all the possible two-colour rosette symmetry groups that are compatible with GHPs. We first provide descriptions and examples of the two-colour rosette symmetry groups.

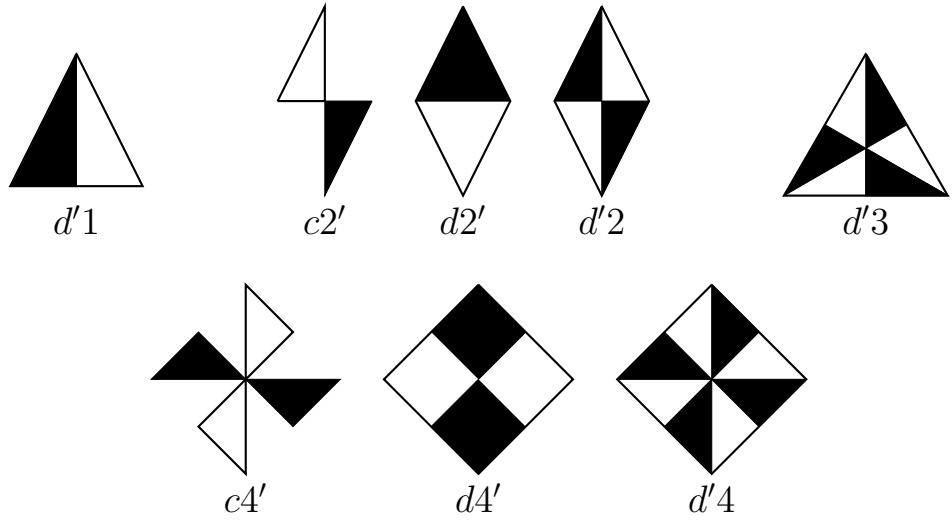
\begin{figure}[h!tbp]
\centering

\begin{tikzpicture}[scale=1.75]
\draw[thick] (0,1)--(-0.5,0)--(0.5,0)--(0,1);
\draw[fill=black] (-0.5,0)--(0,0)--(0,1)--(-0.5,0);
\draw (0,-.25)node{\Large $d'1$};
\end{tikzpicture}
\hspace{15mm}
\begin{tikzpicture}[scale=1.25]
\draw[thick] (-0.5,0)--(0,0)--(0,1)--(-0.5,0);
\draw[fill=black] (0.5,0)--(0,0)--(0,-1)--(0.5,0);
\draw (0,-1.25)node{\Large $c2'$};
\end{tikzpicture}
\hspace{3mm}
\begin{tikzpicture}[scale=1.25]
\draw[thick] (0,1)--(-0.5,0)--(0,-1)--(0.5,0)--(0,1);
\draw[fill=black] (-0.5,0)--(0.5,0)--(0,1)--(-0.5,0);
\draw (0,-1.25)node{\Large $d2'$};
\end{tikzpicture}
\hspace{3mm}
\begin{tikzpicture}[scale=1.25]
\draw[thick] (0,1)--(-0.5,0)--(0,-1)--(0.5,0)--(0,1);
\draw[fill=black] (-0.5,0)--(0,0)--(0,1)--(-0.5,0);
\draw[fill=black] (0.5,0)--(0,0)--(0,-1)--(0.5,0);
\draw (0,-1.25)node{\Large $d'2$};
\end{tikzpicture}
\hspace{15mm}
\begin{tikzpicture}[scale=1.25]
\draw[thick] (0,1)--(0.866,-0.5)--(-0.866,-0.5)--(0,1);
\draw[fill=black] (0,0)--(0,-.5)--(0.866,-0.5)--(0,0);
\draw[fill=black] (0,0)--(0,-.5)--(0.866,-0.5)--(0,0);
\draw[fill=black] (0,0)--(-.433,0.25)--(-0.866,-0.5)--(0,0);
\draw[fill=black] (0,0)--(.433,0.25)--(0,1)--(0,0);
\draw (0,-.75)node{\Large $d'3$};
\end{tikzpicture}

\vspace{5mm}
\begin{tikzpicture}[scale=1.25]
\draw[fill=black] (0,0)--(.5,-.5)--(1,0)--(0,0);
\draw[thick] (0,0)--(-.5,-.5)--(0,-1)--(0,0);
\draw[fill=black] (0,0)--(-.5,.5)--(-1,0)--(0,0);
\draw[thick] (0,0)--(.5,.5)--(0,1)--(0,0);
\draw (0,-1.25)node{\Large $c4'$};
\end{tikzpicture}
\hspace{3mm}
\begin{tikzpicture}[scale=1.25]
\draw[thick] (0,1)--(1,0)--(0,-1)--(-1,0)--(0,1);
\draw[fill=black] (0,0)--(-.5,-.5)--(0,-1)--(.5,-.5);
\draw[fill=black] (0,0)--(.5,.5)--(0,1)--(-.5,.5);
\draw (0,-1.25)node{\Large $d4'$};
\end{tikzpicture}
\hspace{3mm}
\begin{tikzpicture}[scale=1.25]
\draw[thick] (0,1)--(1,0)--(0,-1)--(-1,0)--(0,1);
\draw[fill=black] (0,0)--(.5,-.5)--(1,0)--(0,0);
\draw[fill=black] (0,0)--(-.5,-.5)--(0,-1)--(0,0);
\draw[fill=black] (0,0)--(-.5,.5)--(-1,0)--(0,0);
\draw[fill=black] (0,0)--(.5,.5)--(0,1)--(0,0);
\draw (0,-1.25)node{\Large $d'4$};
\end{tikzpicture}

\caption{Examples of two-colour rosette symmetries} \label{fig:Rosettes}
\end{figure}

Using the notation described by \cite{Washburn_Crowe}, there are three kinds of two-colour rosette symmetry groups: $cn'$, $dn'$, and $d'n$. Examples of the possible two-colour rosette symmetries are given in Figure~\ref{fig:Rosettes} where $n=1,2,3,4$. Both $cn'$ and $dn'$ are symmetry groups that contain centers of rotational symmetry with $n$ rotations where the smallest rotation reverses colors; this means that the different rotations alternately preserve and reverse colours. This alternation means these groups only exist when $n$ is even. While $cn'$ contains no reflections, $dn'$ has $n$ lines of reflection where half preserve colours and half reverse colours. Finally, the symmetry group $d'n$ consists of a center of $n$-fold rotational symmetry (that preserves colour) with $n$ lines of reflection that all reverse colors; this type exists for any positive integer $n$.

\subsection{Describing all possible rosette symmetries}

We can synthesize the results of Theorems~\ref{ros_symm} and \ref{flip_mappings} to not only tell us which rosette groups (one-colour and two-colour) are compatible with GHPs, but also give explicit instructions for how to produce particular symmetries.

Table~\ref{rosettes} provides instructions for producing all the possible rosette symmetry types. While construction of this table is straightforward, it does require some careful consideration. We start by determining the parity of $n=|x|$ and $m=|y|$ and whether these dimensions are equal. From there, we can give requirements on $x$ where $x=x^R$, $x=x^{RC}$, or $x$ has no restriction (and similarly for $y$). We can further dictate how $x$ and $y$ relate to each other: $x=y$, $x=y^R$, $x=y^C$, $x=y^{RC}$ or no relation. However, there are not $3\cdot3\cdot5=45$ different possibilities to consider: if $x=x^R$ and $y=y^{RC}$, for instance, it is impossible to have $x=y$ since this implies $x=x^R=x^{RC}=x^C$. In this situation, we also couldn't have $x=y^R$, as this would imply $x=x^R=y^R$, so $x=y$ (which again implies $x=x^C$). Through careful process of elimination, we significantly reduce the scenarios we need to consider to 24 entries.

\def\arraystretch{1.3}
\begin{table}[h!tbp]
\resizebox{1\textwidth}{!}{
\begin{tabular}{|c|c|c|c|c|c|}
\hline
Properties of $n$ and $m$             & Properties of $x$           & Properties of $y$           & Relationship of $x$ and $y$ & Symmetries                                                                             & Type \\ \hline

None                                  & None                        & None                        & None                        & $e$                                                                                    & $c1$            \\ \hline
\multirow{4}{*}{$n$ and $m$ even}     & None                        & $y=y^R$                     & None                        & $e,R_H$                                                                                & $d1$            \\ \cline{2-6} 
                                      & \multirow{2}{*}{$x=x^R$}    & None                        & None                        & $e,R_V$                                                                                & $d1$            \\ \cline{3-6} 
                                      &                             & $y=y^R$                     & None                        & $e,r_{180},R_H,R_V$                                                                    & $d2$            \\ \cline{2-6} 
                                      & $x=x^{RC}$                  & $y=y^{RC}$                  & None                        & $e,\overline{r_{180}}$                                                                 & $c2'$         \\ \hline
\multirow{8}{*}{$n=m$ even}           & \multirow{4}{*}{None}       & \multirow{4}{*}{None}       & $x=y$                       & $e,R_D$                                                                                & $d1$            \\ \cline{4-6} 
                                      &                             &                             & $x=y^R$                     & $e,R_A$                                                                                & $d1$            \\ \cline{4-6} 
                                      &                             &                             & $x=y^C$                     & $e,\overline{R_D}$                                                                     & $d'1$         \\ \cline{4-6} 
                                      &                             &                             & $x=y^{RC}$                  & $e,\overline{R_A}$                                                                     & $d'1$         \\ \cline{2-6} 
                                      & \multirow{2}{*}{$x=x^R$}    & \multirow{2}{*}{$y=y^R$}    & $x=y=y^R$                   & $e,r_{90},r_{180},r_{270},R_H,R_V,R_D,R_A$                                             & $d4$            \\ \cline{4-6} 
                                      &                             &                             & $x=y^C=y^{RC}$              & $e,\overline{r_{90}},r_{180},\overline{r_{270}},R_H,R_V,\overline{R_D},\overline{R_A}$ & $d4'$         \\ \cline{2-6} 
                                      & \multirow{2}{*}{$x=x^{RC}$} & \multirow{2}{*}{$y=y^{RC}$} & $x=y=y^{RC}$                & $e,\overline{r_{180}},R_D,\overline{R_A}$                                              & $d2'$         \\ \cline{4-6} 
                                      &                             &                             & $x=y^R=y^C$                 & $e,\overline{r_{180}},\overline{R_D},R_A$                                              & $d2'$         \\ \hline
\multirow{2}{*}{$n$ even and $m$ odd} & $x=x^R$                     & None                        & None                        & $e,R_V$                                                                                & $d1$            \\ \cline{2-6} 
                                      & $x=x^{RC}$                  & $y=y^R$                     & None                        & $e,r_{180}$                                                                            & $c2$            \\ \hline
\multirow{2}{*}{$n$ odd and $m$ even} & None                        & $y=y^R$                     & None                        & $e,R_H$                                                                                & $d1$            \\ \cline{2-6} 
                                      & $x=x^R$                     & $y=y^{RC}$                  & None                        & $e,r_{180}$                                                                            & $c2$            \\ \hline
$n$ and $m$ odd                       & $x=x^R$                     & $y=y^R$                     & None                        & $e,\overline{r_{180}}$                                                                 & $c2'$         \\ \hline
\multirow{6}{*}{$n=m$ odd}            & \multirow{4}{*}{None}       & \multirow{4}{*}{None}       & $x=y$                       & $e,R_D$                                                                                & $d1$            \\ \cline{4-6} 
                                      &                             &                             & $x=y^R$                     & $e,\overline{R_A}$                                                                     & $d'1$         \\ \cline{4-6} 
                                      &                             &                             & $x=y^C$                     & $e,\overline{R_D}$                                                                     & $d'1$         \\ \cline{4-6} 
                                      &                             &                             & $x=y^{RC}$                  & $e,R_A$                                                                                & $d1$            \\ \cline{2-6} 
                                      & \multirow{2}{*}{$x=x^R$}    & \multirow{2}{*}{$y=y^R$}    & $x=y=y^R$                   & $e,\overline{r_{180}},R_D,\overline{R_A}$                                              & $d2'$         \\ \cline{4-6} 
                                      &                             &                             & $x=y^C=y^{RC}$              & $e,\overline{r_{180}},\overline{R_D},R_A$                                              & $d2'$         \\ \hline
\end{tabular}
}

\caption{A summary of all the ways to obtain symmetries in a GHP.}
\label{rosettes}
\end{table}

Table~\ref{rosettes} is the result of the careful collation of all these considerations. The first column gives the restrictions on $n$ and $m$. The next three columns give the properties of $x$, $y$, and the relationship between $x$ and $y$. The fifth column lists the symmetries present in a GHP that has the properties of the first four columns. The final column gives the rosette group (either one-colour or two-colour) to which the set of symmetries is associated.

From this table, we arrive at the following theorem:

\begin{theorem}
There are 4 two-colour rosette symmetry groups that are compatible with GHPs. They are:

\begin{itemize}
\item $d'1$: contains one flip-reflection
\item $c2'$: contains a center of flip-two-fold rotation symmetry
\item $d2'$: contains one line of reflection and one perpendicular line of flip-reflection that intersect at a center of flip-two-fold rotation symmetry
\item $d4'$: contains two lines of reflection and two lines of flip-reflection that form $45^{\circ}$ angles, intersecting at a center of flip-four-fold rotation symmetry
\end{itemize}

\end{theorem}

\begin{figure}[h!tbp]
\centering
\subfloat[$d'1$: front stitches and symmetry markings]{%
\resizebox*{4cm}{!}{\includegraphics{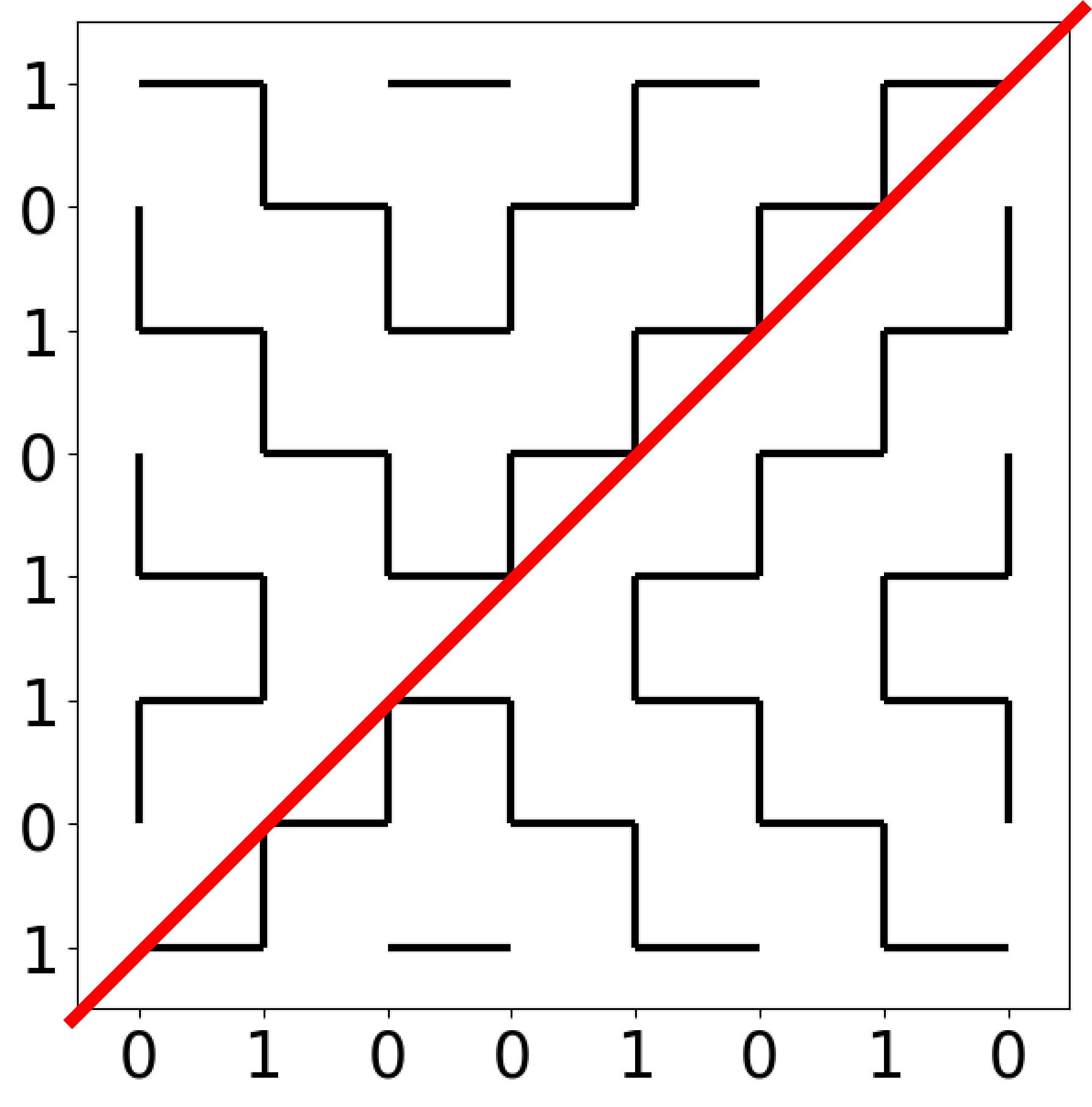}}}\hspace{5pt}
\subfloat[$d'1$: front stitches]{%
\resizebox*{4cm}{!}{\includegraphics{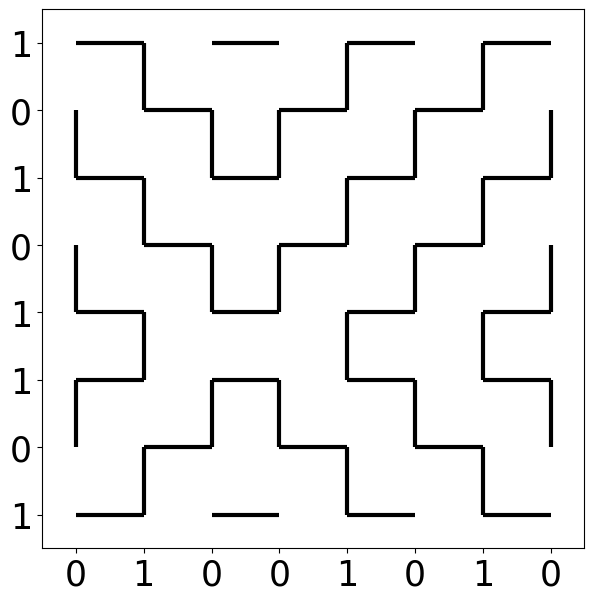}}}\hspace{5pt}
\subfloat[$d'1$: back stitches]{%
\resizebox*{4cm}{!}{\includegraphics{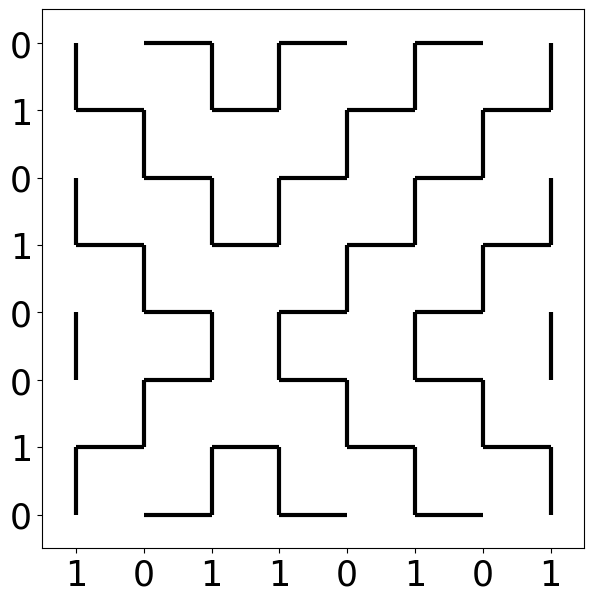}}}

\subfloat[$c2'$: front stitches and symmetry markings]{%
\resizebox*{4cm}{!}{\includegraphics{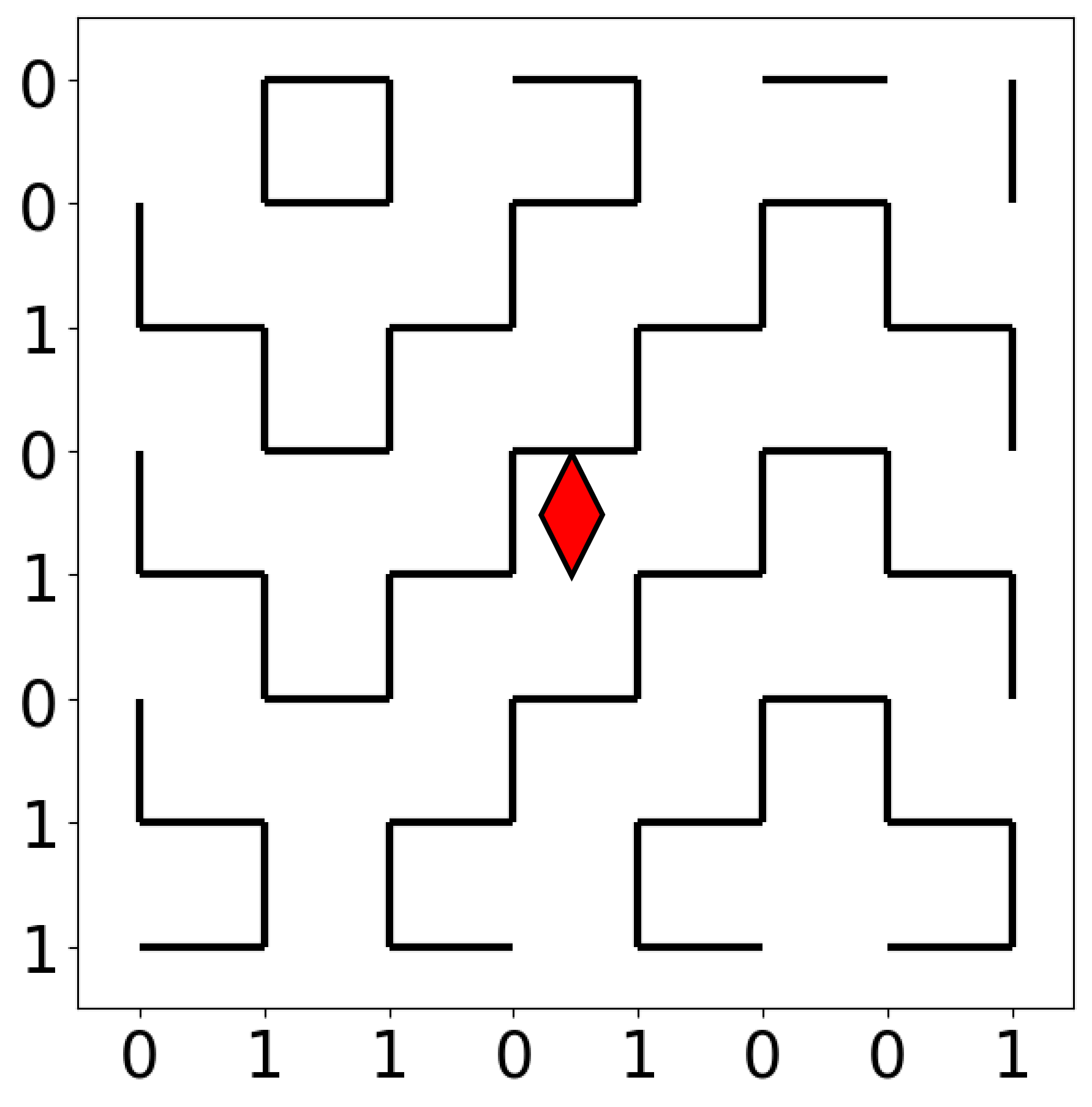}}}\hspace{5pt}
\subfloat[$c2'$: front stitches]{%
\resizebox*{4cm}{!}{\includegraphics{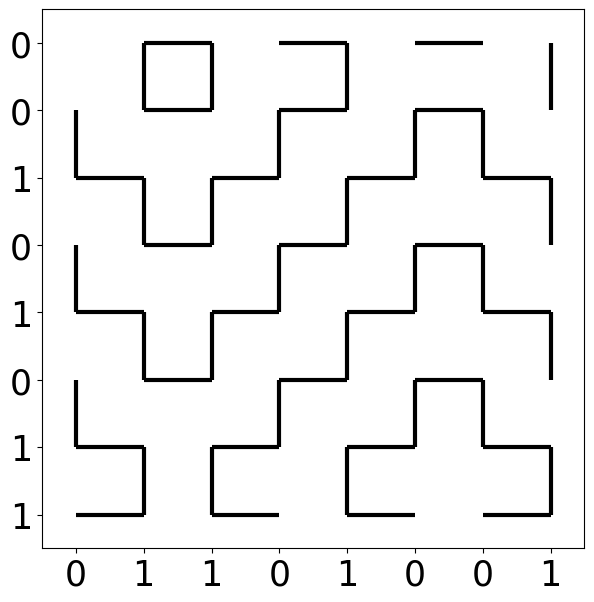}}}\hspace{5pt}
\subfloat[$c2'$: back stitches]{%
\resizebox*{4cm}{!}{\includegraphics{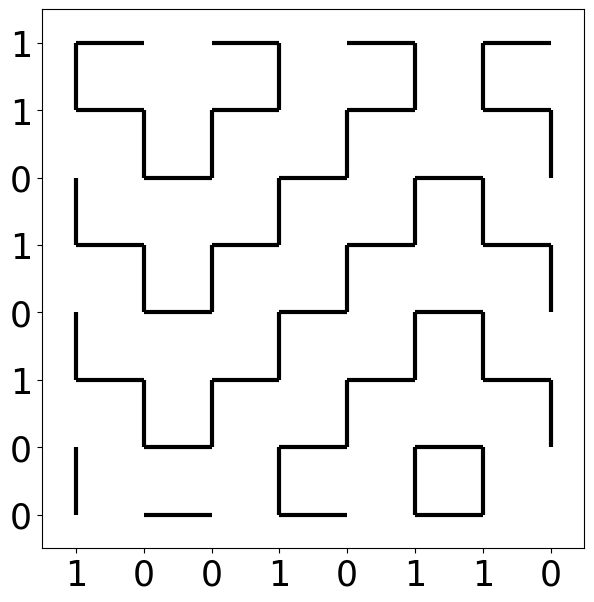}}}

\subfloat[$d2'$: front stitches and symmetry markings]{%
\resizebox*{4cm}{!}{\includegraphics{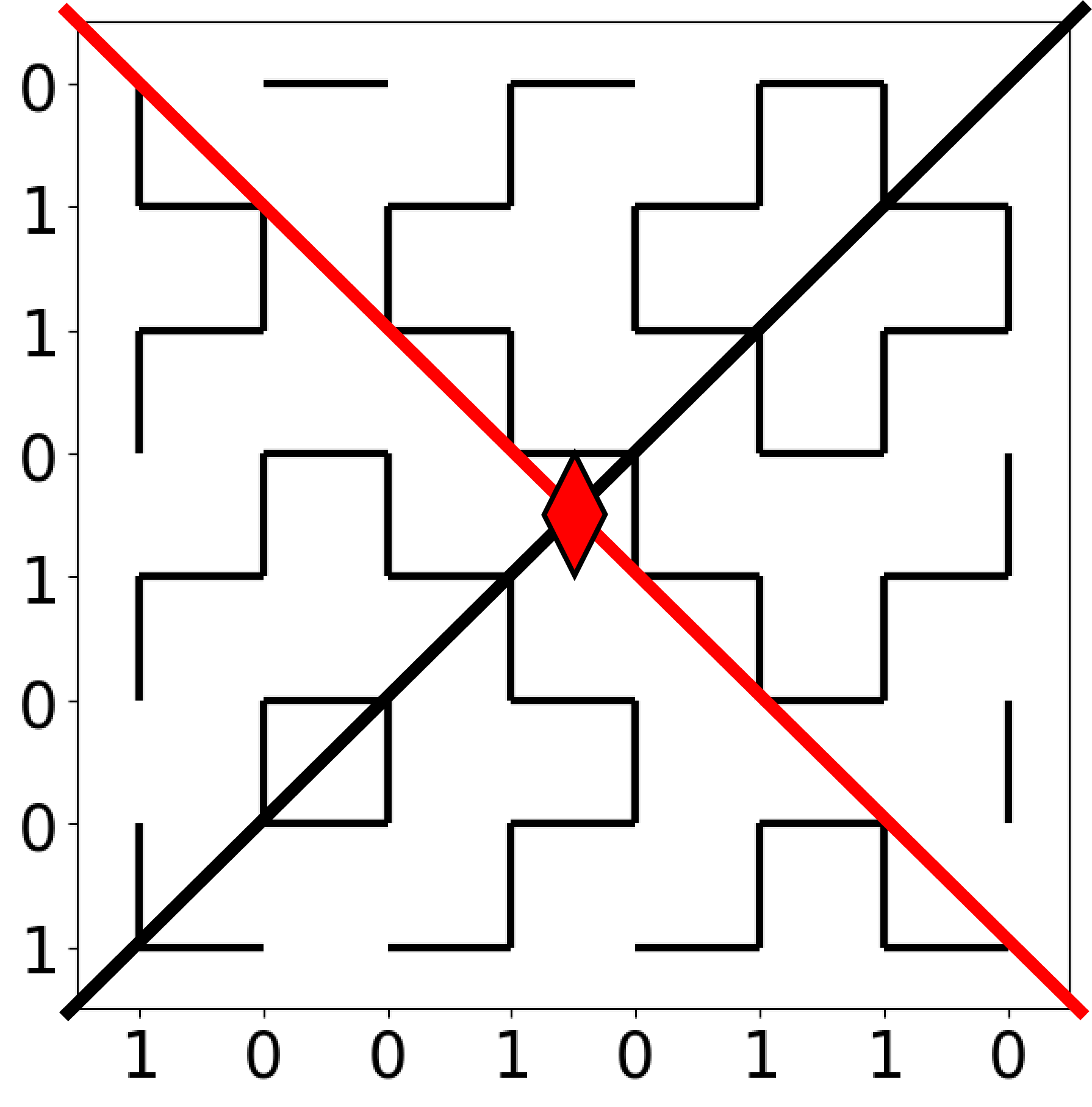}}}\hspace{5pt}
\subfloat[$d2'$: front stitches]{%
\resizebox*{4cm}{!}{\includegraphics{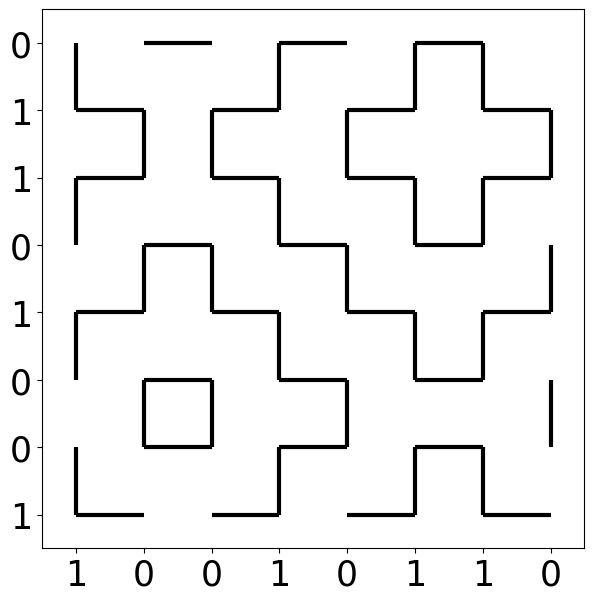}}}\hspace{5pt}
\subfloat[$d2'$: back stitches]{%
\resizebox*{4cm}{!}{\includegraphics{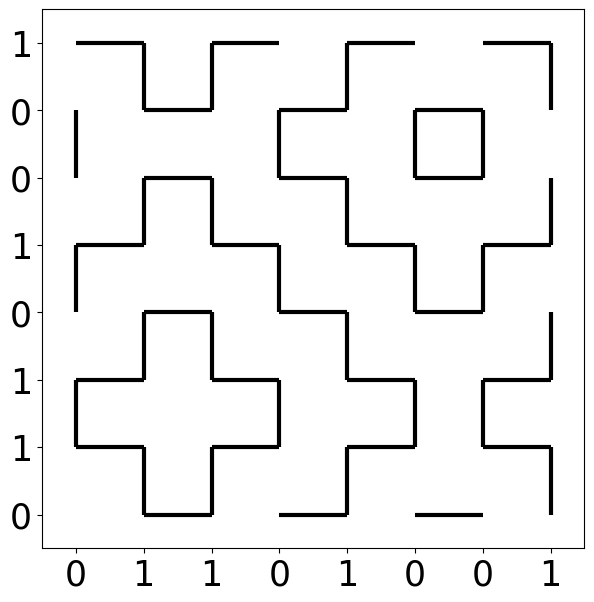}}}

\caption{Examples of GHPs with d'1, c2', and d2' symmetry.} \label{fig:Two_Color_Rosettes}
\end{figure}

Examples of $d'1$, $c2'$, and $d2'$ are provided in Figure \ref{fig:Two_Color_Rosettes}. An example of $d4'$ is pictured in Figure~\ref{fig:Two_Color_Ex}.

Of the possible two-colour rosette groups with $n=1, 2, 4$, we observe that Table~\ref{rosettes} does not include $d'2$, $c4'$, and $d'4$; these types are incompatible with GHPs. This fact will be utilized as we begin to analyse two-colour wallpaper symmetries.

%
%
%

\section{Wallpaper Patterns} \label{wallpaper}

Now we turn our attention to wallpaper patterns, which are unbounded in all directions and characterized by translation symmetry in two distinct directions. We can generate these patterns with two periodic binary words. For instance, we could take $x=\ldots 0011000110\ldots=(00110)^{\infty}$ and $y=\ldots 11011101\ldots=(1101)^{\infty}$, and assign each entry of these periodic words to an integer point on the $x$- or $y$-axis. Stitches can then be populated in all directions following the standard practice of alternating above and below the fabric. Sen \& Martinez called these \emph{periodic generalized hitomezashi patterns} (PGHPs for short).

There is some ambiguity with this setup as shifting $x$ one unit to the left or right while preserving the position of $y$ will create a different design. However, the focus of this section is not to describe the precise conditions on $x$ and $y$ that will produce symmetries; rather, we will focus on describing which wallpaper symmetry groups are compatible with PGHPs. Therefore, this ambiguity does not affect our results.

There are 17 different one-colour wallpaper groups. In this paper, we will use the shortened form of the crystallographic notation; the full crystallographic notation has 4 characters, the shortened form uses fewer. For those looking for a primer to understand the wallpaper symmetry groups, using a flow chart to assist with classification is especially helpful. Washburn \& Crowe have such a flow chart in Section 5.1 of \emph{Symmetries of Culture} \citep{Washburn_Crowe}. They are also easy to find across the internet.

Sen \& Martinez showed that 9 of the possible 17 wallpaper groups are compatible with PGHPs. When considering only wallpapers that have 1-, 2-, or 4-fold rotational symmetry, they found that $pgg$, $p4$, and $p4g$ are incompatible with PGHPs. 

\subsection{The two-colour wallpaper groups}

In total, there are 46 two-colour wallpaper groups. Notating all these groups is somewhat challenging and there are generally two ways the notation is presented. The method we prefer is given in the format $T/T'$. Both $T$ and $T'$ are one of the 17 one-colour wallpaper groups where: $T'$ describes the wallpaper group associated with the symmetries that preserve colour and $T$ describes the complete set of symmetries that are consistent with colour.

There is a pair of distinct two-colour wallpaper groups that can both be labelled $pm/pm$. One such group has colour-reversing lines of reflection while the other does not; to separate these two types, they are notated $pm/pm(m')$ and $pm/pm(m)$ (respectively). Otherwise, the $T/T'$ notation provides no ambiguity for two-colour wallpapers and has the advantage of relying only on the 17 one-colour wallpapers to be able to label all 46 of the two-colour wallpaper groups.

\begin{figure}[h!tbp]
\centering
\subfloat[$p_b'gm$ or $pmm/pmg$]{%
\resizebox*{6.5cm}{!}{
\begin{tikzpicture}
\foreach \j in {0,3,6}{
    \foreach \i in {0,2,4,6,8}{
    	\draw[fill=black] (\i,\j)--(\i+.5,\j+.5)--(\i+1,\j)--(\i+1,\j+1.5)--(\i,\j+1.5);
    	\draw[fill=black] (\i,\j)--(\i-.5,\j-.5)--(\i-1,\j)--(\i-1,\j-1.5)--(\i,\j-1.5);
    	\draw[fill=black] (\i,\j)--(\i+.5,\j-.5)--(\i+1,\j);
    	\draw[fill=black] (\i,\j)--(\i-.5,\j+.5)--(\i-1,\j);
    	}
	}
\end{tikzpicture}
}} \hspace{1cm}
\subfloat[$p_b'mg$ or $pmg/pmg$]{%
\resizebox*{6.5cm}{!}{
\begin{tikzpicture}
\foreach \i in {0,4,8}{
	\foreach \j in {0,2,4,6,8}{
		\draw[fill=black] (\i,\j)--(\i+2,\j+1)--(\i+2,\j)--(\i,\j-1);
		\draw[fill=black] (\i+2,\j+1)--(\i+4,\j)--(\i+4,\j-1)--(\i+2,\j);
		}
	}
\end{tikzpicture}
}}

\caption{Examples of two-colour wallpapers}
\label{fig:2color_wallpapers}
\end{figure}
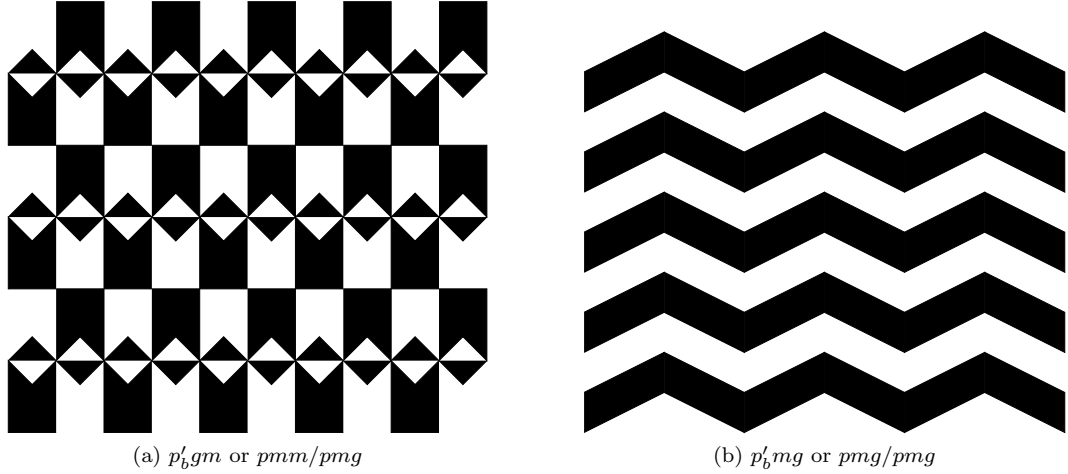

There is an additional method of notation described by Washburn \& Crowe using a string of characters; this is the method of notation used in Baloglou's \emph{Isometrica}. We provide an example of that notation in Figure~\ref{fig:2color_wallpapers}, but in general will opt to use only the $T/T'$ notation in this paper. In what follows, we will often present the symmetry diagrams for different two-colour wallpaper groups; these are based on those found in \emph{Isometrica}.

\subsection{Compatible two-colour wallpapers}

As we move through this section, we will be proving a series of lemmas indicating which two-colour wallpaper groups are incompatible with PGHPs. At the end of this section, Table~\ref{table:wallpapers} is provided to show the complete picture of compatible and incompatible two-colour wallpapers, along with the lemmas where incompatibility is addressed.

We know that any wallpaper with 3-fold or 6-fold symmetry is naturally incompatible, so we may flag these at the start. 

\begin{lemma} \label{lem:3_fold_wallpaper}
The two-colour wallpaper groups that include 3-fold or 6-fold rotational symmetry are incompatible with PGHPs; this includes $p31m/p3$,  $p3m1/p3$, $p6/p3$, $p6m/p3m1$, $p6m/p31m$, and $p6m/p6$.
\end{lemma}

Another quick way to identify incompatible two-colour wallpaper groups is to focus on the one-colour wallpaper results of Sen \& Martinez: that show $p4$, $p4g$, and $pgg$ are incompatible with PGHPs.

\begin{lemma} \label{lem:one_color_results}
Any two-colour wallpaper group whose colour-preserving wallpaper group is incompatible with PGHPs is also incompatible; this includes $pmg/pgg$, $cmm/pgg$, $p4/p4$, $p4m/p4$, $p4m/p4g$, $p4g/pgg$, and $p4g/p4$.
\end{lemma}

As we move forward, we will consider \emph{patches} we can obtain by cutting out a piece of a wallpaper pattern; ostensibly, these patches could be of any shape or size. The patch allows us a tool to look at a bounded section of a pattern. For our purposes, whenever we cut a patch, we will be cutting a square patch whose edges fall on lines of stitches; this allows us to describe the patch with two binary words and utilize our results on rosette symmetries. If a particular two-colour wallpaper group requires a patch to exist whose symmetries are unattainable by the results in Section~\ref{flip_symm_rosettes}, then we can eliminate that two-colour wallpaper group from the list of compatible wallpaper groups.

\begin{lemma} \label{lem:c4'}
The two-colour wallpaper groups $p4/p2$ and $p4g/cmm$ are incompatible with PGHPs.
\end{lemma}

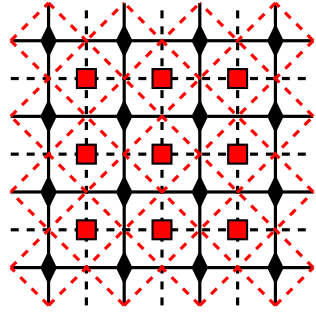
\begin{figure}[h!tbp]
\centering
\begin{tikzpicture}[scale=.5]
\foreach \i in {0,2,4,6}{
	\draw[very thick] (\i,-1)--(\i,7);
	\draw[very thick] (-1,\i)--(7,\i);
	}
\foreach \i in {1,3,5}{
	\draw[very thick,dashed] (\i,-1)--(\i,7);
	\draw[very thick, dashed](-1,\i)--(7,\i);
	}
\foreach \i in {0,2,4,6}{
\draw[red, very thick, dashed] (-1,\i)--(6-\i,7);
\draw[red, very thick, dashed] (\i,-1)--(7,6-\i);
\draw[red, very thick, dashed] (7,\i)--(\i,7);
\draw[red, very thick, dashed] (\i,-1)--(-1,\i);
}
\foreach \k in {0,2,4,6}{
	\foreach \j in {0,2,4,6}{
		\draw (\j,\k)node[2fold]{};
		}
	}
	
\foreach \k in {1,3,5}{
	\foreach \j in {1,3,5}{
		\draw (\j,\k)node[4foldcolorswap]{};
		}
	}

\end{tikzpicture}
\caption{The symmetry diagram for $p4g/cmm$}
\label{fig:p4gcmm}
\end{figure}

\begin{proof}
These two-colour wallpaper groups have points of flip-four-fold rotational symmetry with no intersecting lines of reflection or flip-reflection. A patch that is centered at this kind of flip-four-fold center will have rosette symmetry $c4'$ which is not a possible symmetry type in Table~\ref{rosettes}. 

For $p4/p2$, it is quick to see from notation alone that there is flip-four-fold rotation with no lines of reflection. For $p4g/cmm$, we provide the symmetry diagram in Figure~\ref{fig:p4gcmm}.
\end{proof}

%

\begin{lemma} \label{lem:d'2}
The two-colour wallpaper groups $pmm/p2$, $pmm/pmg$, $pmm/cmm$, $cmm/p2$, and $p4m/cmm$ are incompatible with PGHPs.
\end{lemma}

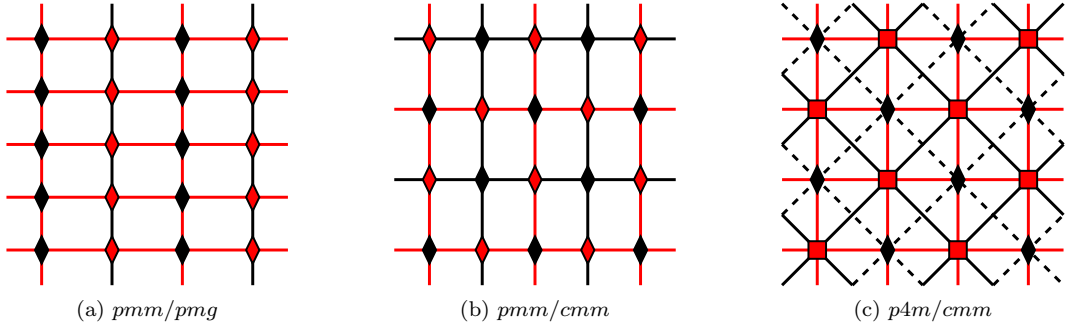
\begin{figure}[h!tbp]
\centering

\subfloat[$pmm/pmg$]{%
\resizebox*{4cm}{!}{
\begin{tikzpicture}[scale=.5]
\foreach \i in {0,4}{
	\draw[very thick,red] (\i,-1)--(\i,7);
	\draw[very thick] (\i+2,-1)--(\i+2,7);
	}
\foreach \j in {0,1.5,3,4.5,6}{
	\draw[very thick,red] (-1,\j)--(7,\j);
	}

\foreach \k in {0,1.5,3,4.5,6}{
	\foreach \j in {0,4}{
		\draw (\j+2,\k)node[2foldcolorswap]{};
		\draw (\j,\k)node[2fold]{};
		}
	}

\end{tikzpicture}
}} 
\hspace{1cm}
\subfloat[$pmm/cmm$]{%
\resizebox*{4cm}{!}{
\begin{tikzpicture}[scale=.5]
\foreach \i in {0,3}{
	\draw[very thick,red] (\i,-1)--(\i,7);
	\draw[very thick] (\i+1.5,-1)--(\i+1.5,7);
	}
\draw[very thick,red] (6,-1)--(6,7);
\foreach \i in {0,4}{
	\draw[very thick,red] (-1,\i)--(7,\i);
	\draw[very thick] (-1,\i+2)--(7,\i+2);
	}

\foreach \k in {0,4}{
	\foreach \j in {0,3}{
		\draw (\j+1.5,\k)node[2foldcolorswap]{};
		\draw (\j,\k)node[2fold]{};
		\draw (\j+1.5,\k+2)node[2fold]{};
		\draw (\j,\k+2)node[2foldcolorswap]{};
		}
	\draw (6,\k+2)node[2foldcolorswap]{};
	\draw (6,\k)node[2fold]{};
	}

\end{tikzpicture}
}} 
\hspace{1cm}
\subfloat[$p4m/cmm$]{%
\resizebox*{4cm}{!}{
\begin{tikzpicture}[scale=.5]
\foreach \i in {0,2,4,6}{
	\draw[very thick,red] (\i,-1)--(\i,7);
	}
\foreach \j in {0,2,4,6}{
	\draw[very thick,red] (-1,\j)--(7,\j);
	}
\draw[very thick] (-1,-1)--(7,7);
\draw[very thick] (3,-1)--(7,3);
\draw[very thick] (-1,3)--(3,7);
\draw[very thick,dashed] (1,-1)--(7,5);
\draw[very thick,dashed] (-1,1)--(5,7);
\draw[very thick,dashed] (5,-1)--(7,1);
\draw[very thick,dashed] (-1,5)--(1,7);

\draw[very thick, dashed] (-1,7)--(7,-1);
\draw[very thick, dashed] (-1,3)--(3,-1);
\draw[very thick, dashed] (3,7)--(7,3);
\draw[very thick] (1,7)--(7,1);
\draw[very thick] (-1,5)--(5,-1);
\draw[very thick] (-1,1)--(1,-1);
\draw[very thick] (5,7)--(7,5);

\foreach \k in {0,4}{
	\foreach \j in {0,4}{
		\draw (\k,\j)node[4foldcolorswap]{};
		\draw (\k+2,\j+2)node[4foldcolorswap]{};
		\draw (\k,\j+2)node[2fold]{};
		\draw (\k+2,\j)node[2fold]{};
		}
	}

\end{tikzpicture}
}} 
\caption{The symmetry diagrams for the proof of Lemma~\ref{lem:d'2}.}
\label{fig:d'2}
\end{figure}

\begin{proof}
Each of these wallpaper groups include a center of two-fold rotational symmetry with two flip-reflection lines passing through them; a patch centered at such a location would have rosette type $d'2$, which is incompatible with GHPs (Table~\ref{rosettes}).

In the cases of $pmm/p2$ and $cmm/p2$, the presence of flip-reflections in two directions is evident since lines of reflection are contained in the symmetries consistent with colour but not with the symmetries that preserve colour. For $pmm/pmg$, $pmm/cmm$, and $p4m/cmm$, this fact is less obvious; therefore, the symmetry diagrams for each type have been provided in Figure~\ref{fig:d'2}.
\end{proof}

\begin{lemma} \label{lem:pmmpmm}
The two-colour wallpaper group $pmm/pmm$ is incompatible with PGHPs.
\end{lemma}

\begin{proof}
The symmetry diagram for $pmm/pmm$ is provided in Figure~\ref{fig:pmmpmm}. Notice that there exist intersections with $d2$ symmetry; Table~\ref{rosettes} shows that a patch centered at such a location must have $\{e,r_{180}, R_H, R_V\}$ as the set of symmetries.

However, $pmm/pmm$ also has lines of flip-reflective symmetry that are parallel to one direction of reflection lines. This creates a contradiction by Lemma~\ref{not_symm}, since $\overline{R_H}$ and $\overline{R_V}$ can never be symmetries.
\end{proof}
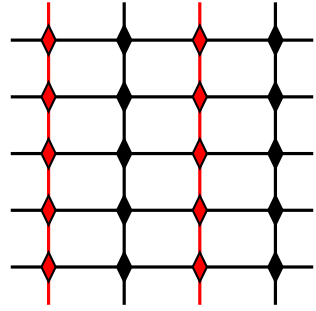
\begin{figure}[h!tbp]
\centering
\begin{tikzpicture}[scale=.5]
\foreach \i in {0,4}{
	\draw[very thick,red] (\i,-1)--(\i,7);
	\draw[very thick] (\i+2,-1)--(\i+2,7);
	}
\foreach \j in {0,1.5,3,4.5,6}{
	\draw[very thick] (-1,\j)--(7,\j);
	}

\foreach \k in {0,1.5,3,4.5,6}{
	\foreach \j in {0,4}{
		\draw (\j+2,\k)node[2fold]{};
		\draw (\j,\k)node[2foldcolorswap]{};
		}
	}

\end{tikzpicture}
\caption{The symmetry diagram for $pmm/pmm$}
\label{fig:pmmpmm}
\end{figure}

Moving forward, we will need a result from Sen \& Martinez concerning the partitioning of PGHPs into identical tiles:

\begin{lemma} \label{tiles}
\citep{Sen_Martinez} A periodic generalized hitomezashi pattern can be viewed as an edge-to-edge tiling of identical, rectangular tiles with even-length dimensions. Further, there exists a smallest such tile so that every other such tiling will have rectangles with length and width greater than or equal to the smallest tile.
\end{lemma}

\begin{figure}[h!tbp]
\centering
\subfloat[]{%
\resizebox*{6cm}{!}{\includegraphics{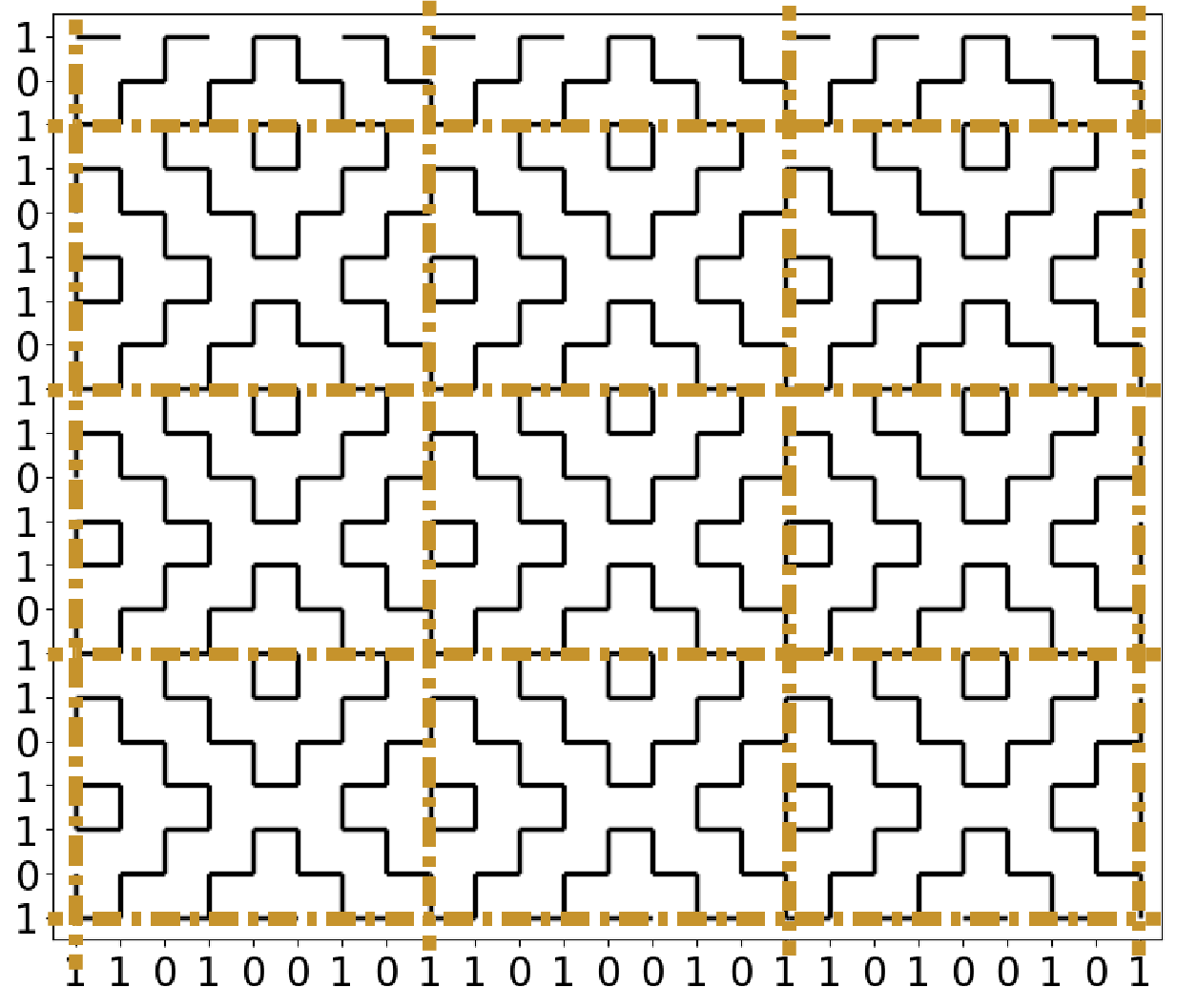}}} \hspace{1cm}
\subfloat[]{%
\resizebox*{6cm}{!}{\includegraphics{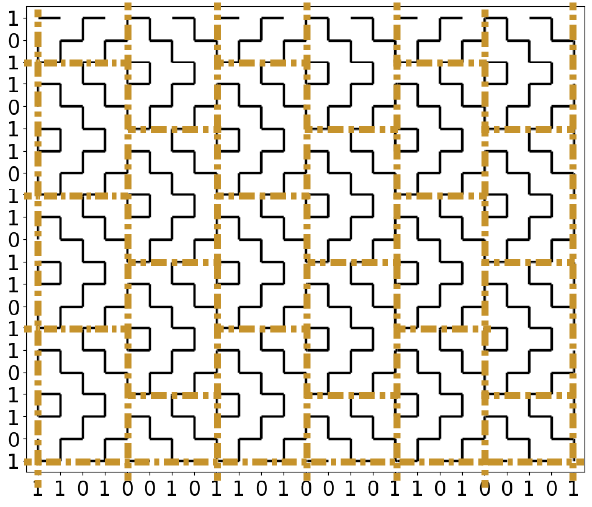}}}
\caption{Figure (a) shows the tiles of a PGHP using Lemma~\ref{tiles}, while Figure (b) shows an alternate way the PGHP could be partitioned, but in such a way that the tiling is not edge-to-edge.} 
\label{tiles_ex}
\end{figure}

We can apply Lemma~\ref{tiles} to the persimmon design that is a traditional hitomezashi design. The persimmon wallpaper in Figure~\ref{tiles_ex} is generated by $x=(11010010)^{\infty}$ and $y=(101)^{\infty}$. Note that the period of $x$ is 8 and the period of $y$ is 3. By Lemma~\ref{tiles}, the dimension of the edge-to-edge tiling is $8 \times 6$, which is shown in Figure~\ref{tiles_ex}(a). Notice that we can partition the persimmon design into smaller tiles, as in Figure~\ref{tiles_ex}(b). However, this tiling is not edge-to-edge, and is not an example of the kind of tiling covered by Lemma~\ref{tiles}.

The alternating nature of hitomezashi stitching leads to a series of restrictions for how reflections and glide-reflections relate to each other in PGHPs. We will develop a series of lemmas detailing these restrictions and use them to identify incompatible two-colour wallpaper groups. In Sen \& Martinez, the following lemma was proven.

\begin{lemma} \label{diagonal_glide}
\citep{Sen_Martinez} If a PGHP has a line of glide-reflective symmetry with slope $\pm 1$, then it also has lines of reflection with the same slope.
\end{lemma}

One caveat is that the proof of this lemma only included glide-reflection lines that pass through the intersections of the stitching grid. There is an additional case to consider where those lines do \emph{not} pass through intersections. We are able to rectify this oversight here. We note that while this case was not considered by Sen \& Martinez, it also doesn't affect the ultimate result this lemma was used for: that $pgg$ is incompatible with PGHPs.

The correctly restated version of the Sen \& Martinez lemma is included as part of the following lemma.

\begin{lemma} \label{diagonal_glide2}
\begin{enumerate}
\item If a PGHP has a line of glide-reflective symmetry with slope $\pm 1$ that passes through intersections of the stitching grid, then it also has lines of reflection with the same slope.
\item If a PGHP has a line of flip-glide-reflective symmetry with slope $\pm 1$ that passes through intersections of the stitching grid, then it also has lines of flip-reflection with the same slope.
\end{enumerate}
\end{lemma}

Before proving this lemma, we introduce a fact that will be useful in any of the subsequent lemmas featuring diagonal reflections or glide-reflections (flipped or not).

\begin{lemma} \label{lem:squares}
If a PGHP has a line of glide-reflective, flip-glide-reflective, reflective, or flip-reflective symmetry whose slope is $\pm 1$, then the tiling with the smallest tiles as described in Lemma~\ref{tiles} must be an edge-to-edge tiling of $2k \times 2k$ squares for some integer $k$.
\end{lemma}

\begin{proof}
The square dimensions follow since a reflection or flip-reflection over a diagonal line will exchange the length and width of the tiles. Since the dimensions of the tiling are unique, the tiles must be squares.
\end{proof}

We proceed to the proof of Lemma~\ref{diagonal_glide2}.

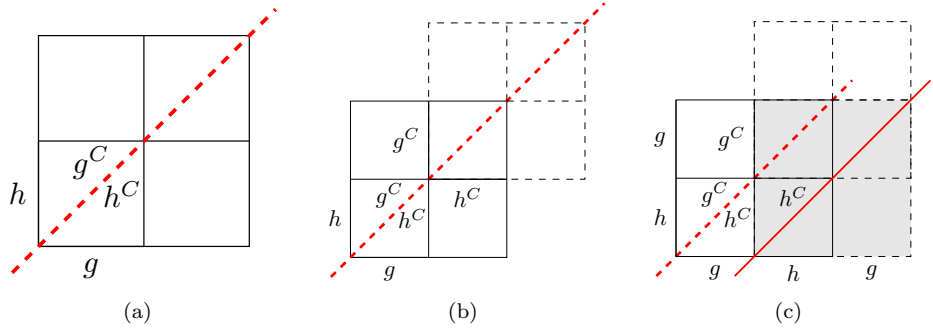
\begin{figure}[h!tbp]
\centering
\subfloat[]{%
\resizebox*{4cm}{!}{
\begin{tikzpicture}[scale=1.2]
\draw(0,0) rectangle (2,2);
\draw (0,0)--(0.5,0)node[below]{\small$g$}--(1,0);
\draw (0,0)--(0,0.5)node[left]{\small$h$}--(0,1);
\draw (0,1)--(0.5,1)node[yshift=-7]{\small$g^C$}--(2,1);
\draw (1,0)--(1,0.5)node[xshift=-7]{\small$h^C$}--(1,2);
\draw[very thick,dashed,red] (-0.25,-0.25)--(2.25,2.25);
\end{tikzpicture}
}}\hspace{5pt}
\subfloat[]{%
\resizebox*{4cm}{!}{
\begin{tikzpicture}[scale=1.2]
\draw(0,0) rectangle (2,2);
\draw (0,0)--(0.5,0)node[below]{\small$g$}--(1,0);
\draw (0,0)--(0,0.5)node[left]{\small$h$}--(0,1);
\draw (0,1)--(0.5,1)node[yshift=-7]{\small$g^C$}--(2,1);
\draw (1,0)--(1,0.5)node[xshift=-7]{\small$h^C$}--(1,2);
\draw[dashed] (1,1) rectangle (3,3);
\draw[very thick,red,dashed] (-0.25,-0.25)--(3.25,3.25);
\draw (1,1)--(1.5,1)node[below]{\small$h^C$}--(2,1);
\draw (1,1)--(1,1.5)node[left]{\small$g^C$}--(1,2);
\draw[dashed] (2,2)--(2,3);
\draw[dashed] (2,2)--(3,2);
\end{tikzpicture}
}}\hspace{5pt}
\subfloat[]{%
\resizebox*{4cm}{!}{
\begin{tikzpicture}[scale=1.2]
\draw[dashed,fill=black!10!white] (1,0) rectangle (3,2);
\draw[very thick,red,dashed] (-0.25,-0.25)--(2.25,2.25);
\draw(0,0) rectangle (2,2);

\draw (0,0)--(0.5,0)node[below]{\small$g$}--(1.5,0)node[yshift=-7]{\small$h$}--(2,0);
\draw (0,0)--(0,0.5)node[left]{\small$h$}--(0,1.5)node[xshift=-7]{\small$g$}--(0,2);
\draw (0,1)--(0.5,1)node[yshift=-7]{\small$g^C$}--(1.5,1)node[below]{\small$h^C$}--(2,1);
\draw (1,0)--(1,0.5)node[xshift=-7]{\small$h^C$}--(1,1.5)node[left]{\small$g^C$}--(1,2);
\draw[dashed] (1,1) rectangle (3,3);
\draw[thick,red] (.75,-.25)--(3.25,2.25);
\draw[dashed] (2.45,0)--(2.5,0)node[below]{\small$g$};
\draw[dashed] (2,2)--(2,3);
\draw[dashed] (2,2)--(3,2);
\end{tikzpicture}
}}

\caption{Images of one of our $2k \times 2k$ lattice squares for the proof of Lemma~\ref{diagonal_glide}. This displays the case where $k$ is odd.}
\label{fig:flip-glide_odd}
\end{figure}

\begin{proof}
Part (1) was already proven by Sen \& Martinez, so we focus on arguing (2). Note that the proofs of (1) and (2) are almost identical with small label changes, so the reader can quickly recreate the Sen \& Martinez argument from the work here.

We only consider the case when a PGHP has a line of flip-glide-reflection or flip-reflection with slope $+1$, since the case with slope equal to $-1$ is symmetric.

Consider one of the $2k \times 2k$ square tiles guaranteed by Lemma~\ref{lem:squares}, positioning it so that the line of flip-glide-reflection is its diagonal. The flip-glide-reflection will have glide distance that is half the length of the diagonal of the square. For now, we will assume this fact; at the end of this proof we return to the issue of glide distance.

Our argument is illustrated in Figure~\ref{fig:flip-glide_odd}. We break our $2k \times 2k$ tile into four $k \times k$ squares that share borders; our flip-glide-reflection will map the lower left square onto the upper right square. Each quarter can be generated with two binary strings of length $k+1$. Let the bottom-left square be generated by $g,h \in \{0,1\}^{k+1}$.

By the alternating nature of our stitches, we can label the edges of the other squares that lie parallel. This analysis must be broken into two cases depending on the parity of $k$. When $k$ is even, the label remains unchanged on the parallels; when $k$ is odd, the labels are now the complement word. We will continue this proof assuming $k$ is odd; the proof when $k$ is even proceeds similarly.

Figure~\ref{fig:flip-glide_odd}(a) shows the labeling for the case when $k$ is odd. When we apply the flip-glide-reflection symmetry, we can find more labels for our square as shown in Figure~\ref{fig:flip-glide_odd}(b); after the translation along the diagonal, we exchange the vertical and horizontal words and then complement them. Now we can see that the $2k \times 2k$ square is generated with the words $gh$ (horizontal word) and $hg$ (vertical word); it must be the case that the last entry of $g$ is the first entry of $h$ (and vice-versa), since these entries overlap. If we now shift our focus to the $2k \times 2k$ square that is $k$ units to the right (shaded in Figure~\ref{fig:flip-glide_odd}(c)), we are considering a square generated by $hg$ and $h^Cg^C=(hg)^C$. This indicates there is a line of flip-reflective symmetry passing through the diagonal of the shaded square. This flip-reflection will continue through the entire wallpaper design, since the shaded square tiles the whole wallpaper.

Now we return to the issue of glide-distance. Composing the glide-reflection with itself is a translation. We already know that shifting the design up $2k$ units and right $2k$ units is a translation symmetry; the smallest translation symmetry that exists in this direction will divide $2k$. Having a translation smaller than shifting up $k$ units and right $k$ units will contradict the minimality of the $2k \times 2k$ tile, so we don't need to consider these cases. If we suppose we have a translation symmetry the shifts $k$ units up and to the right then: when $k$ is even we have $h=g$ and find that our tile was not minimal; or, when $k$ is odd, we have $h=g^C$ and a line of flip-reflective symmetry instead of flip-glide-reflective symmetry on the diagonal. Neither of these options is allowable, so the translation $k$ units up and to the right is not a symmetry. Hence our flip-glide-reflection has glide distance equal to half the distance of the diagonal of the $2k \times 2k$ square.
\end{proof}

Now we work on the case where the glide-reflection lines do not pass through intersections of the square grid.

\begin{lemma} \label{diagonal_glide3}
\begin{enumerate}
\item If a PGHP has a line of glide-reflective symmetry with slope $\pm 1$ that does not pass through intersections of the stitching grid, then it also has lines of flip-reflection with the same slope.
\item If a PGHP has a line of flip-glide-reflective symmetry with slope $\pm 1$ that does not pass through intersections of the stitching grid, then it also has lines of reflection with the same slope.
\end{enumerate}
\end{lemma}

\begin{figure}[h!tbp]
\centering
\subfloat[Illustration of the glide-reflection process]{%
\resizebox*{12cm}{!}{
\begin{tikzpicture}[scale=1.2]
\draw(0,0) rectangle (2,2);
\draw (0,0)--(1,0);
\draw (0,0)--(0,1);
\draw (0,1)--(2,1);
\draw (1,0)--(1,2);
\draw[fill=gray!40!white] (0,0) rectangle (1,1);
\draw[->,ultra thick,dotted] (0,0)--(0.5,0);
\draw[->,ultra thick] (0,0)--(0,0.5);
\draw[very thick,dashed] (0.25,-0.25)--(2.25,1.75);
\draw (2.5,1)node{\large$\mapsto$};
\end{tikzpicture}
\begin{tikzpicture}[scale=1.2]
\draw(0,0) rectangle (2,2);
\draw (0,0)--(1,0);
\draw (0,0)--(0,1);
\draw (0,1)--(2,1);
\draw (1,0)--(1,2);
\draw[fill=gray!40!white] (0.5,0.5) rectangle (1.5,1.5);
\draw[->,ultra thick,dotted] (0.5,0.5)--(1,0.5);
\draw[->,ultra thick] (0.5,0.5)--(0.5,1);
\draw[very thick,dashed] (0.25,-0.25)--(2.25,1.75);
\draw (2.5,1)node{\large$\mapsto$};
\end{tikzpicture}
\begin{tikzpicture}[scale=1.2]
\draw(0,0) rectangle (2,2);
\draw (0,0)--(1,0);
\draw (0,0)--(0,1);
\draw (0,1)--(2,1);
\draw (1,0)--(1,2);
\draw[fill=gray!40!white] (1,0) rectangle (2,1);
\draw[->,ultra thick,dotted] (1,0)--(1,0.5);
\draw[->,ultra thick] (1,0)--(1.5,0);
\draw[very thick,dashed] (0.25,-0.25)--(2.25,1.75);
\end{tikzpicture}
}}

\subfloat[Set-up of labels]{%
\resizebox*{4cm}{!}{
\begin{tikzpicture}[scale=1.2]
\draw(0,0) rectangle (2,2);
\draw (0,0)--(0.5,0)node[below]{\small$g$}--(2,0);
\draw (0,0)--(0,0.5)node[left]{\small$h$}--(0,1);
\draw (0,1)--(0.5,1)node[yshift=-7]{\small$g^C$}--(2,1);
\draw (1,0)--(1,0.5)node[xshift=-7]{\small$h^C$}--(1,2);
\draw[very thick,dashed] (0.25,-0.25)--(2.25,1.75);
\end{tikzpicture}
}}\hspace{5pt}
\subfloat[Labels required for glide-reflection symmetry]{%
\resizebox*{4cm}{!}{
\begin{tikzpicture}[scale=1.2]
\draw(0,0) rectangle (2,2);
\draw (0,0)--(0.5,0)node[yshift=-8]{\small$g$}--(1.5,0)node[yshift=-7]{\small$g^C$}--(2,0);
\draw (0,0)--(0,0.5)node[left]{\small$g^C$}--(0,1.5)node[left]{\small$g$}--(0,2);
\draw (0,1)--(0.5,1)node[yshift=-7]{\small$g^C$}--(1.5,1)node[yshift=-5]{\small$g$}--(2,1);
\draw (1,0)--(1,0.5)node[xshift=-5]{\small$g$}--(1,1.5)node[xshift=-7]{\small$g^C$}--(1,2);
\draw[very thick,dashed] (0.25,-0.25)--(2.25,1.75);
\end{tikzpicture}
}}\hspace{5pt}
\subfloat[Existence of flip-reflection]{%
\resizebox*{4cm}{!}{
\begin{tikzpicture}[scale=1.2]
\draw(0,0) rectangle (2,2);
\draw (0,0)--(0.5,0)node[yshift=-8]{\small$g$}--(1.5,0)node[yshift=-7]{\small$g^C$}--(2,0);
\draw (0,0)--(0,0.5)node[left]{\small$g^C$}--(0,1.5)node[left]{\small$g$}--(0,2);
\draw (0,1)--(2,1);
\draw (1,0)--(1,2);
\draw[very thick,dashed] (0.25,-0.25)--(2.25,1.75);
\draw[thick,red] (-.25,-.25)--(2.25,2.25);
\end{tikzpicture}
}}

\subfloat[Labels required for flip-glide-reflection symmetry]{%
\resizebox*{4cm}{!}{
\begin{tikzpicture}[scale=1.2]
\draw(0,0) rectangle (2,2);
\draw (0,0)--(0.5,0)node[yshift=-8]{\small$g$}--(1.5,0)node[yshift=-7]{\small$g^C$}--(2,0);
\draw (0,0)--(0,0.5)node[left]{\small$g$}--(0,1.5)node[left]{\small$g^C$}--(0,2);
\draw (0,1)--(0.5,1)node[yshift=-7]{\small$g^C$}--(1.5,1)node[yshift=-5]{\small$g$}--(2,1);
\draw (1,0)--(1,0.5)node[xshift=-7]{\small$g^C$}--(1,1.5)node[xshift=-7]{\small$g$}--(1,2);
\draw[very thick,dashed,red] (0.25,-0.25)--(2.25,1.75);
\end{tikzpicture}
}}\hspace{5pt}
\subfloat[Existence of reflection]{%
\resizebox*{4cm}{!}{
\begin{tikzpicture}[scale=1.2]
\draw(0,0) rectangle (2,2);
\draw (0,0)--(0.5,0)node[yshift=-8]{\small$g$}--(1.5,0)node[yshift=-7]{\small$g^C$}--(2,0);
\draw (0,0)--(0,0.5)node[left]{\small$g$}--(0,1.5)node[left]{\small$g^C$}--(0,2);
\draw (0,1)--(2,1);
\draw (1,0)--(1,2);
\draw[very thick,dashed,red] (0.25,-0.25)--(2.25,1.75);
\draw[thick] (-.25,-.25)--(2.25,2.25);
\end{tikzpicture}
}}

\caption{Images illustrating the proof of Lemma~\ref{diagonal_glide3}.}
\label{fig:glide_off-grid}
\end{figure}
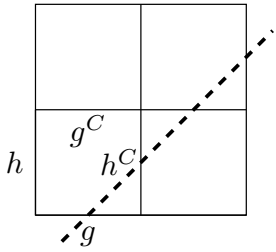
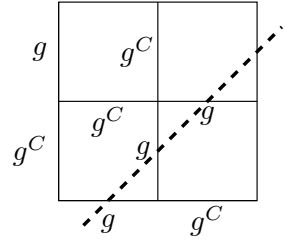
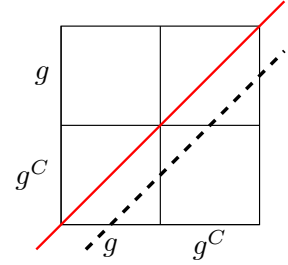
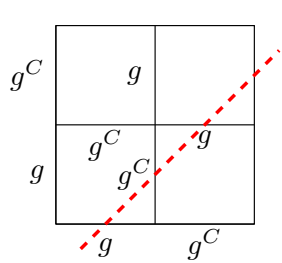
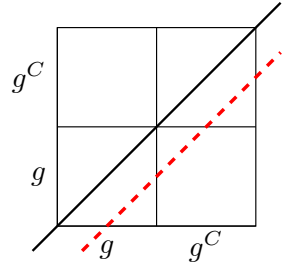

\begin{proof}
We consider the case where the slope is $+1$, since the $-1$ case is symmetric. We know the minimal tile for this PGHP has dimensions $2k \times 2k$; we will further argue that $k$ must be odd.

Since the glide-reflection does not pass through grid intersections, it must map a point at the top/bottom of a grid square to a point on the side of a grid square before the reflection is performed. This will ensure that the square grid is realigned after the glide-reflection is completed. Therefore the glide distance is $j+0.5$ units up and to the right for some integer $j$. Since composing a glide-reflection with itself is a translation, we also know that the translation $2j+1$ up and $2j+1$ units to the right is a symmetry. We know the square tile cannot have odd dimensions by Lemma~\ref{tiles}, so it is the case that $2j+1=k$ and the dimension of the square tile is $2(2j+1) \times 2(2j+1)$. Hence $k$ is odd. This argument works for the flip-glide-reflection as well.

Now consider the diagram created when our square tiles are positioned so that a line of glide-reflection passes through the bottom of the tile $j+0.5$ units from the left (or rather, so that 1/4 of the bottom border is to the left). We can divide our square tile into four $k \times k$ squares. The glide reflection maps the bottom left square onto the bottom right square where the left border becomes the bottom border (and vice-versa), which is illustrated in Figure~\ref{fig:glide_off-grid}(a) (the arrows in this figure are used to mark the sides of the shaded square). 

Let the bottom-left square be generated by $g,h \in \{0,1\}^{k+1}$. By the alternating nature of our stitches, we can label the edges of the other squares that lie parallel. Since $k$ is odd, the adjacent, parallel edges will be labeled with the complement word, as seen in Figure~\ref{fig:glide_off-grid}(b). In order for the glide-reflection to be a symmetry, we must have $h^C=g$, so we can edit our diagram to the one in Figure~\ref{fig:glide_off-grid}(c). In this diagram, we have also labeled the top right square with the appropriate labels to make the inverse of the glide-reflection a symmetry and populated the parallel edges with labels. One subtlety of this set-up is that the last entry of $g$ must be the first entry of $g^C$, and vice-versa, since these entries overlap.

Now we can see that the $2k \times 2k$ square is generated by $gg^C$ and $g^Cg=(gg^C)^C$. This indicates there is a line of flip-reflective symmetry passing through the diagonal of the square tile.

For the case of the flip-glide-reflection, we can use the same diagrams with small alterations. When the flip-glide-reflection is performed, we must complement the words. The updated figures are in Figure~\ref{fig:glide_off-grid}(e) and (f), showing the existence of the line of reflection.
\end{proof}

The above lemma means that it is possible to obtain a type $pg$ wallpaper in a PGHP where the glide-reflection lines are diagonal. This calls into question the assertion of Sen \& Martinez, that you cannot have diagonal glide-reflection lines in two directions without reflection lines. To address this, we prove the following lemma, making use of our two-colour rosette results.

\begin{lemma} \label{diagonal_glides4}
If a PGHP has lines of glide-reflective symmetry with both slope 1 and slope -1, then there must also exist lines of reflection.
\end{lemma}

\begin{proof}
If the glide-reflections pass through the intersections of the grid, we automatically have lines of reflection by Lemma~\ref{diagonal_glide2}. If the glide-reflections do not pass through the intersections of the grid, then by Lemma~\ref{diagonal_glide3}, there must be flip-reflections with slope 1 and slope -1 as well. Where these lines intersect, there will be $180^{\circ}$-rotational symmetry. If there were no intersecting lines of reflection at this point, then a patch centered at the intersection would have type $d'2$. We know this is not possible by Table~\ref{rosettes}. Therefore there must also be lines of vertical and horizontal reflection at this point.
\end{proof}

With Lemma~\ref{diagonal_glides4}, the proof that $pgg$ is incompatible with PGHPs in Sen \& Martinez remains sound, since it bases its argument on horizontal and vertical lines of glide-reflection, asserting that diagonal lines of glide-reflection will also require the existence of lines of reflection.

The preceding lemmas can be used to identify two incompatible two-colour wallpaper types:

\begin{lemma} \label{diagonal_glide_ref_flip}
The two-colour wallpaper groups $pmg/pg$ and $pmg/p2$ are incompatible with PGHPs.
\end{lemma}

\begin{figure}[h!tbp]
\centering

\subfloat[$pmg/pg$]{%
\resizebox*{3.5cm}{!}{
\begin{tikzpicture}[scale=.5]
\foreach \i in {0,3,6}{
	\draw[very thick,red] (-1,\i-1)--(7-\i,7);
	\draw[very thick,dashed] (7,\i-1)--(\i-1,7);
	}
\foreach \i in {1,4}{
	\draw[very thick,red] (\i+1,-1)--(7,5-\i);
	\draw[very thick,dashed] (\i,-1)--(-1,\i);
	}
\foreach \i in {0,3,6}{
	\foreach \j in {0,3,6}{
		\draw (\i+.75,\j-.75)node[2foldcolorswap]{};
		}
		}
\foreach \i in {-2,1,4}{
	\foreach \j in {1,4,7}{
		\draw (\i+1.25,\j-.25)node[2foldcolorswap]{};
		}
		}
\end{tikzpicture}
}} \hspace{1cm}
\subfloat[$pmg/p2$]{%
\resizebox*{3.5cm}{!}{
\begin{tikzpicture}[scale=.5]
\foreach \i in {0,3,6}{
	\draw[very thick,red] (-1,\i-1)--(7-\i,7);
	\draw[very thick,dashed,red] (7,\i-1)--(\i-1,7);
	}
\foreach \i in {1,4}{
	\draw[very thick,red] (\i+1,-1)--(7,5-\i);
	\draw[very thick,dashed,red] (\i,-1)--(-1,\i);
	}
\foreach \i in {0,3,6}{
	\foreach \j in {0,3,6}{
		\draw (\i+.75,\j-.75)node[2fold]{};
		}
		}
\foreach \i in {-2,1,4}{
	\foreach \j in {1,4,7}{
		\draw (\i+1.25,\j-.25)node[2fold]{};
		}
		}
\end{tikzpicture}}}
\caption{Symmetry diagrams for Lemma~\ref{diagonal_glide_ref_flip}}
\label{fig:flip_glide_lemma_results}
\end{figure}
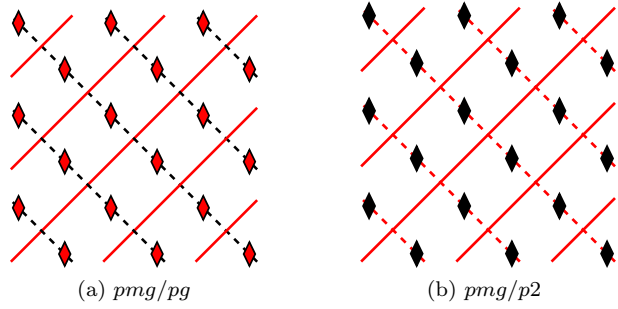

\begin{proof}

Any wallpaper group that includes flip-reflective symmetry needs those lines to have slope $\pm 1$. The symmetry diagrams for $pmg/pg$ and $pmg/p2$ are provided in Figure~\ref{fig:flip_glide_lemma_results}. Both these types require lines of glide-reflection or flip-glide-reflection on the diagonal with no parallel lines of reflection or flip-reflection. This is impossible by the results of Lemmas~\ref{diagonal_glide2} and \ref{diagonal_glide3}.
\end{proof}

We will be able to identify more incompatible two-colour wallpaper groups by analysing the behaviour surrounding flip-reflections. To that end, we include the following lemma.

\begin{lemma}
\label{diagonal_reflection_flip}
If a PGHP has a line of flip-reflective symmetry with slope $\pm 1$ then one of the following is true:
\begin{enumerate}
\item the PGHP has parallel lines of flip-glide-reflective symmetry,
\item the PGHP has parallel lines of reflective symmetry and a flip-translation symmetry in a parallel direction, or
\item the PGHP has parallel lines of glide-reflective symmetry.
\end{enumerate}
\end{lemma}

\begin{figure}[h!tbp]
\centering
\subfloat[Odd case]{%
\resizebox*{5cm}{!}{
\begin{tikzpicture}[scale=1.2]
\draw[dashed,fill=black!10!white] (1,0) rectangle (3,2);
\draw(0,0) rectangle (2,2);
\draw (0,0)--(0.5,0)node[below]{\small$g$}--(1.5,0)node[below]{\small$h$}--(1,0);
\draw (0,0)--(0,0.5)node[left]{\small$g^C$}--(0,1.5)node[left]{\small$h^C$}--(0,1);
\draw (0,1)--(2,1);
\draw (1,0)--(1,0.5)node[xshift=-7]{\small$g$}--(1,1.5)node[xshift=-7]{\small$h$}--(1,2);
\draw (2,0)--(2,1.5)node[xshift=-7]{\small$h^C$};
\draw[very thick,red] (-0.25,-0.25)--(3.25,3.25);
\draw[dashed] (2.45,0)--(2.5,0)node[yshift=-7]{\small$g$};
\draw[dashed] (2,1)--(2.5,1)node[below]{\small$g^C$}--(3.5,1)node[below]{\small$h^C$}--(4,1);
\draw[dashed] (2,2)--(2,2.5)node[xshift=-7]{\small$g^C$}--(2,3);
\draw[dashed] (3,1)--(3,3);
\draw[thick,red,dashed] (.75,-.25)--(4.25,3.25);
\draw[dashed] (2,3)--(4,3)--(4,1);
\draw[dashed] (3,2)--(4,2);
\end{tikzpicture}
}}
\subfloat[Odd case with $h=g^C$]{%
\resizebox*{5cm}{!}{
\begin{tikzpicture}[scale=1.2]
\draw[dashed,fill=black!10!white] (1,0) rectangle (3,2);
\draw(0,0) rectangle (2,2);
\draw (0,0)--(0.5,0)node[below]{\small$g$}--(1.5,0)node[below]{\small$g^C$}--(1,0);
\draw (0,0)--(0,0.5)node[left]{\small$g^C$}--(0,1.5)node[left]{\small$g$}--(0,1);
\draw (0,1)--(2,1);
\draw (1,0)--(1,0.5)node[xshift=-7]{\small$g$}--(1,1.5)node[xshift=-7]{\small$g^C$}--(1,2);
\draw (2,0)--(2,1.5)node[xshift=-7]{\small$g$};
\draw[very thick,red] (-0.25,-0.25)--(3.25,3.25);
\draw[dashed] (2.45,0)--(2.5,0)node[yshift=-7]{\small$g$};
\draw[dashed] (2,1)--(2.5,1)node[below]{\small$g^C$}--(3.5,1)node[below]{\small$g$}--(4,1);
\draw[dashed] (2,2)--(2,2.5)node[xshift=-7]{\small$g^C$}--(2,3);
\draw[dashed] (3,1)--(3,3);
\draw[very thick,red] (.75,-.25)--(4.25,3.25);
\draw[dashed] (2,3)--(4,3)--(4,1);
\draw[dashed] (3,2)--(4,2);
\draw[ultra thick, ->] (1,0)--(2,1);

\draw[very thick,dashed] (0.25,-0.25)--(3.75,3.25);
\end{tikzpicture}
}}
\subfloat[Odd case with $g=h$]{%
\resizebox*{5cm}{!}{
\begin{tikzpicture}[scale=1.2]
\draw[dashed,fill=black!10!white] (1,0) rectangle (3,2);
\draw(0,0) rectangle (2,2);
\draw (0,0)--(0.5,0)node[below]{\small$g$}--(1.5,0)node[below]{\small$g$}--(1,0);
\draw (0,0)--(0,0.5)node[left]{\small$g^C$}--(0,1.5)node[left]{\small$g^C$}--(0,1);
\draw (0,1)--(2,1);
\draw (1,0)--(1,0.5)node[xshift=-7]{\small$g$}--(1,1.5)node[xshift=-7]{\small$g$}--(1,2);
\draw (2,0)--(2,1.5)node[xshift=-7]{\small$g^C$};
\draw[very thick,red] (-0.25,-0.25)--(3.25,3.25);
\draw[dashed] (2.45,0)--(2.5,0)node[yshift=-7]{\small$g$};
\draw[dashed] (2,1)--(2.5,1)node[below]{\small$g^C$}--(3.5,1)node[below]{\small$g^C$}--(4,1);
\draw[dashed] (2,2)--(2,2.5)node[xshift=-7]{\small$g^C$}--(2,3);
\draw[dashed] (3,1)--(3,3);
\draw[very thick] (.75,-.25)--(4.25,3.25);
\draw[dashed] (2,3)--(4,3)--(4,1);
\draw[dashed] (3,2)--(4,2);
\draw[ultra thick,red,->] (1,0)--(2,1);
\end{tikzpicture}
}}

\subfloat[Even case]{%
\resizebox*{5cm}{!}{
\begin{tikzpicture}[scale=1.2]
\draw[dashed,fill=black!10!white] (1,0) rectangle (3,2);
\draw(0,0) rectangle (2,2);
\draw (0,0)--(0.5,0)node[below]{\small$g$}--(1.5,0)node[below]{\small$h$}--(1,0);
\draw (0,0)--(0,0.5)node[left]{\small$g^C$}--(0,1.5)node[left]{\small$h^C$}--(0,1);
\draw (0,1)--(2,1);
\draw (1,0)--(1,0.5)node[xshift=-7]{\small$g^C$}--(1,1.5)node[xshift=-7]{\small$h^C$}--(1,2);
\draw (2,0)--(2,1.5)node[xshift=-7]{\small$h^C$};
\draw[very thick,red] (-0.25,-0.25)--(3.25,3.25);
\draw[dashed] (2.45,0)--(2.5,0)node[yshift=-7]{\small$g$};
\draw[dashed] (2,1)--(2.5,1)node[below]{\small$g$}--(3.5,1)node[below]{\small$h$}--(4,1);
\draw[dashed] (2,2)--(2,2.5)node[xshift=-7]{\small$g^C$}--(2,3);
\draw[dashed] (3,1)--(3,3);
\draw[thick,red,dashed] (.75,-.25)--(4.25,3.25);
\draw[dashed] (2,3)--(4,3)--(4,1);
\draw[dashed] (3,2)--(4,2);
\end{tikzpicture}
}}
\subfloat[Even case with $h=g^C$]{%
\resizebox*{5cm}{!}{
\begin{tikzpicture}[scale=1.2]
\draw[dashed,fill=black!10!white] (1,0) rectangle (3,2);
\draw(0,0) rectangle (2,2);
\draw (0,0)--(0.5,0)node[below]{\small$g$}--(1.5,0)node[below]{\small$g^C$}--(1,0);
\draw (0,0)--(0,0.5)node[left]{\small$g^C$}--(0,1.5)node[left]{\small$g$}--(0,1);
\draw (0,1)--(2,1);
\draw (1,0)--(1,0.5)node[xshift=-7]{\small$g^C$}--(1,1.5)node[xshift=-7]{\small$g$}--(1,2);
\draw (2,0)--(2,1.5)node[xshift=-7]{\small$g$};
\draw[very thick,red] (-0.25,-0.25)--(3.25,3.25);
\draw[dashed] (2.45,0)--(2.5,0)node[yshift=-7]{\small$g$};
\draw[dashed] (2,1)--(2.5,1)node[below]{\small$g$}--(3.5,1)node[below]{\small$g^C$}--(4,1);
\draw[dashed] (2,2)--(2,2.5)node[xshift=-7]{\small$g^C$}--(2,3);
\draw[dashed] (3,1)--(3,3);
\draw[very thick] (.75,-.25)--(4.25,3.25);
\draw[ultra thick,red,->] (1,0)--(2,1);
\draw[dashed] (2,3)--(4,3)--(4,1);
\draw[dashed] (3,2)--(4,2);
\end{tikzpicture}
}}

\caption{Images of one of our $2k \times 2k$ lattice squares for the proof of Lemma~\ref{diagonal_reflection_flip}.}
\label{fig:flip-reflection}
\end{figure}
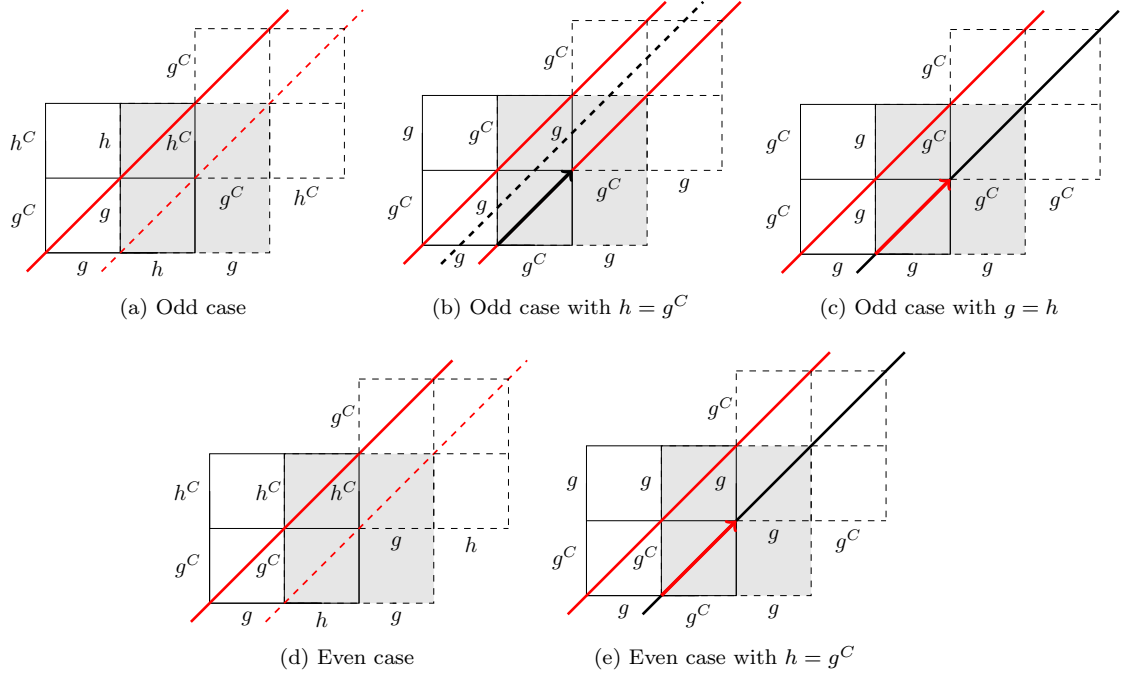

\begin{proof}
We again assume our line of flip-reflective symmetry has slope $+1$, since the $-1$ case is symmetric. We also consider the $2k \times 2k$ square that tiles the design by Lemma~\ref{tiles}. 

We proceed under the assumption that there is flip-reflective symmetry. To accomplish this, our $2k \times 2k$ tiles must be generated by $x$ and $y$ such that $x=y^C$. So, assume the square is generated by $gh$ and $(gh)^C=g^Ch^C$ where $g,h \in \{0,1\}^{k+1}$ and the first entry of $h$ is equal to the last entry of $g$ (and vice-versa).

The general case for $k$ odd is pictured in Figure~\ref{fig:flip-reflection}(a). In the case where $k$ is odd, the labels on parallel lines running through the $2k \times 2k$ square tile are the complement of what they were on the border. If we shift our focus to the $2k \times 2k$ square that is $k$ units to the right (the shaded square in Figure~\ref{fig:flip-reflection}(a)), we have a square generated by $hg$ and $gh$. This indicates there is a line of flip-glide-reflective symmetry passing through the diagonal of the shaded square---this situation aligns with that presented in Figure~\ref{fig:flip-glide_odd}. This is the situation described in part (1) of the lemma.

We can obtain stronger symmetries if we require additional restrictions on $h$ and $g$. There are two ways to strengthen the symmetries: the line of flip-glide-reflection could either be a line of flip-reflection, or a line of reflection (note that it could never be a line of glide-reflection, since this would require $h=h^C$). To produce the line of flip-reflection, we need $h=g^C$. To produce a line of reflection, we need $g=h$.
\begin{itemize}
\item If $h=g^C$, we have the scenario described in (3) of the lemma. This is pictured in Figure~\ref{fig:flip-reflection}(b), where the flip-glide-reflection is another line of flip-reflection. This situation aligns with the one pictured earlier in Figure~\ref{fig:glide_off-grid}(d), so there is additionally a line of glide-reflective symmetry and a translation symmetry shifting the design $k$ units up and to the right. Note that the glide-reflection doesn't pass through the intersections of the grid, so there is no scenario where this line is a reflection of any kind. It is also impossible to have a flip-translation symmetry parallel to the flip-glide-reflection lines. Such a translation could theoretically be halg the distance of the diagonal translation symmetry, but since $k$ is odd, this would not realign the square stitching grid.

\item If we instead set $g=h$, we get the scenario described in (2) of the lemma. The line of flip-glide-reflection is actually a line of reflection and there is a flip-translation symmetry shifting the design $k$ units to the right and $k$ units up as shown in Figure~\ref{fig:flip-reflection}(c).
\end{itemize}

If $k$ is instead even, we can again see the line of flip-glide-reflection that occurs in Figure~\ref{fig:flip-reflection}(d) that corresponds to (1) of the lemma. If we want the line of flip-glide reflection to be a regular line of reflection, we can set $h=g^C$, as pictured in Figure~\ref{fig:flip-reflection}(e). This yields parallel lines of reflection with a flip-translation symmetry $k$ units up and to the right and corresponds to part (2) of the lemma. To get another line of flip-reflection we would need $h=g$, which is not a possible scenario because then the minimal tile by Lemma~\ref{tiles} would be the $k \times k$ square.
\end{proof}

Lemma~\ref{diagonal_reflection_flip} can be used to identify some incompatible two-colour wallpapers.

\begin{lemma}  \label{lem:flip_reflection_lemma_results}
The two-colour wallpaper groups $pm/p1$, $pm/pm(m')$, $pm/pg$, and $pmm/pm$ are incompatible with PGHPs.
\end{lemma}

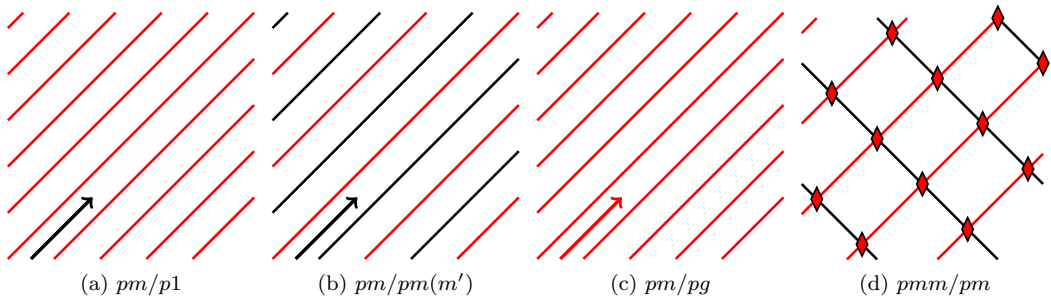
\begin{figure}[h!tbp]
\centering
\subfloat[$pm/p1$]{%
\resizebox*{3.5cm}{!}{
\begin{tikzpicture}[scale=.5]
\foreach \i in {0,3,6}{
	\draw[very thick,red] (-1,\i-1)--(7-\i,7);
	\draw[very thick,red] (-1,\i+.5)--(5.5-\i,7);
	}
\foreach \i in {1,4}{
	\draw[very thick,red] (\i-.5,-1)--(7,6.5-\i);
	\draw[very thick,red] (\i+1,-1)--(7,5-\i);
	}
\draw[ultra thick, ->] (-.25,-1)--(1.75,1);
\end{tikzpicture}
}} 
\subfloat[$pm/pm(m')$]{%
\resizebox*{3.5cm}{!}{
\begin{tikzpicture}[scale=.5]
\foreach \i in {0,3,6}{
	\draw[very thick,red] (-1,\i-1)--(7-\i,7);
	\draw[very thick] (-1,\i+.5)--(5.5-\i,7);
	}
\foreach \i in {1,4}{
	\draw[very thick] (\i-.5,-1)--(7,6.5-\i);
	\draw[very thick,red] (\i+1,-1)--(7,5-\i);
	}
\draw[ultra thick, ->] (-.25,-1)--(1.75,1);
\end{tikzpicture}
}} 
\subfloat[$pm/pg$]{%
\resizebox*{3.5cm}{!}{
\begin{tikzpicture}[scale=.5]
\foreach \i in {0,3,6}{
	\draw[very thick,red] (-1,\i-1)--(7-\i,7);
	\draw[very thick,red] (-1,\i+.5)--(5.5-\i,7);
	}
\foreach \i in {1,4}{
	\draw[very thick,red] (\i-.5,-1)--(7,6.5-\i);
	\draw[very thick,red] (\i+1,-1)--(7,5-\i);
	}
\draw[ultra thick, ->,red] (-.25,-1)--(1.75,1);
\end{tikzpicture}
}}
\subfloat[$pmm/pm$]{%
\resizebox*{3.5cm}{!}{
\begin{tikzpicture}[scale=.5]
\foreach \i in {0,3,6}{
	\draw[very thick,red] (-1,\i+.5)--(5.5-\i,7);	
	}
\foreach \i in {1,4}{
	\draw[very thick,red] (\i-.5,-1)--(7,6.5-\i);
	}
\foreach \i in {1,5}{
	\draw[very thick] (-1,6.5-\i)--(6.5-\i,-1);
	\draw[very thick] (7,\i+.5)--(.5+\i,7);
	}
\foreach \i in {0,2,4,6}{
	\draw (\i+1,\i-.5)node[2foldcolorswap]{};
	\draw (\i-.5,\i+1)node[2foldcolorswap]{};
	}
\foreach \i in {2,4}{
	\draw (\i-2,\i+2.5)node[2foldcolorswap]{};
	\draw (\i+2.5,\i-2)node[2foldcolorswap]{};
	}
\end{tikzpicture}
}}

\caption{Symmetry diagrams for Lemma~\ref{diagonal_reflection_flip}}
\label{fig:flip_reflection_lemma_results}
\end{figure}

\begin{proof}

Any wallpaper group that includes flip-reflective symmetry needs those lines to have slope $\pm 1$. This will mean the symmetry diagram must exhibit one of the arrangements detailed in Lemma~\ref{diagonal_reflection_flip}. Each wallpaper group included in this lemma has lines of flip-reflective symmetry, and each of their symmetry diagrams have been provided in Figure~\ref{fig:flip_reflection_lemma_results}. For each, the arrangement of flip-reflections with other symmetries are not allowable in PGHPs by Lemma~\ref{diagonal_reflection_flip}.

In particular, $pm/p1$, $pm/pg$, and $pmm/pm$ requires parallel lines of flip-reflection with no parallel glide-reflections. The type $pm/pm(m')$ requires parallel lines of reflection and flip-reflection, but doesn't have a flip-translation symmetry.
\end{proof}

There are two more two-colour wallpaper groups that are incompatible with PGHPs. Each of these require somewhat involved arguments that are similar to the one showing that the one-colour wallpaper type $pgg$ is incompatible with PGHPs in Sen \& Martinez. To be able to make these arguments, we expand upon a lemma from Sen \& Martinez concerning the distance of a glide in various glide-reflections.

\begin{lemma} \label{Glide_Parity}
\begin{enumerate}
\item \citep{Sen_Martinez} A horizontal or vertical line of glide-reflective symmetry lies on (resp.\ off) a gridline if and only if the distance of glide is even (resp.\ odd).
\item A horizontal or vertical line of flip-glide-reflective symmetry lies on (resp.\ off) a gridline if and only if the distance of glide is odd (resp.\ even).
\end{enumerate}
\end{lemma}

\begin{proof}

Part (1) was already proven in Sen \& Martinez. So, we prove (2). We will consider a line of flip-glide-reflective symmetry that is horizontal; the vertical case is symmetric. 

If the line of flip-glide-reflective symmetry lies off a gridline, then this line is a reflection symmetry when considering only the vertical stitches. It follows that the flip-translation associated with the flip-glide-reflection must also be a symmetry of the vertical stitches. 

In order for this flip-glide-reflection to not just be a reflection symmetry, the associated flip-translation cannot be a symmetry for horizontal stitches. Due to the alternating nature of the stitches, a flip-translation of odd length will be a symmetry for the horizontal stitches, so this flip-glide-reflection must have an even-length glide distance.

Now suppose there exists a horizontal line of flip-glide-reflection with an even distance glide. If this line were to lie on a gridline, it passes through a set of horizontal stitches. These horizontal stitches are translated by some distance, swapped to the back, and are fixed by the reflection part of the flip-glide-reflection. Hence, the flip-translation associated with the flip-glide-reflection is a symmetry for this line of horizontal stitches. However, a flip-translation along a horizontal line must have odd-length distance by the alternating nature of the stitches. This implies that the line of flip-glide-reflection must actually lie off the gridlines.
\end{proof}

\begin{figure}[h!tbp]
\centering
\subfloat[Odd length glide-reflection]{%
\resizebox*{6cm}{!}{\includegraphics{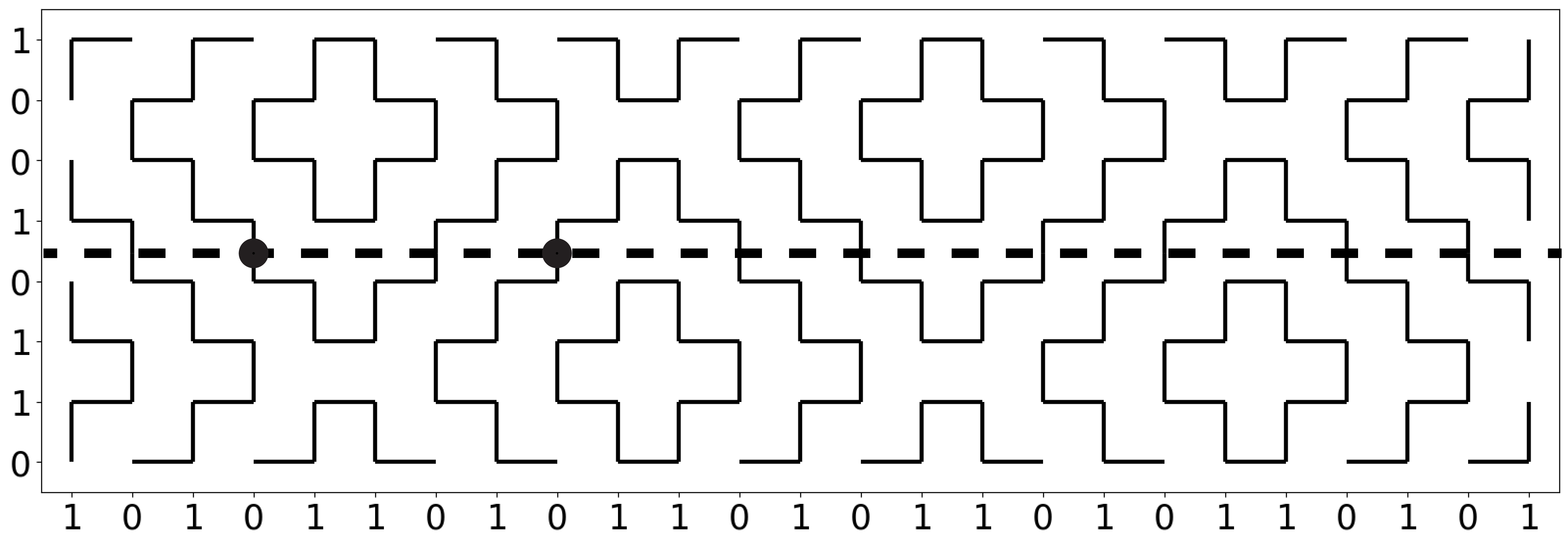}}}\hspace{10pt}
\subfloat[Even length flip-glide-reflection]{%
\resizebox*{6cm}{!}{\includegraphics{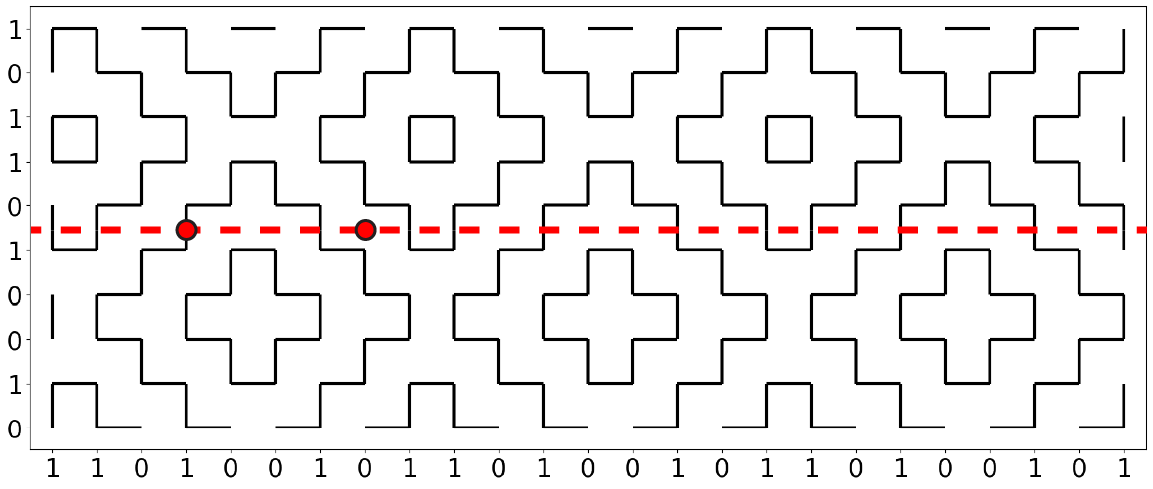}}}

\subfloat[Odd length flip-glide-reflection]{%
\resizebox*{6cm}{!}{\includegraphics{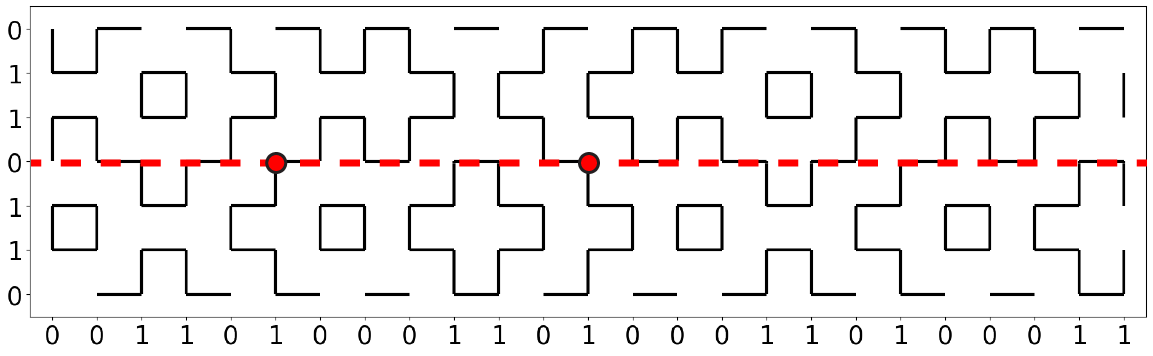}}}\hspace{10pt}
\subfloat[Even length glide-reflection]{%
\resizebox*{6cm}{!}{\includegraphics{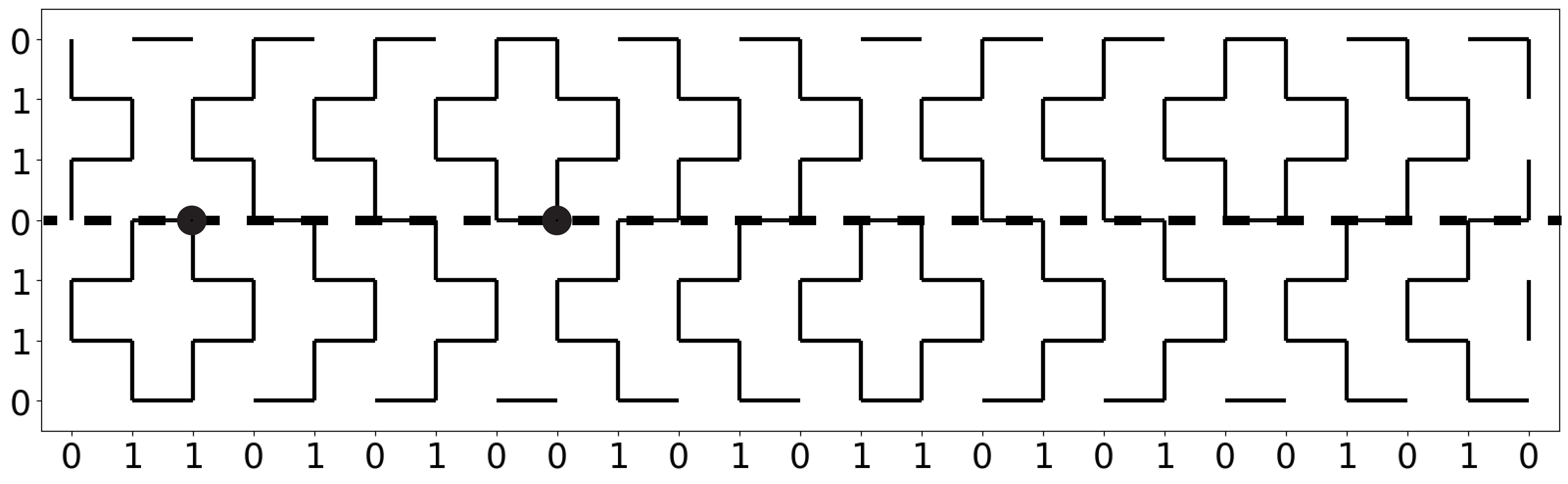}}}

\caption{Examples of glide-reflection and flip-glide-reflection symmetries with different parity glide distances.} \label{fig:Glides}
\end{figure}

Examples illustrating Lemma~\ref{Glide_Parity} are in Figure~\ref{fig:Glides}; each example has dots marking the smallest glide distance. We now have all the tools we need to deal with the two remaining incompatible wallpaper types.

\begin{lemma} \label{lem:pgg/pg}
The two-colour wallpaper groups $pgg/pg$ and $pgg/p2$ are incompatible with PGHPs.
\end{lemma}

\begin{proof}
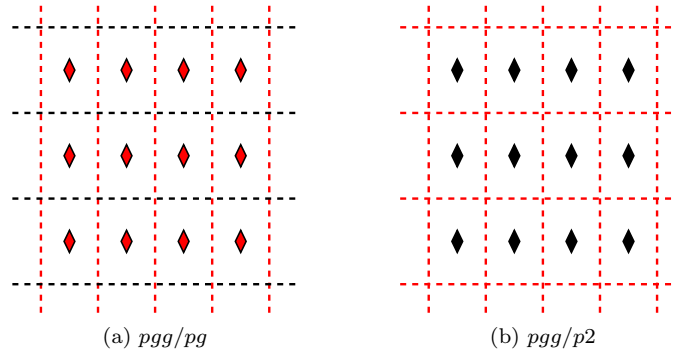
\begin{figure}[h!tbp]
\centering
\subfloat[$pgg/pg$]{%
\resizebox*{4cm}{!}{
\begin{tikzpicture}[scale=.5]
\foreach \i in {1,3,5,7,9}
	\draw[very thick,dashed,red] (\i,0)--(\i,11);
\foreach \j in {1,4,7,10}
	\draw[very thick,dashed] (0,\j)--(10,\j);
\foreach \k in {2,4,6,8}{
	\foreach \l in {2.5, 5.5, 8.5}{
		\node[draw,2foldcolorswap] at (\k,\l){}; 
		}
	}
\end{tikzpicture}
}}
\hspace{1cm}
\subfloat[$pgg/p2$]{%
\resizebox*{4cm}{!}{
\begin{tikzpicture}[scale=.5]
\foreach \i in {1,3,5,7,9}
	\draw[very thick,dashed,red] (\i,0)--(\i,11);
\foreach \j in {1,4,7,10}
	\draw[very thick,dashed,red] (0,\j)--(10,\j);
\foreach \k in {2,4,6,8}{
	\foreach \l in {2.5, 5.5, 8.5}{
		\node[draw,2fold] at (\k,\l){}; 
		}
	}
\end{tikzpicture}
}}
\caption{The symmetry diagram for Lemma~\ref{lem:pgg/pg}} \label{fig:pgg/pg}
\end{figure}

The symmetry diagrams for $pgg/pg$ and $pgg/p2$ are provided in Figure~\ref{fig:pgg/pg}---note that in these diagrams the distances between parallel lines of glide-reflection or flip-glide-reflection are arbitrary and the lines form rectangles. 

For both types, the lines of glide-reflection and flip-glide-reflection must be in the horizontal/vertical orientation, since Lemmas~\ref{diagonal_glide2} and \ref{diagonal_glide3} show there would be parallel lines of reflection or flip-reflection when they have slope $\pm 1$.

It is the case that the glide distances for these wallpaper types correspond to the lengths of the sides of the rectangles formed by the lines of glide-reflective or flip-glide-reflective symmetry. We will need to make a careful argument considering the spacing between the glide-reflection lines to show that these symmetry diagrams are not possible for a PGHP. 

First consider type $pgg/pg$. For this proof, we are going to assume the lines of flip-glide-reflection are vertical and the glide-reflections are horizontal. Ostensibly, the symmetry diagram in Figure~\ref{fig:pgg/pg}(a) could be rotated by 90 degrees, but it is a symmetric argument to show that diagram is also incompatible.

The flip-two-fold rotation centers must be either at the intersection of gridlines or centered in a square of the grid. This follows from Table~\ref{rosettes}, where we see that $\{e,\overline{r_{180}}\}$ is only possible when $|x|$ and $|y|$ have the same parity. 

Suppose the flip-two fold rotation is centered in a square of the grid.  If the horizontal glide-reflection has even glide distance, this forces the vertical flip-glide-reflection lines to be off-grid, as shown in Figure~\ref{fig:pgg/pg_spacing}(a). Additionally, by Lemma~\ref{Glide_Parity}, we know the horizontal glide-reflection lines must be on the gridlines. However, this implies that the vertical flip-glide-reflection has odd-length glide. This is a contradiction, since by Lemma~\ref{Glide_Parity}, odd-length glide would imply the flip-glide-reflection line is on a gridline. The case where the horizontal glide-reflection has odd glide distance yields a similar contradiction, as illustrated in Figure~\ref{fig:pgg/pg_spacing}(b).

\begin{figure}[h!tbp] 
\centering
\subfloat[]{%
\resizebox*{3.4cm}{!}{
\begin{tikzpicture}[scale=0.35]
\draw[dotted, thick] (0,0)grid(7,7);
\draw[ultra thick,red,dashed] (1.5,-.25)--(1.5,7.25);
\draw[ultra thick,red,dashed] (5.5,-.25)--(5.5,7.25);
\node[2foldcolorswap] (d) at (3.5,3.5) {};
\draw [decorate,decoration={brace,amplitude=10pt,mirror},xshift=0pt,yshift=0pt]
(1.5,-.25) -- (5.5,-.25) node [black,midway,xshift=0cm,yshift=-5mm] 
{even};
\draw[ultra thick,black,dashed] (-.25,1)--(7.25,1);
\draw[ultra thick,black,dashed] (-.25,6)--(7.25,6);
\draw [decorate,decoration={brace,amplitude=10pt},xshift=0pt,yshift=0pt]
(-.25,1) -- (-.25,6) node [black,midway,xshift=-9mm,yshift=0mm] 
{odd};
\end{tikzpicture}
}}\hspace{2pt}
\subfloat[]{%
\resizebox*{3.4cm}{!}{
\begin{tikzpicture}[scale=0.35]
\draw[dotted, thick] (0,0)grid(7,7);
\draw[ultra thick,red,dashed] (1,-.25)--(1,7.25);
\draw[ultra thick,red,dashed] (6,-.25)--(6,7.25);
\node[2foldcolorswap] (d) at (3.5,3.5) {};
\draw [decorate,decoration={brace,amplitude=10pt,mirror},xshift=0pt,yshift=0pt]
(1,-.25) -- (6,-.25) node [black,midway,xshift=0cm,yshift=-5mm] 
{odd};
\draw[ultra thick,black,dashed] (-.25,1.5)--(7.25,1.5);
\draw[ultra thick,black,dashed] (-.25,5.5)--(7.25,5.5);
\draw [decorate,decoration={brace,amplitude=10pt},xshift=0pt,yshift=0pt]
(-.25,1.5) -- (-.25,5.5) node [black,midway,xshift=-9mm,yshift=0mm] 
{even};
\end{tikzpicture}
}}
\hspace{2pt}
\subfloat[]{%
\resizebox*{3.4cm}{!}{
\begin{tikzpicture}[scale=0.35]
\draw[dotted, thick] (0,0)grid(8,8);
\draw[ultra thick,red,dashed] (2,-.25)--(2,8.25);
\draw[ultra thick,red,dashed] (6,-.25)--(6,8.25);
\node[2foldcolorswap] (d) at (4,4) {};
\draw [decorate,decoration={brace,amplitude=10pt,mirror},xshift=0pt,yshift=0pt]
(2,-.25) -- (6,-.25) node [black,midway,xshift=0cm,yshift=-5mm] 
{even};
\draw[ultra thick,black,dashed] (-.25,1)--(8.25,1);
\draw[ultra thick,black,dashed] (-.25,7)--(8.25,7);
\draw [decorate,decoration={brace,amplitude=10pt},xshift=0pt,yshift=0pt]
(-.25,1) -- (-.25,7) node [black,midway,xshift=-9mm,yshift=0mm] 
{even};
\end{tikzpicture}
}}\hspace{2pt}
\subfloat[]{%
\resizebox*{3.4cm}{!}{
\begin{tikzpicture}[scale=0.35]
\draw[dotted, thick] (0,0)grid(8,8);
\draw[ultra thick,red,dashed] (1.5,-.25)--(1.5,8.25);
\draw[ultra thick,red,dashed] (6.5,-.25)--(6.5,8.25);
\node[2foldcolorswap] (d) at (4,4) {};
\draw [decorate,decoration={brace,amplitude=10pt,mirror},xshift=0pt,yshift=0pt]
(1.5,-.25) -- (6.5,-.25) node [black,midway,xshift=0cm,yshift=-5mm] 
{odd};
\draw[ultra thick,black,dashed] (-.25,1.5)--(8.25,1.5);
\draw[ultra thick,black,dashed] (-.25,6.5)--(8.25,6.5);
\draw [decorate,decoration={brace,amplitude=10pt},xshift=0pt,yshift=0pt]
(-.25,1.5) -- (-.25,6.5) node [black,midway,xshift=-9mm,yshift=0mm] 
{odd};
\end{tikzpicture}
}}\hspace{3pt}

\caption{An illustration of how lines of glide-reflection would relate to a flip-two-fold center of symmetry in $pgg/pg$.}
\label{fig:pgg/pg_spacing}
\end{figure}
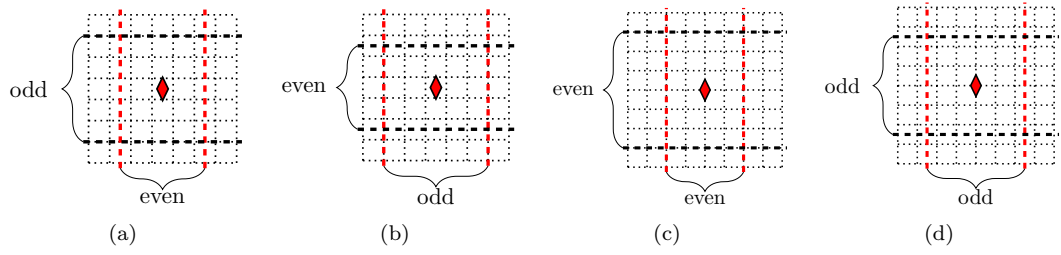

When the center of 2-fold flip-rotation is at the intersection of gridlines, we can obtain similar contradictions. When we assume the horizontal glide-reflection has even distance, we get the diagram in Figure~\ref{fig:pgg/pg_spacing}(c), where the flip-glide-reflection is forced to be on a gridline but have even glide distance. In the case of odd glide distance, we get the diagram in Figure~\ref{fig:pgg/pg_spacing}(d), where the flip-glide-reflection is forced to be off a gridline but have odd glide distance. 

It follows that $pgg/pg$ is incompatible with PGHPs.

The proof for $pgg/p2$ proceeds similarly. Since $pgg/p2$ has centers of two-fold rotational symmetry with no lines of reflection, the centers cannot be off the gridlines or at the intersection of gridlines by Lemma~\ref{cor:4Fold}. Assume the center of rotation lies on a horizontal gridline as shown in Figure~\ref{fig:pgg/p2_spacing}; the case of the center being on a vertical gridline is symmetric.

Now suppose the horizontal line of flip-glide-reflection has even-length glide. By Lemma~\ref{Glide_Parity}, this line of flip-glide-reflection must lie off gridlines. This forces the glide-distance for the vertical flip-glide-reflection lines to be odd-length. However, those lines must also lie off gridlines (see Figure~\ref{fig:pgg/p2_spacing}(a)), so we have obtained a contradiction. A similar contradiction can be obtained with odd-length glide distance for the horizontal flip-glide-reflection; a supporting figure is provided in Figure~\ref{fig:pgg/p2_spacing}(b).

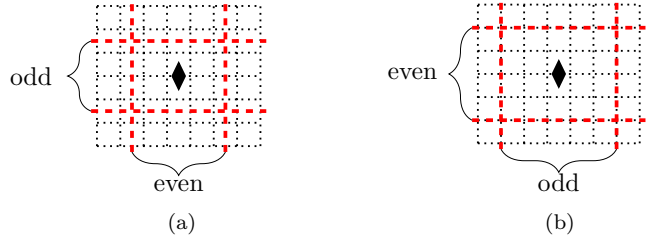
\begin{figure}[h!tbp] 
\centering
\subfloat[]{%
\resizebox*{5cm}{!}{
\begin{tikzpicture}[scale=0.35]
\draw[dotted, thick] (0,0)grid(7,6);
\draw[ultra thick,red,dashed] (1.5,-.25)--(1.5,6.25);
\draw[ultra thick,red,dashed] (5.5,-.25)--(5.5,6.25);
\node[2fold] (d) at (3.5,3) {};
\draw [decorate,decoration={brace,amplitude=10pt,mirror},xshift=0pt,yshift=0pt]
(1.5,-.25) -- (5.5,-.25) node [black,midway,xshift=0cm,yshift=-5mm] 
{even};
\draw[ultra thick,red,dashed] (-.25,1.5)--(7.25,1.5);
\draw[ultra thick,red,dashed] (-.25,4.5)--(7.25,4.5);
\draw [decorate,decoration={brace,amplitude=10pt},xshift=0pt,yshift=0pt]
(-.25,1.5) -- (-.25,4.5) node [black,midway,xshift=-9mm,yshift=0mm] 
{odd};
\node[] (d) at (11,0) {};
\end{tikzpicture}
}}
\subfloat[]{%
\resizebox*{5cm}{!}{
\begin{tikzpicture}[scale=0.35]
\draw[dotted, thick] (0,0)grid(7,6);
\draw[ultra thick,red,dashed] (1,-.25)--(1,6.25);
\draw[ultra thick,red,dashed] (6,-.25)--(6,6.25);
\node[2fold] (d) at (3.5,3) {};
\draw [decorate,decoration={brace,amplitude=10pt,mirror},xshift=0pt,yshift=0pt]
(1,-.25) -- (6,-.25) node [black,midway,xshift=0cm,yshift=-5mm] 
{odd};
\draw[ultra thick,red,dashed] (-.25,1)--(7.25,1);
\draw[ultra thick,red,dashed] (-.25,5)--(7.25,5);
\draw [decorate,decoration={brace,amplitude=10pt},xshift=0pt,yshift=0pt]
(-.25,1) -- (-.25,5) node [black,midway,xshift=-9mm,yshift=0mm] 
{even};
\node[] (d) at (11,0) {};
\end{tikzpicture}
}}

\caption{An illustration of how lines of flip-glide-reflection would relate to a two-fold center of symmetry in $pgg/p2$.}
\label{fig:pgg/p2_spacing}
\end{figure}

It follows that $pgg/p2$ is incompatible with PGHPs.

\end{proof}

Of the 46 two-colour wallpaper groups, we are left with 17 types compatible with PGHPs. The final step is to produce PGHPs that exhibit the symmetries in each of these wallpaper groups.

\begin{theorem}
There are 17 two-colour wallpaper groups that are compatible with PGHPs. These types are: $p1/p1$, $p2/p1$, $p2/p2$, $pm/pm(m)$, $pm/cm$, $pg/p1$, $pg/pg$, $pmg/pm$, $pmg/pmg$, $cm/pm$, $cm/pg$, $cm/p1$, $cmm/pmm$, $cmm/pmg$, $cmm/cm$, $p4m/pmm$, and $p4m/p4m$.
\end{theorem}

\begin{figure}[h!tbp] 
\centering

\includegraphics[width=4.5cm]{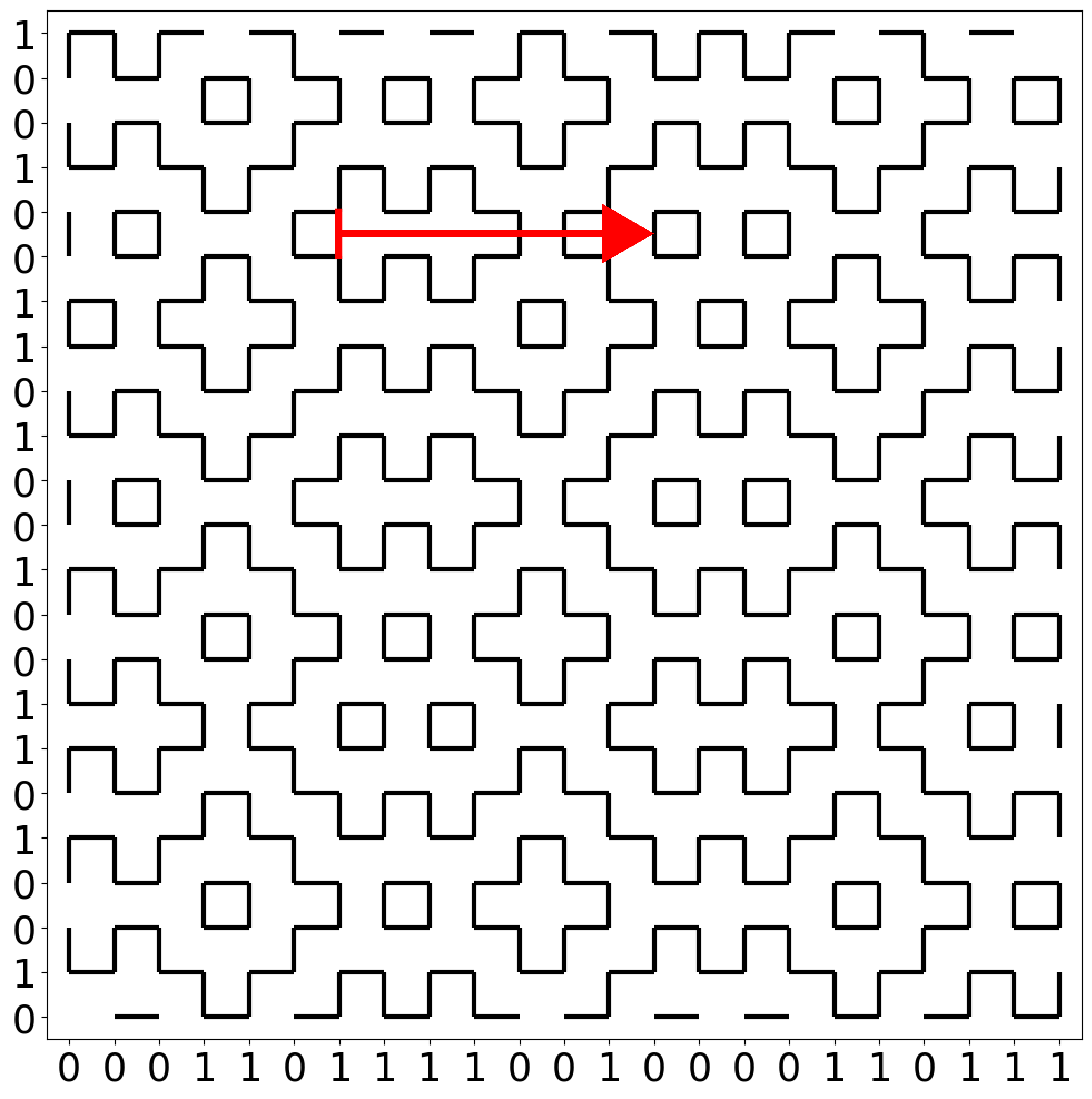}
\includegraphics[width=4.5cm]{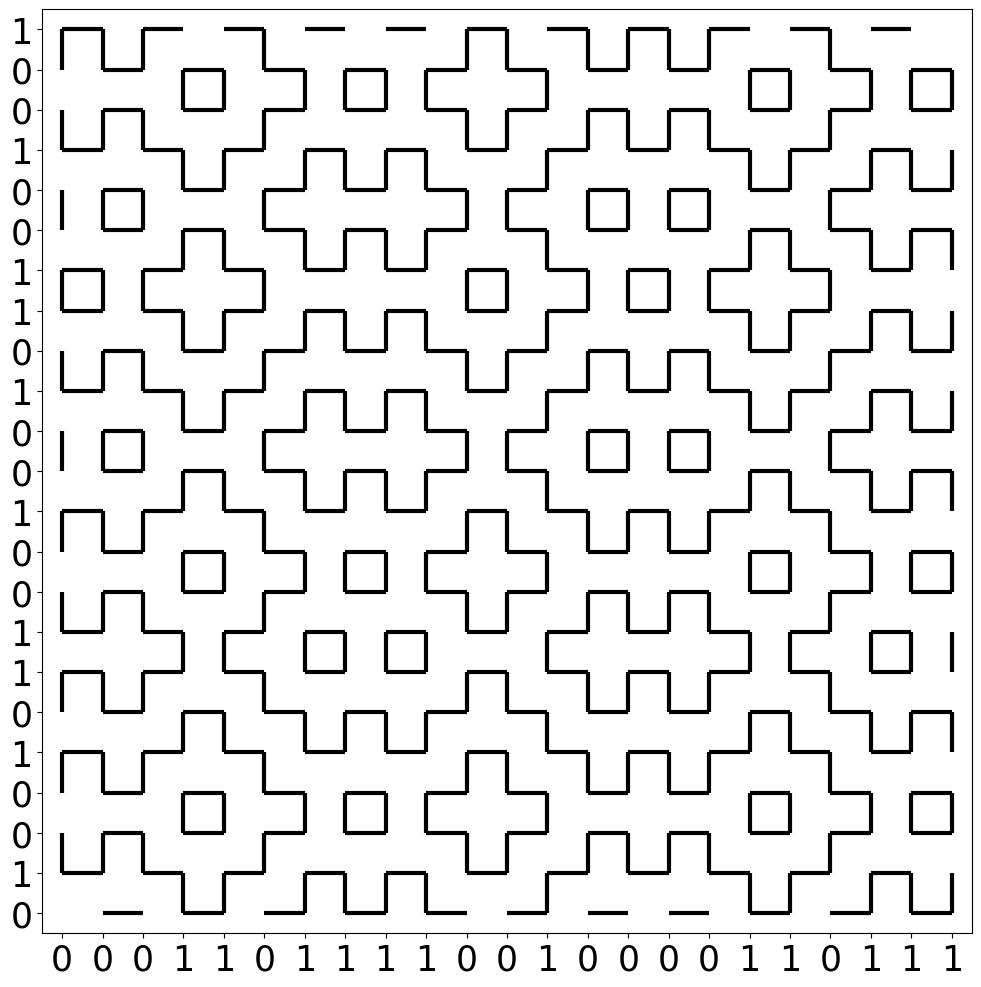}
\includegraphics[width=4.5cm]{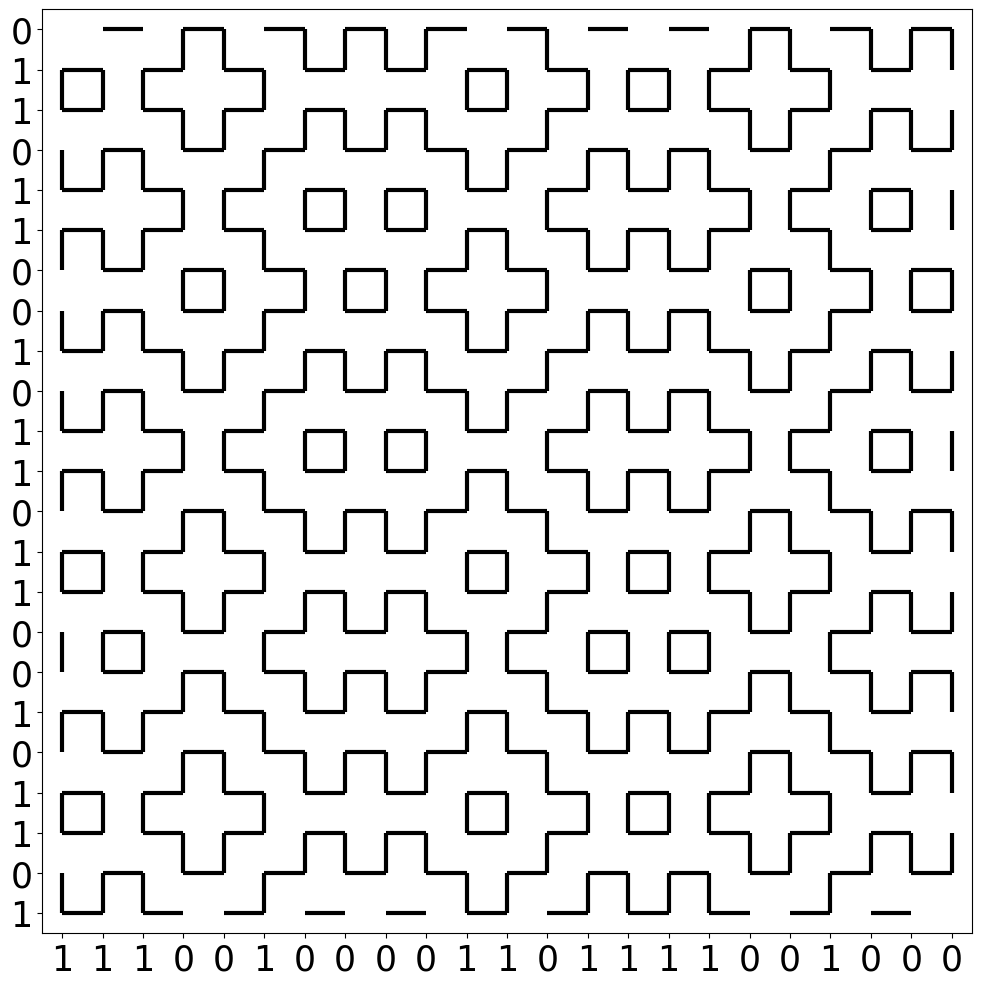}

\caption{Type $p1/p1$ where $x=(00011011110010)^{\infty}$ and $y=(010010110)^{\infty}$.}
\label{fig:p1/p1}
\end{figure}

\begin{figure}[h!tbp] 
\centering

\includegraphics[width=4.5cm]{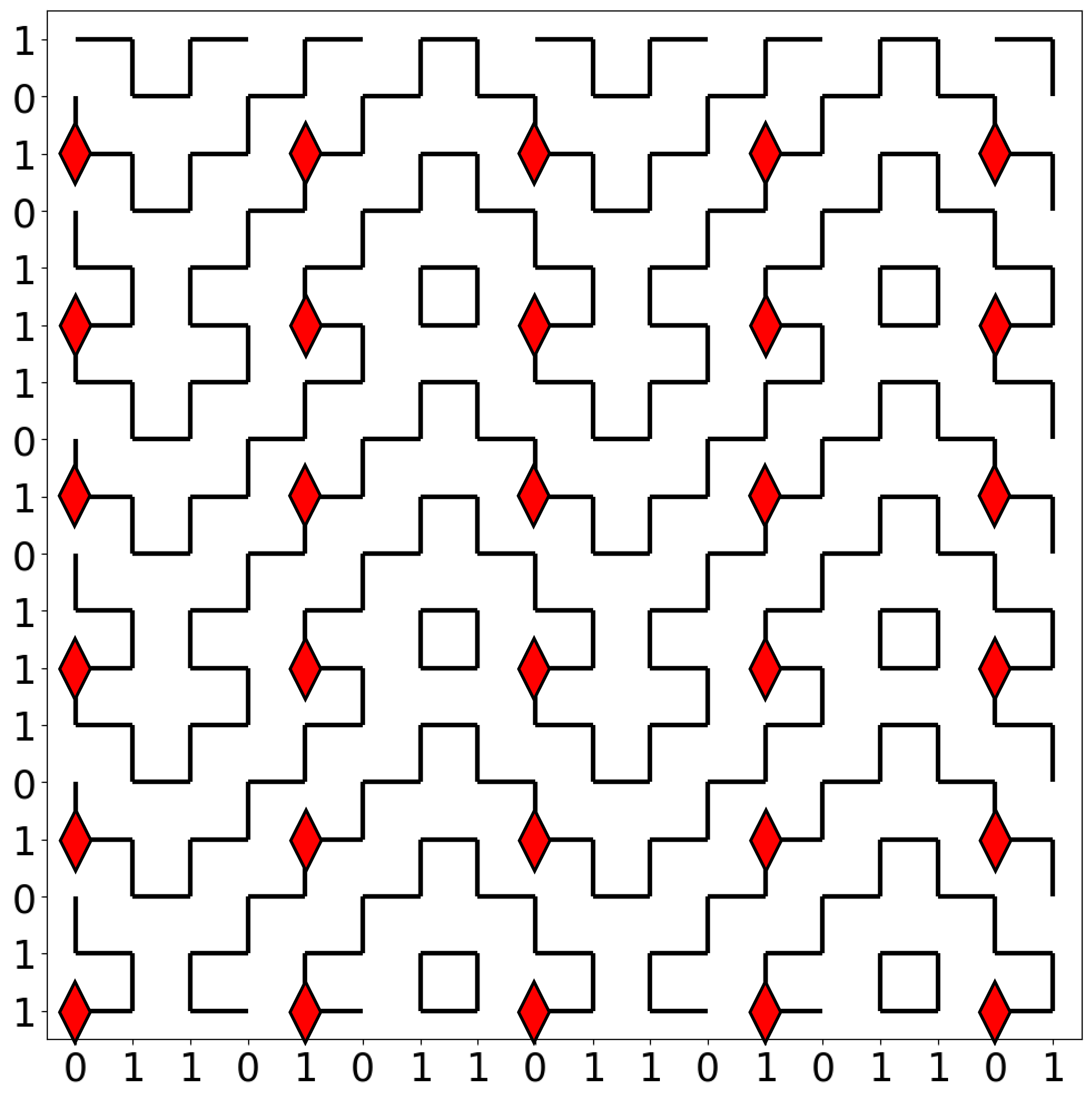}
\includegraphics[width=4.5cm]{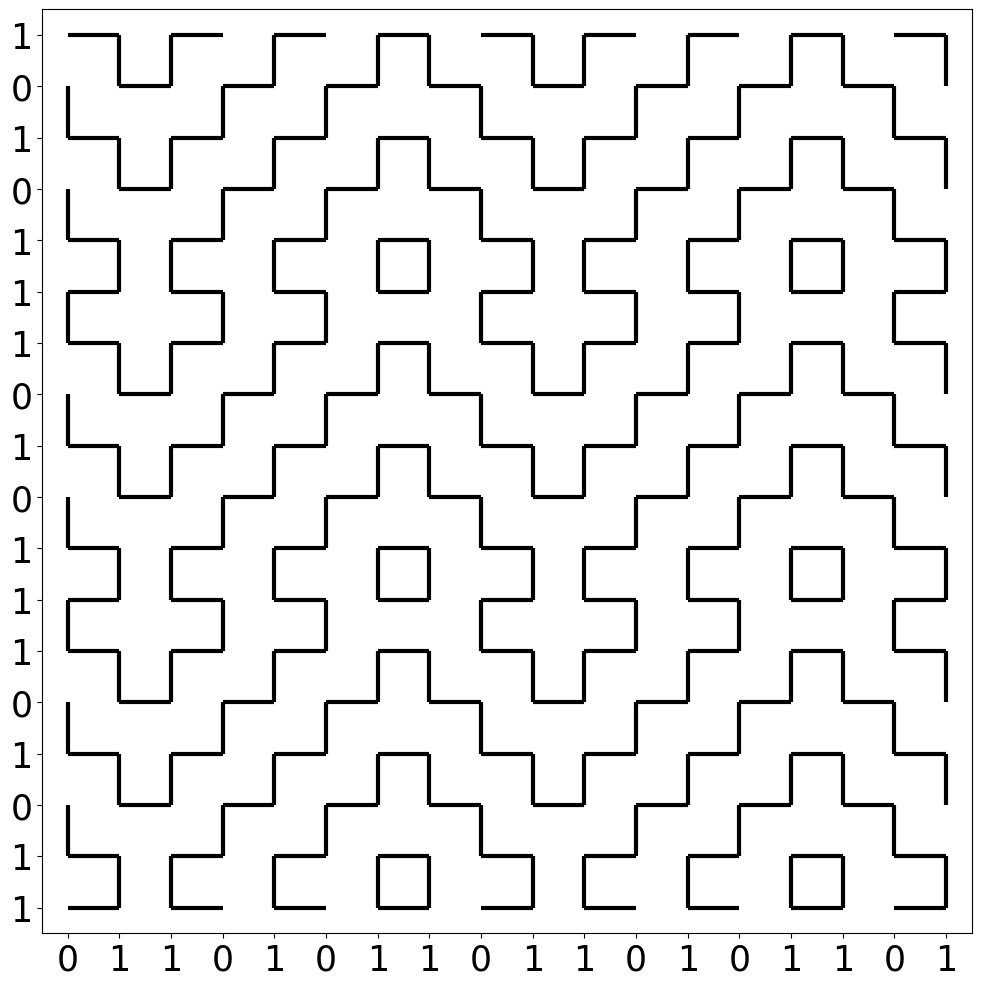}
\includegraphics[width=4.5cm]{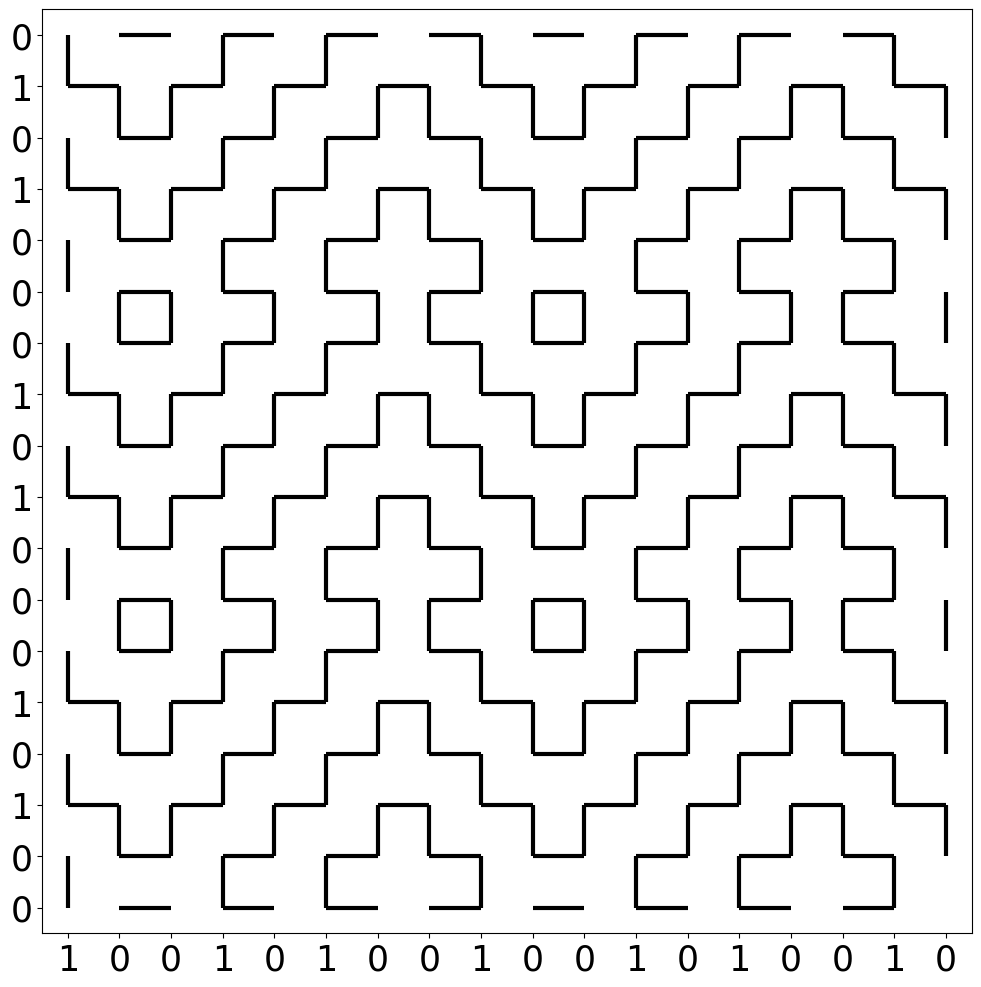}

\caption{Type $p2/p1$ where $x=(01101011)^{\infty}$ and $y=(110101)^{\infty}$}
\label{fig:p2/p1}
\end{figure}

\begin{figure}[h!tbp] 
\centering

\includegraphics[width=4.5cm]{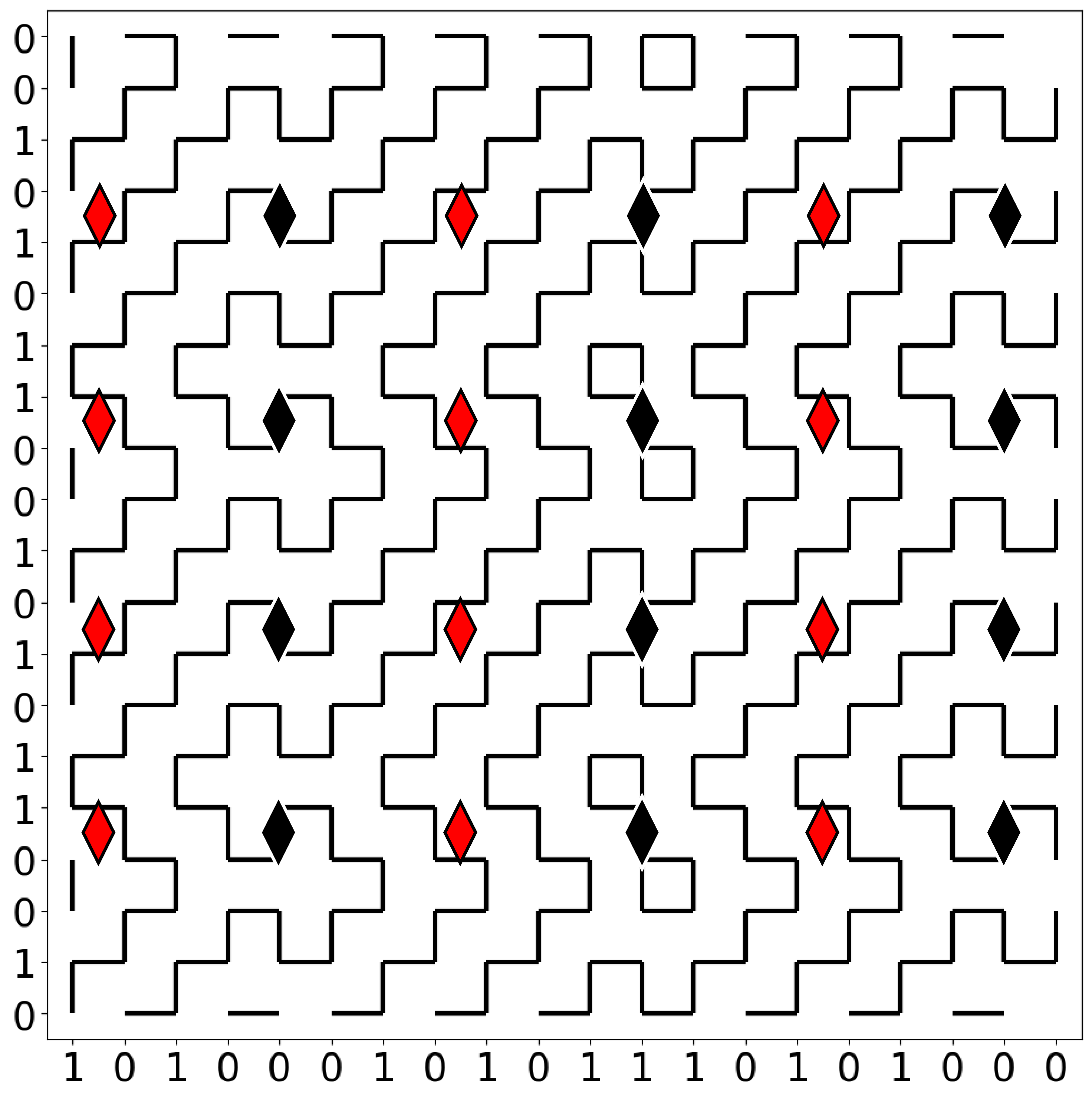}
\includegraphics[width=4.5cm]{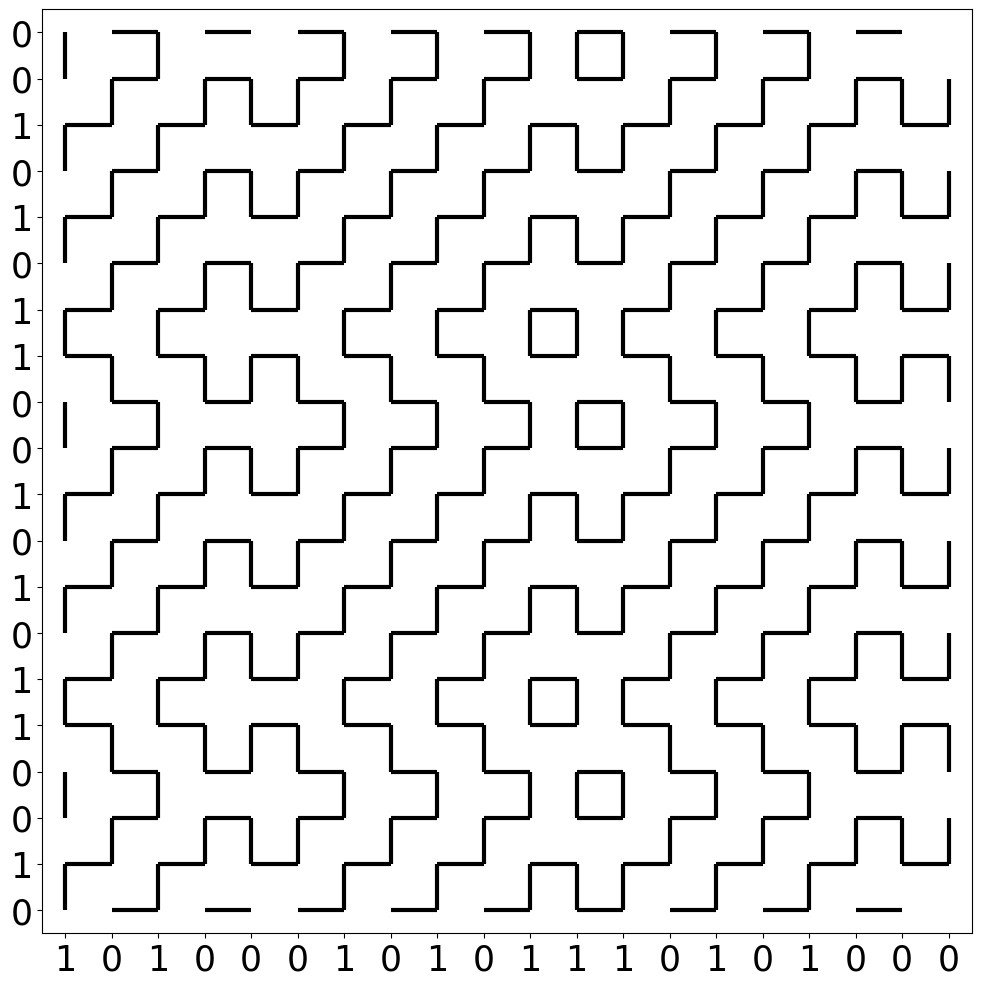}
\includegraphics[width=4.5cm]{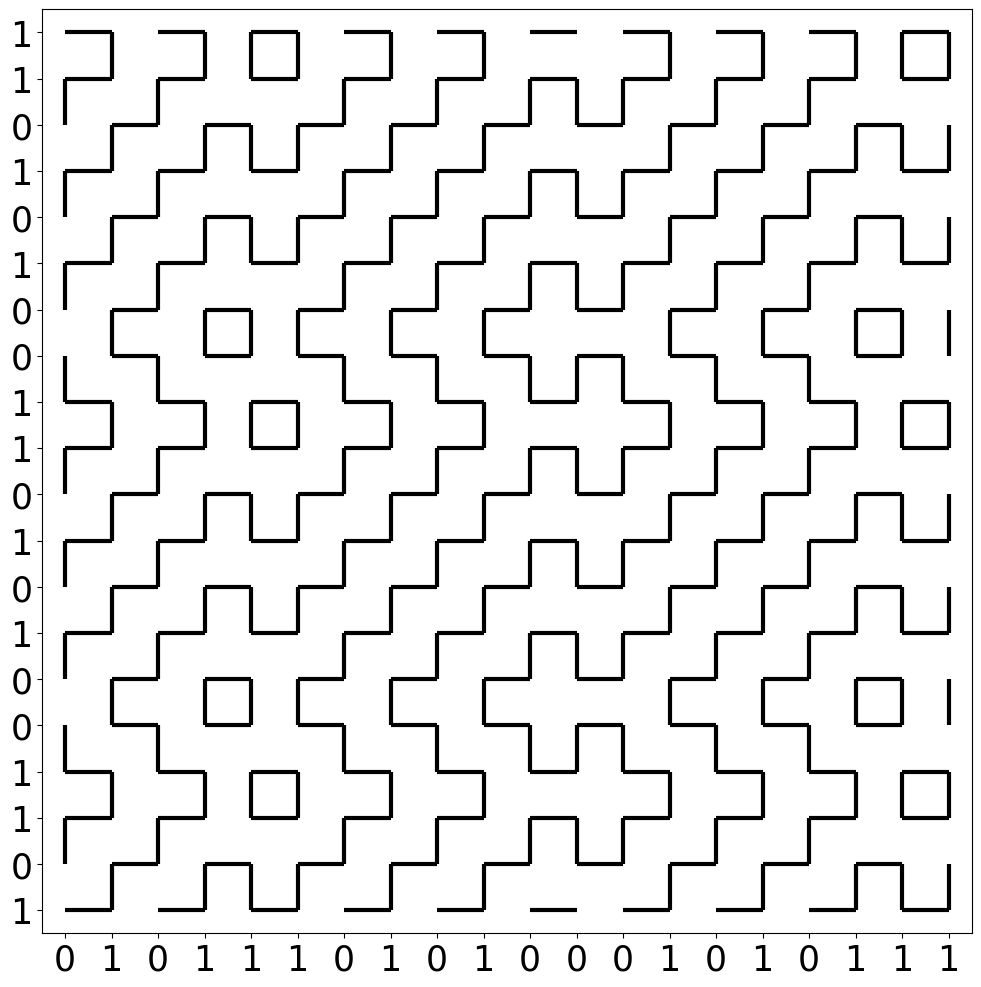}

\caption{Type $p2/p2$ where $x=(10100010101110)^{\infty}$ and $y=(01001101)^{\infty}$}
\label{fig:p2/p2}
\end{figure}

\begin{figure}[h!tbp] 
\centering

\includegraphics[width=4.5cm]{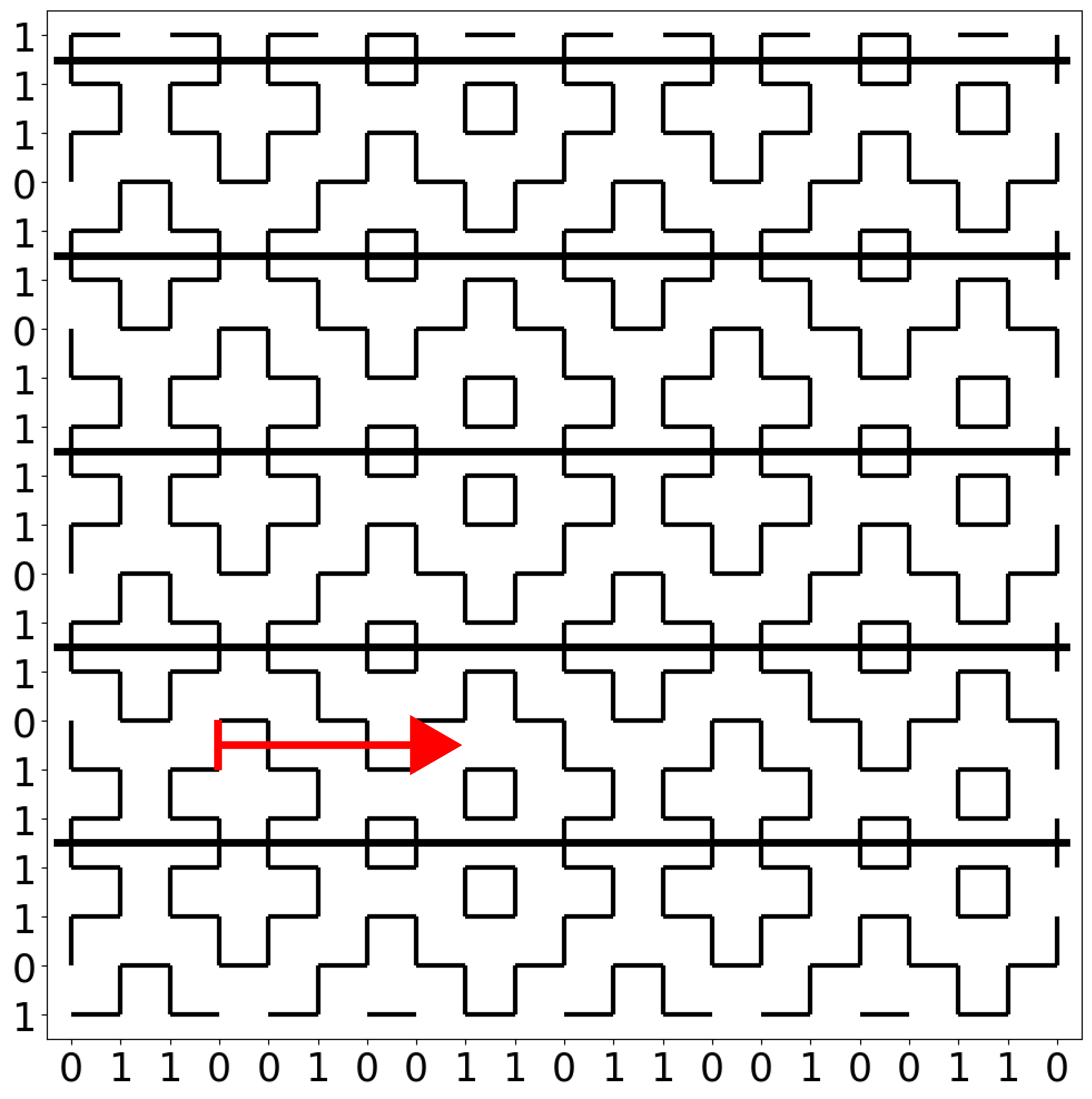}
\includegraphics[width=4.5cm]{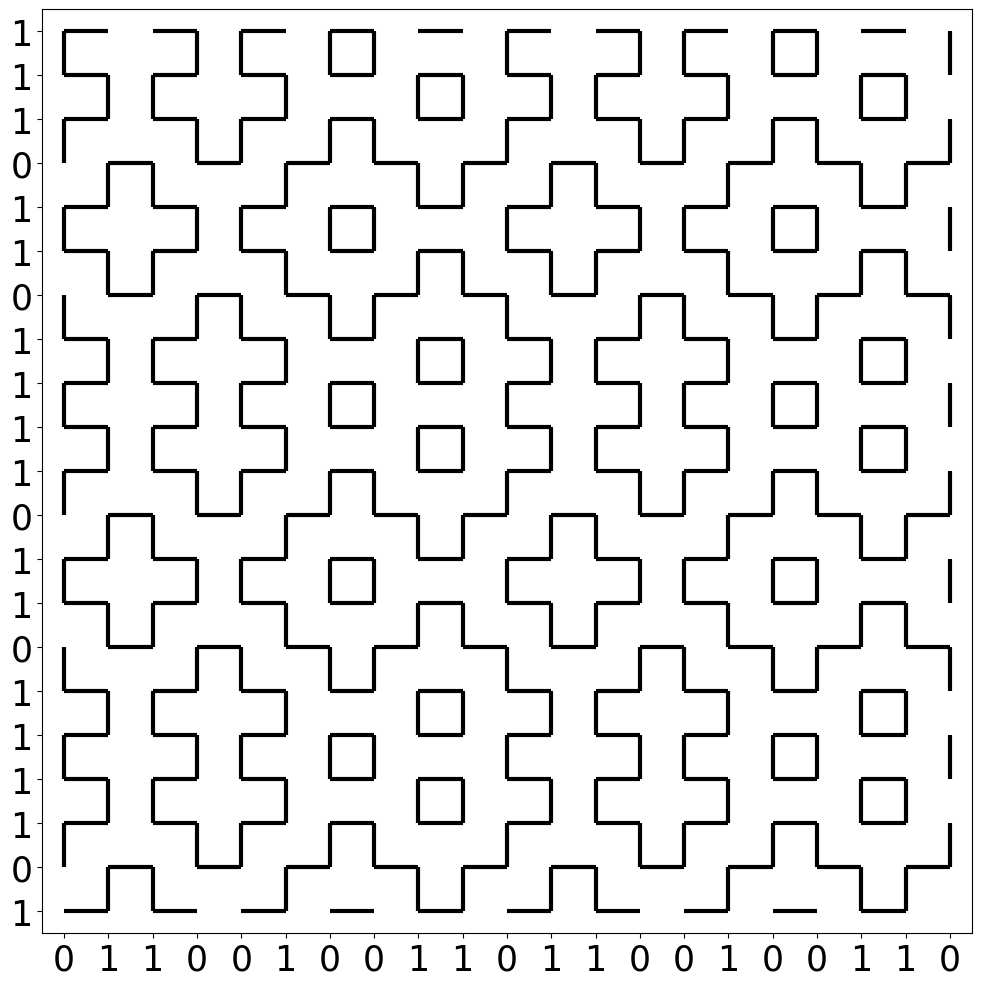}
\includegraphics[width=4.5cm]{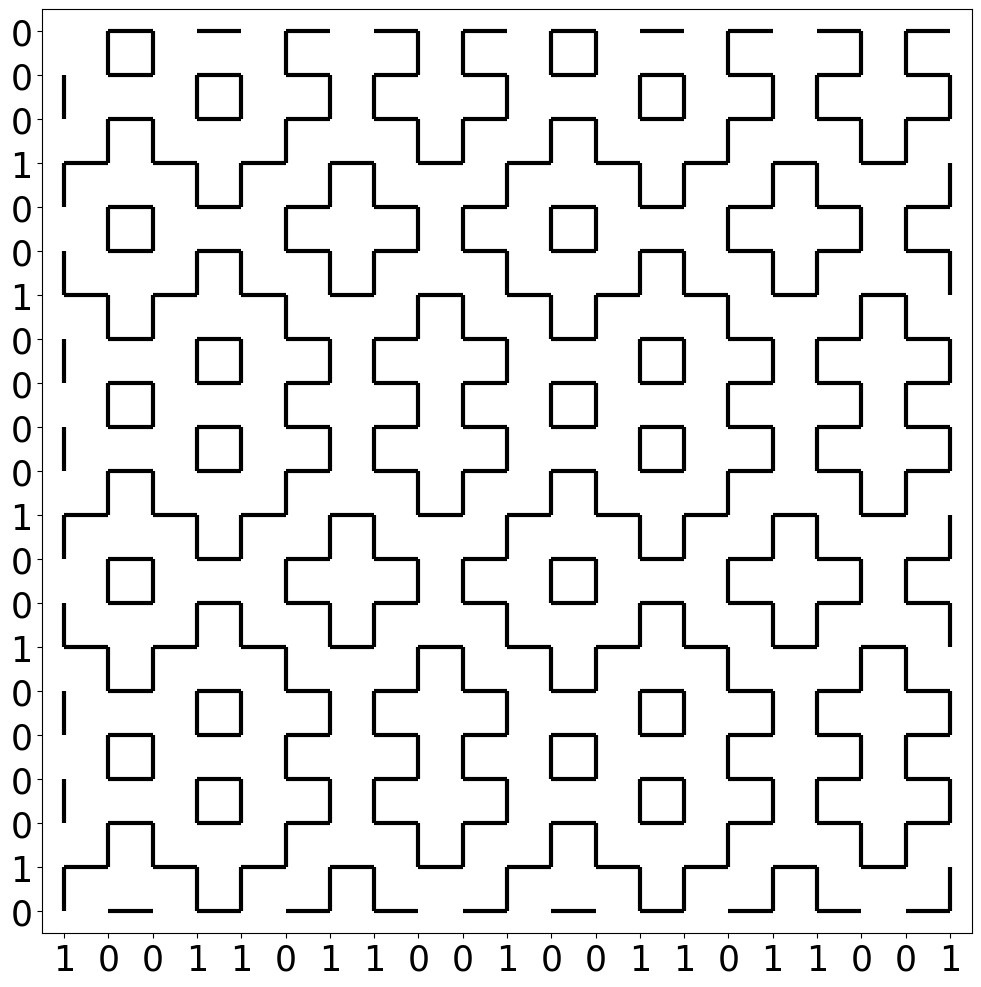}

\caption{Type $pm/pm(m)$ where $x=(0110010011)^{\infty}$ and $y=(10111101)^{\infty}$.}
\label{fig:pm/pm(m)}
\end{figure}

\begin{figure}[h!tbp] 
\centering

\includegraphics[width=4.5cm]{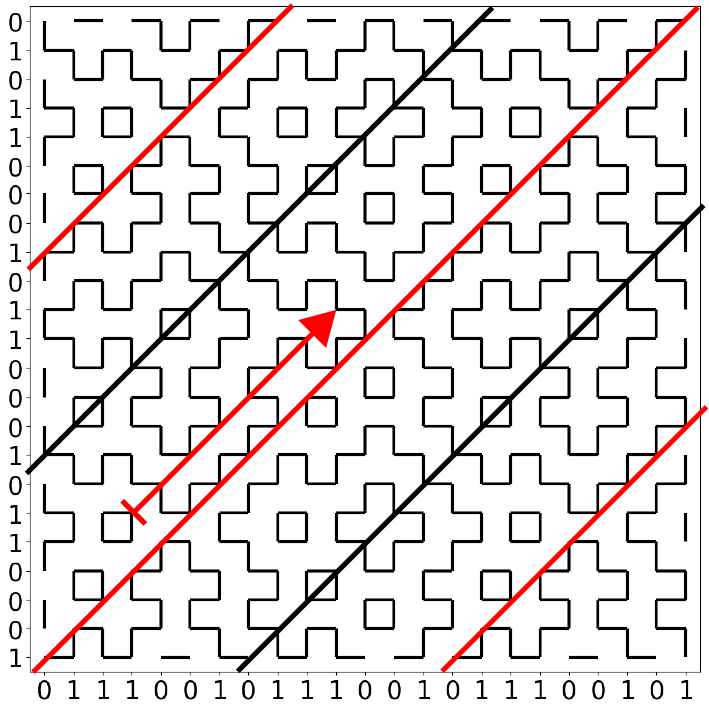}
\includegraphics[width=4.5cm]{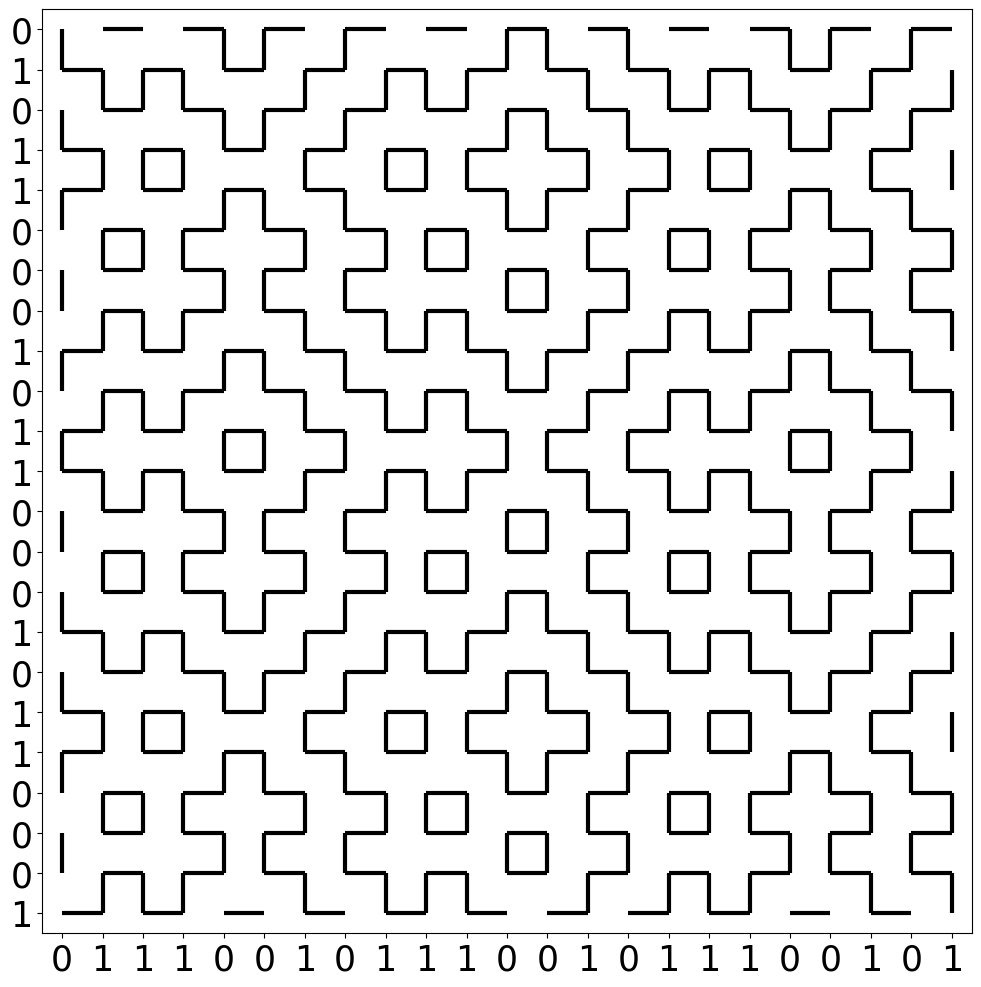}
\includegraphics[width=4.5cm]{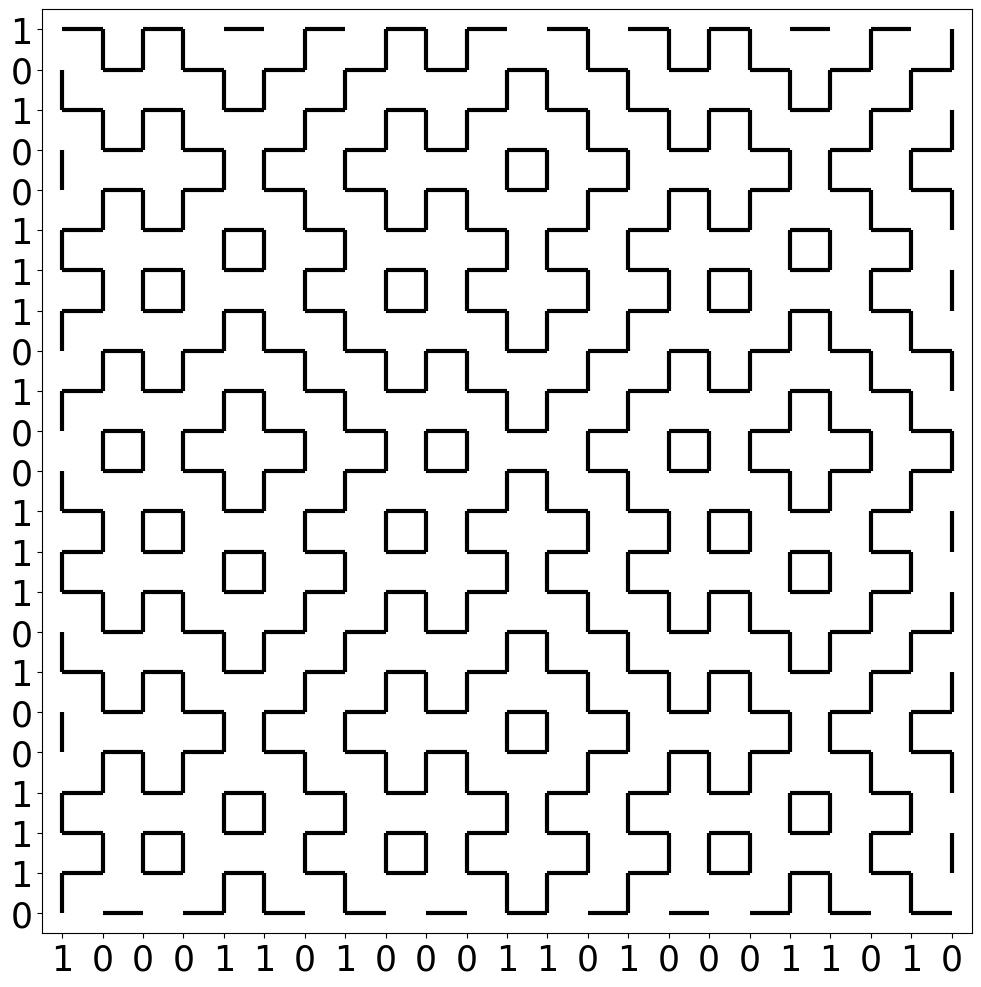}

\caption{Type $pm/cm$ where $x=(0111001)^{\infty}$ and $y=(1000110)^{\infty}$.}
\label{fig:pm/cm}
\end{figure}

\begin{figure}[h!tbp] 
\centering

\includegraphics[width=4.5cm]{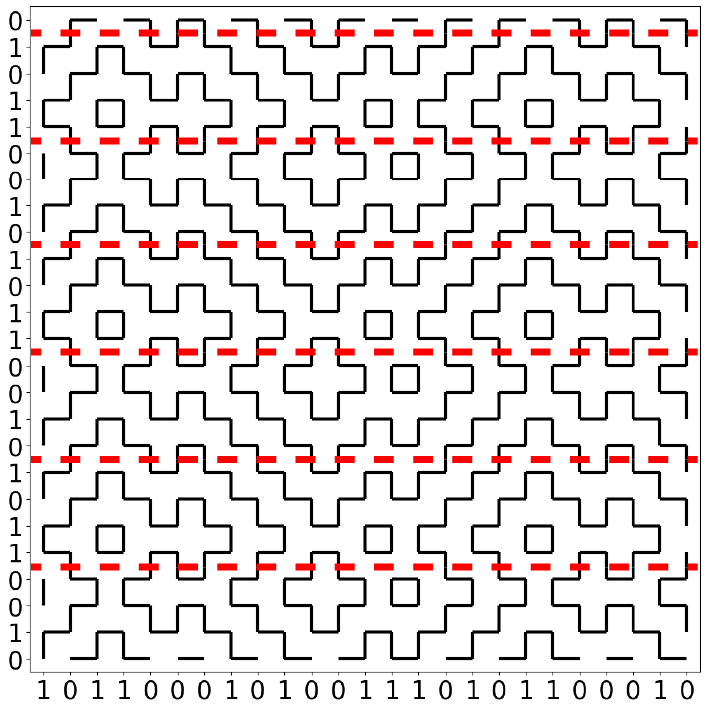}
\includegraphics[width=4.5cm]{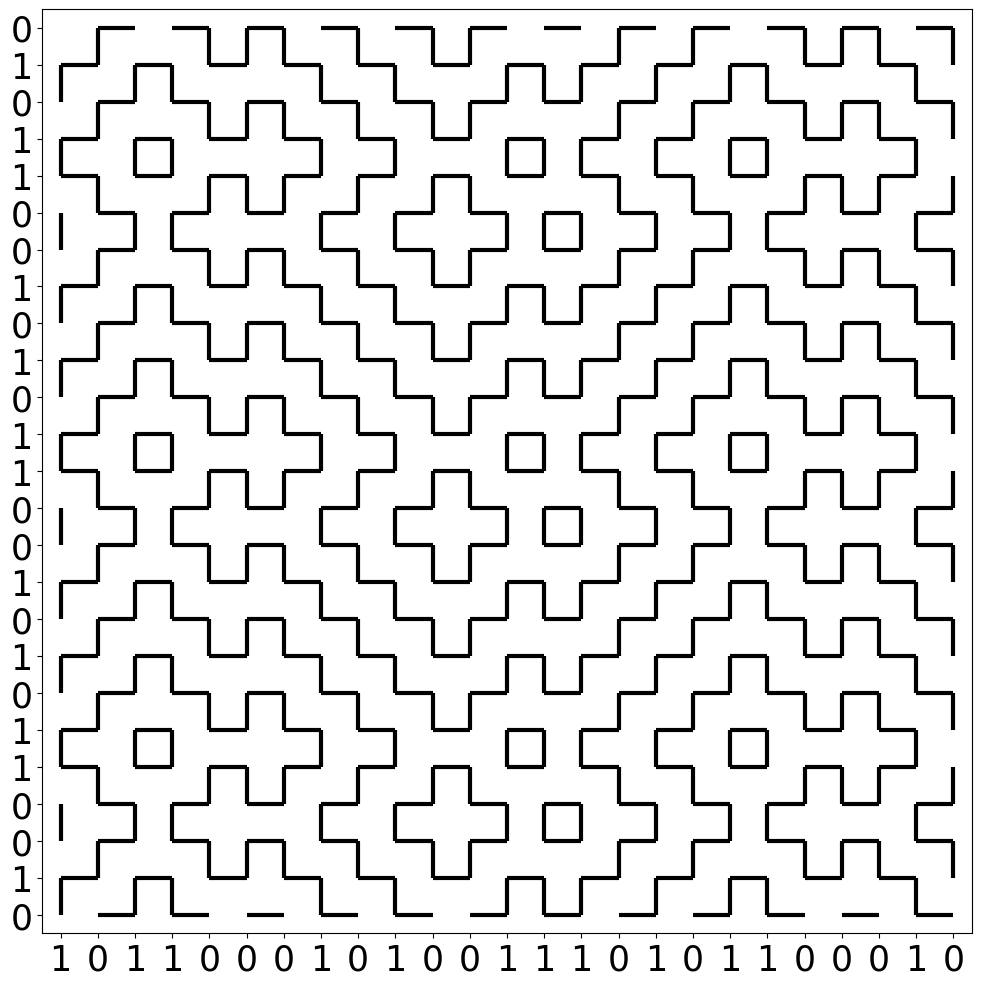}
\includegraphics[width=4.5cm]{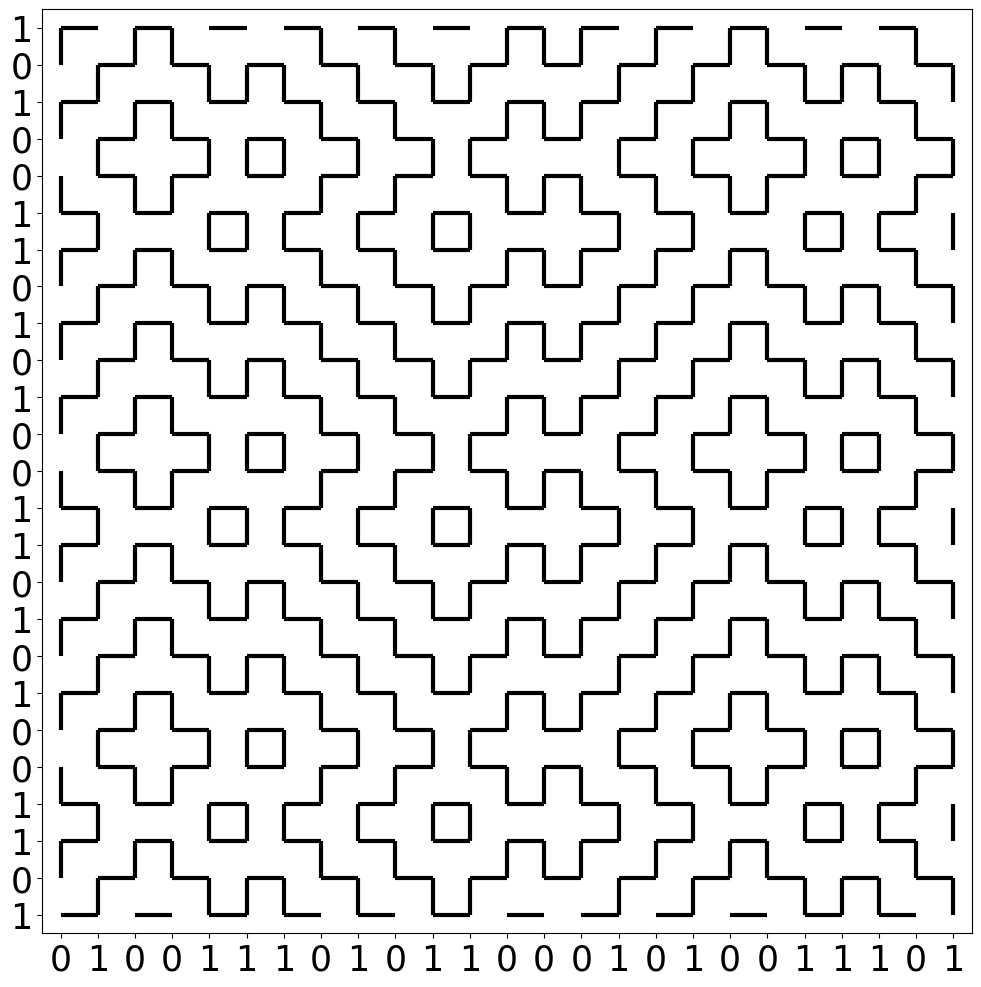}

\caption{Type $pg/p1$ where $x=(1011000101001110)^{\infty}$ and $y=(01001101)^{\infty}$.}
\label{fig:pg/p1}
\end{figure}

\begin{figure}[h!tbp] 
\centering

\includegraphics[width=4.5cm]{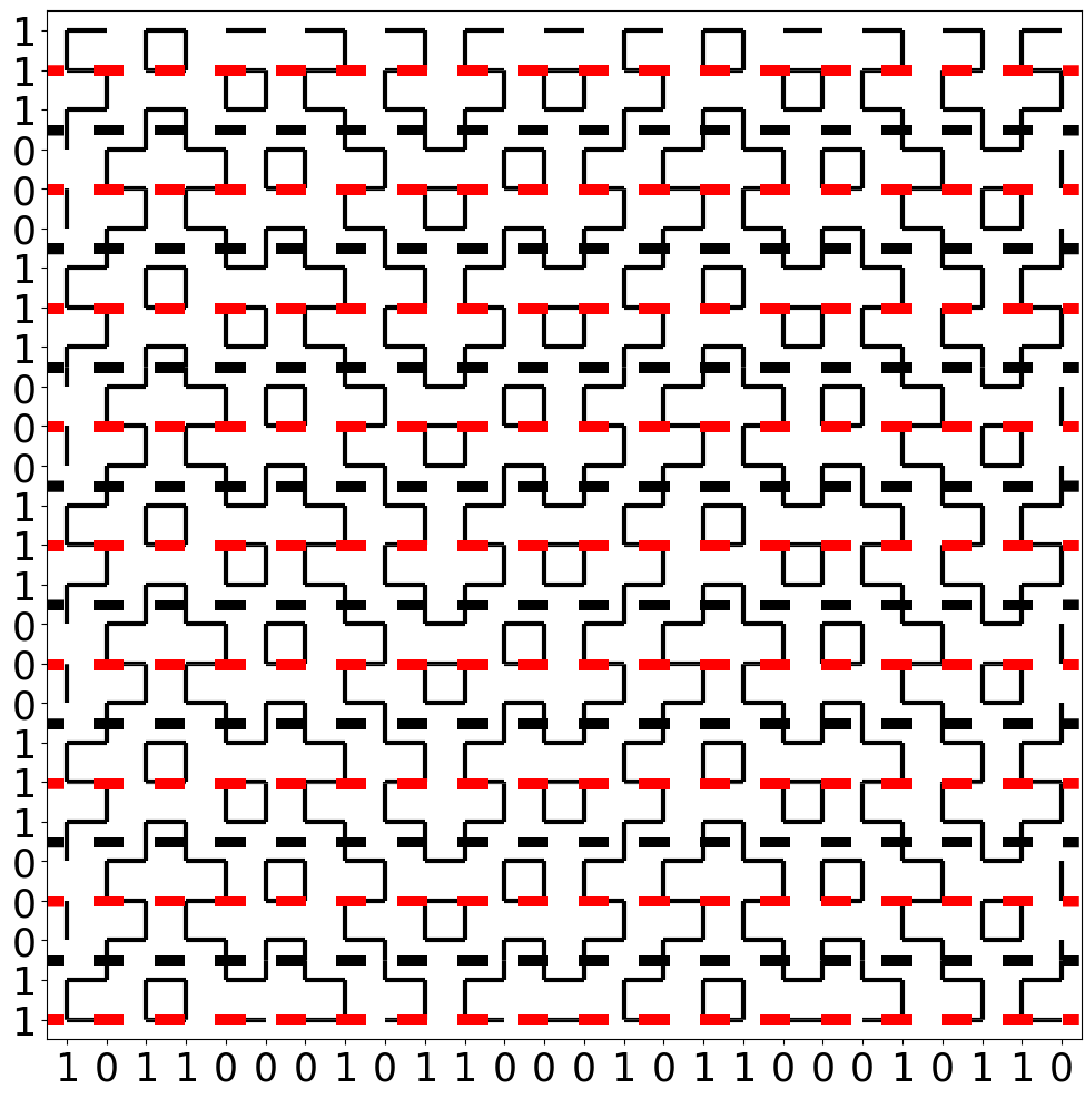}
\includegraphics[width=4.5cm]{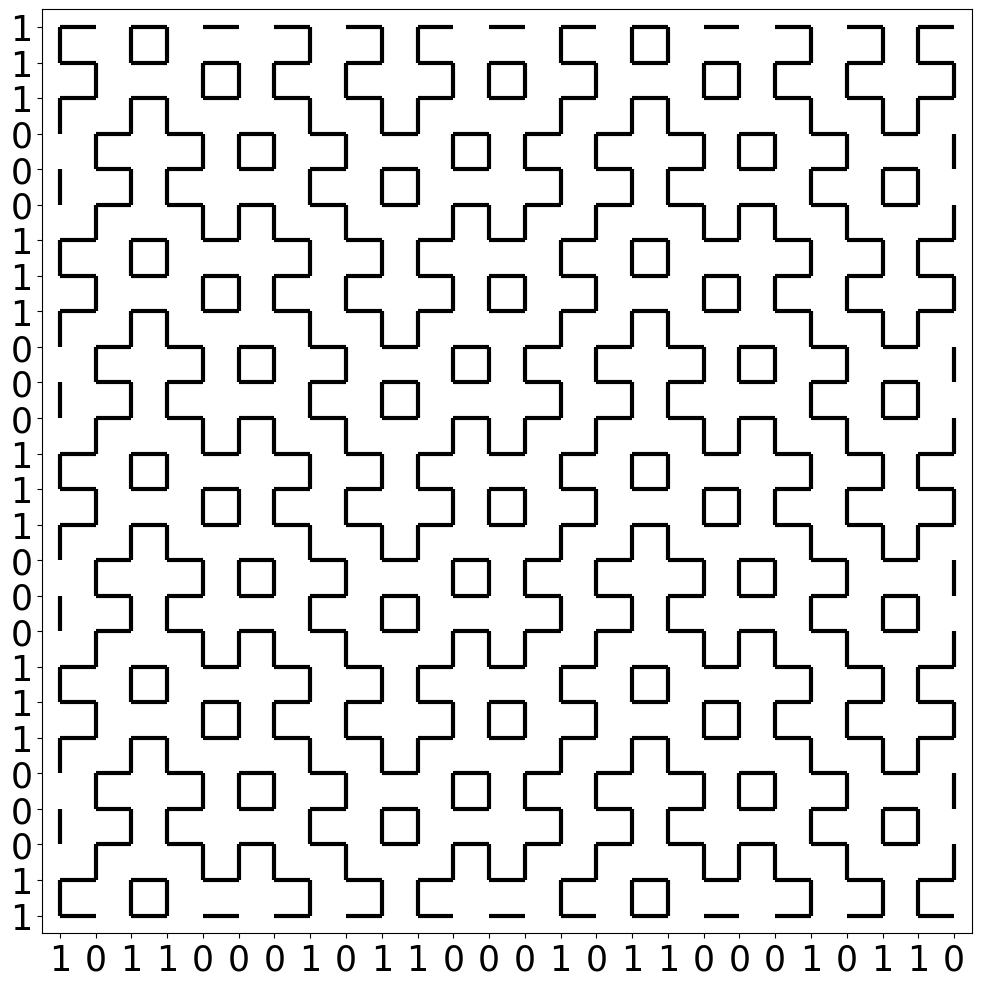}
\includegraphics[width=4.5cm]{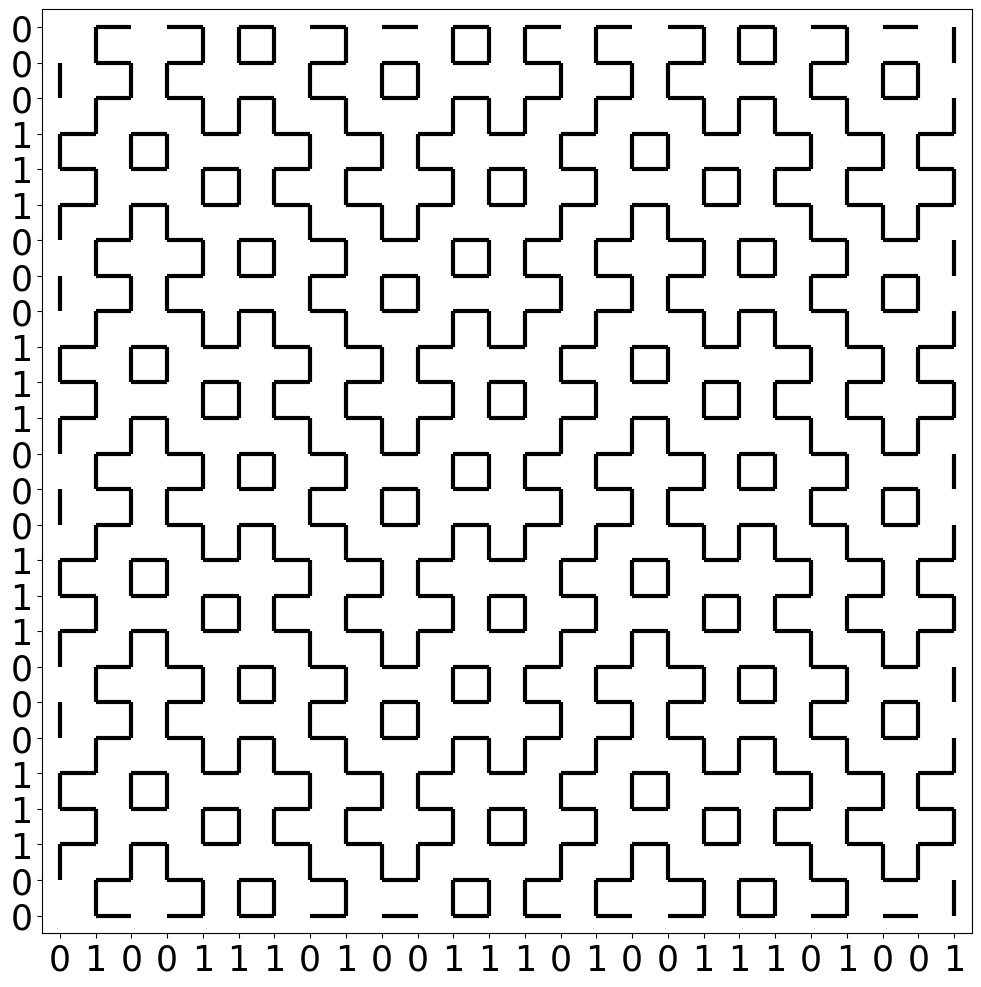}

\caption{Type $pg/pg$ where $x=(1011000)^{\infty}$ and $y=(110001)^{\infty}$.}
\label{fig:pg/pg}
\end{figure}

\begin{figure}[h!tbp] 
\centering

\includegraphics[width=4.5cm]{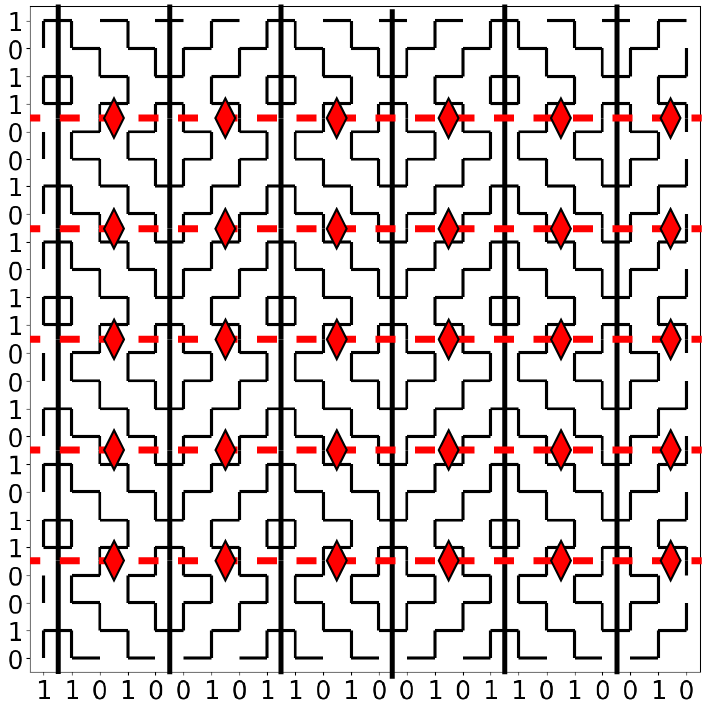}
\includegraphics[width=4.5cm]{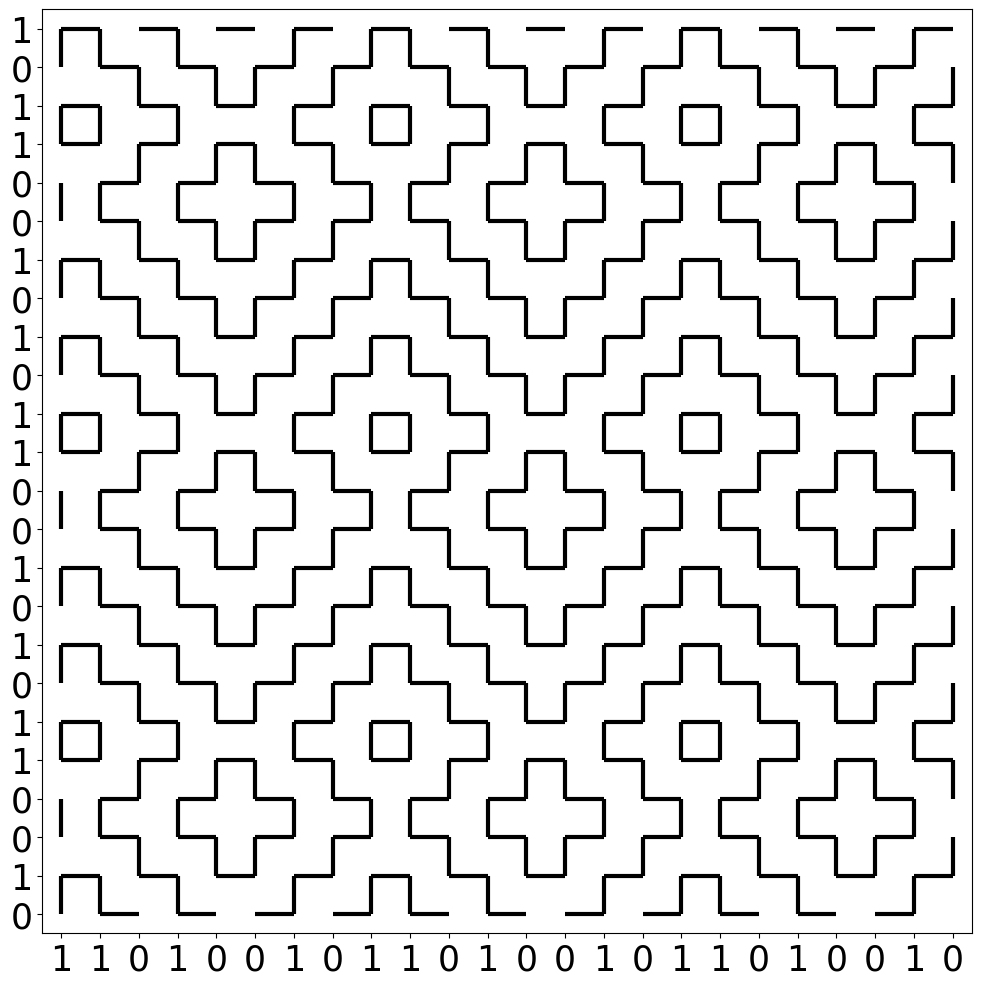}
\includegraphics[width=4.5cm]{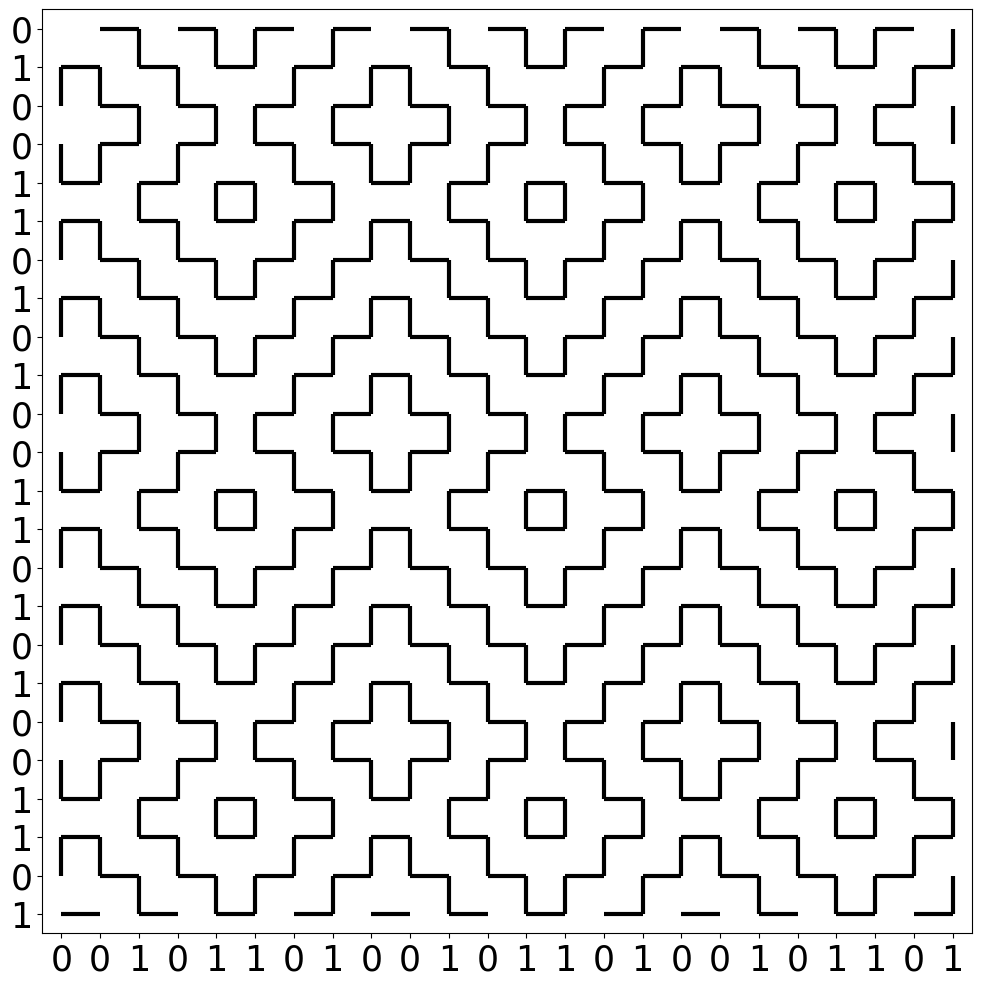}

\caption{Type $pmg/pm$ where $x=(11010010)^{\infty}$ and $y=(01001101)^{\infty}$.}
\label{fig:pmg/pm}
\end{figure}

\begin{figure}[h!tbp] 
\centering

\includegraphics[width=4.5cm]{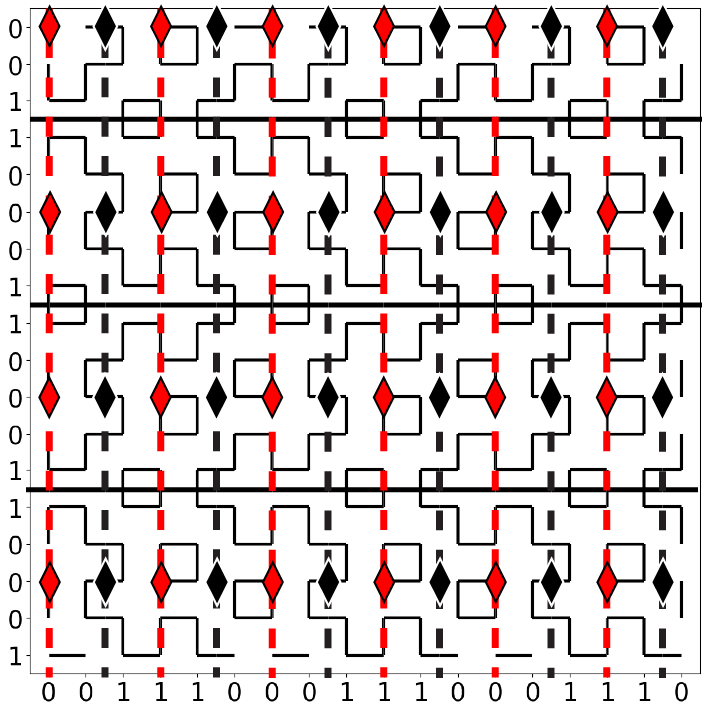}
\includegraphics[width=4.5cm]{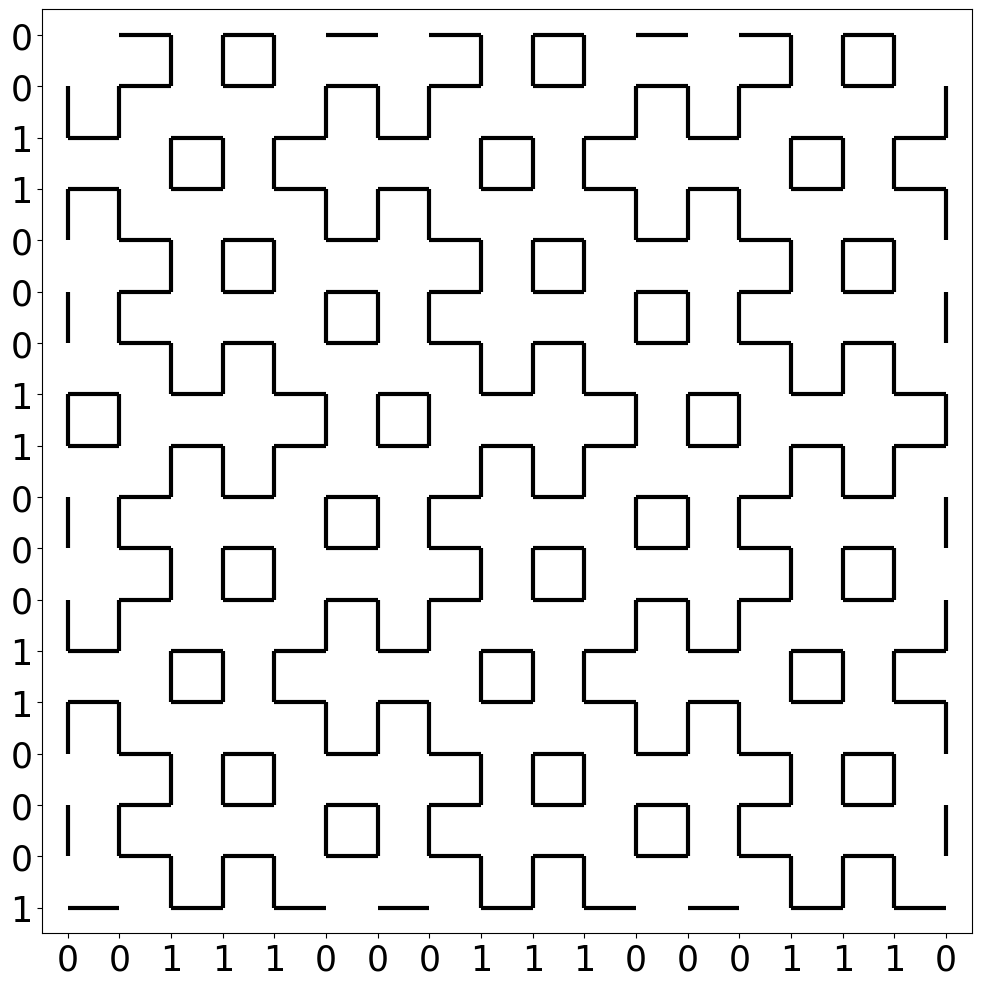}
\includegraphics[width=4.5cm]{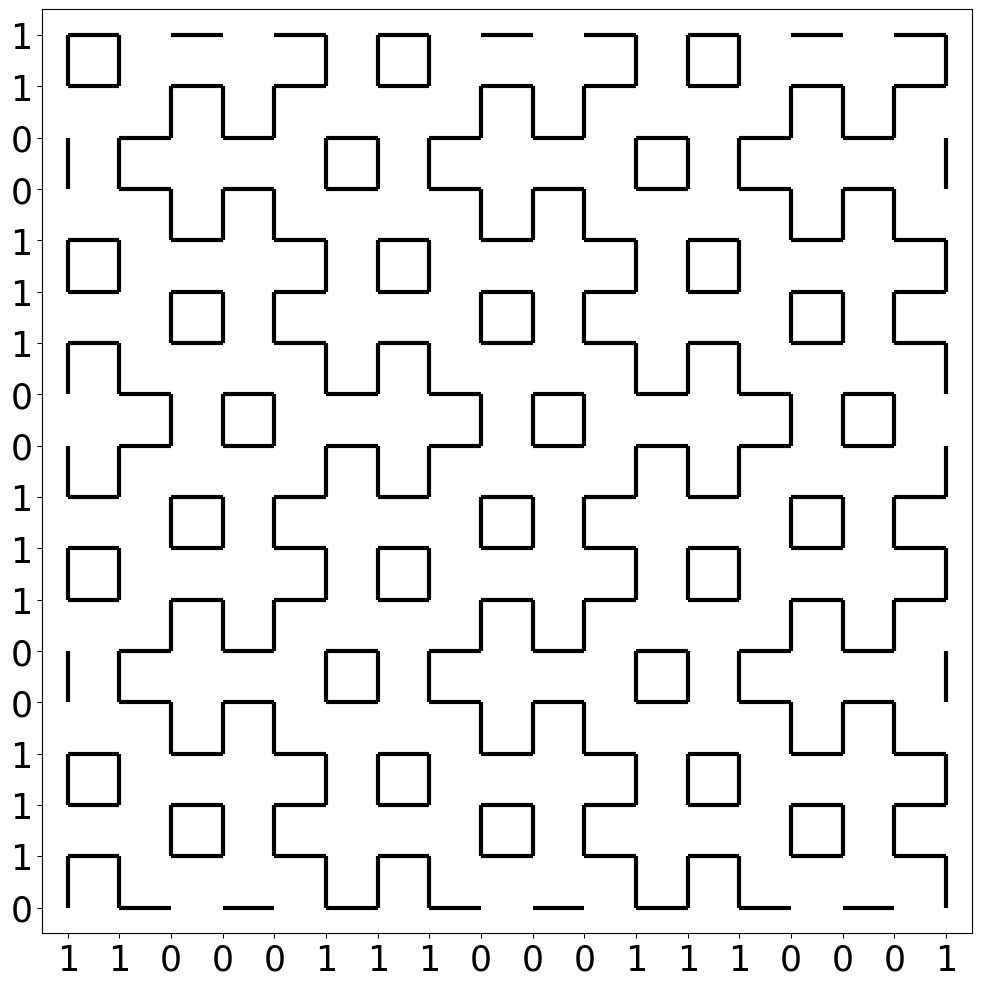}

\caption{Type $pmg/pmg$ where $x=(001110)^{\infty}$ and $y=(10001)^{\infty}$.}
\label{fig:pmg/pmg}
\end{figure}

\begin{figure}[h!tbp] 
\centering

\includegraphics[width=4.5cm]{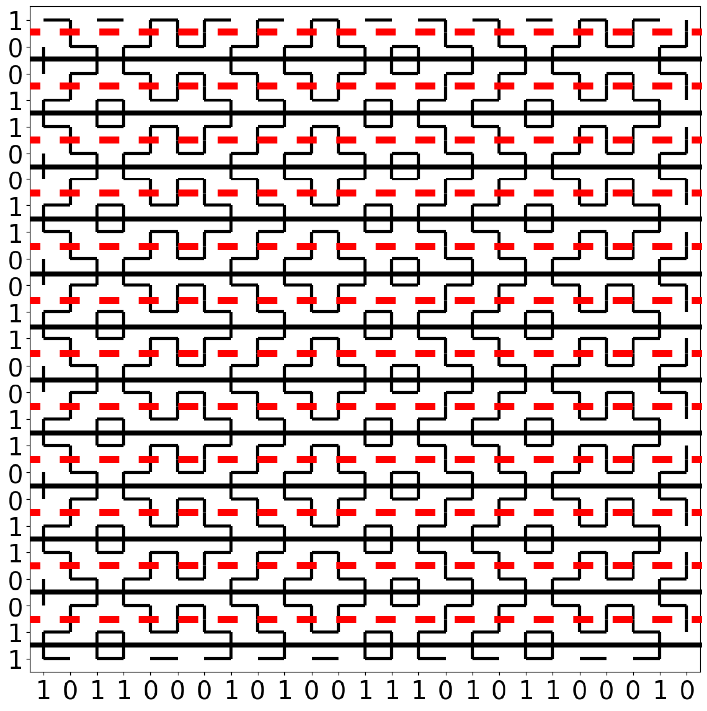}
\includegraphics[width=4.5cm]{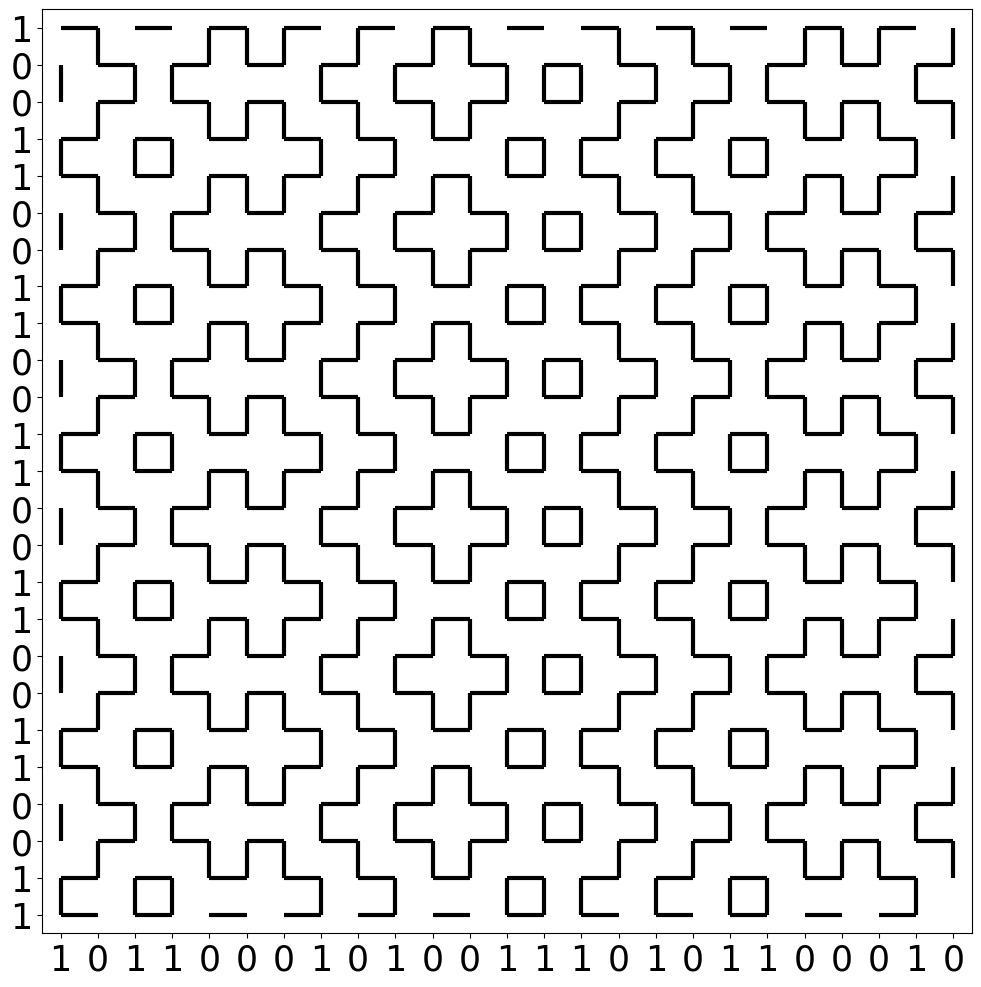}
\includegraphics[width=4.5cm]{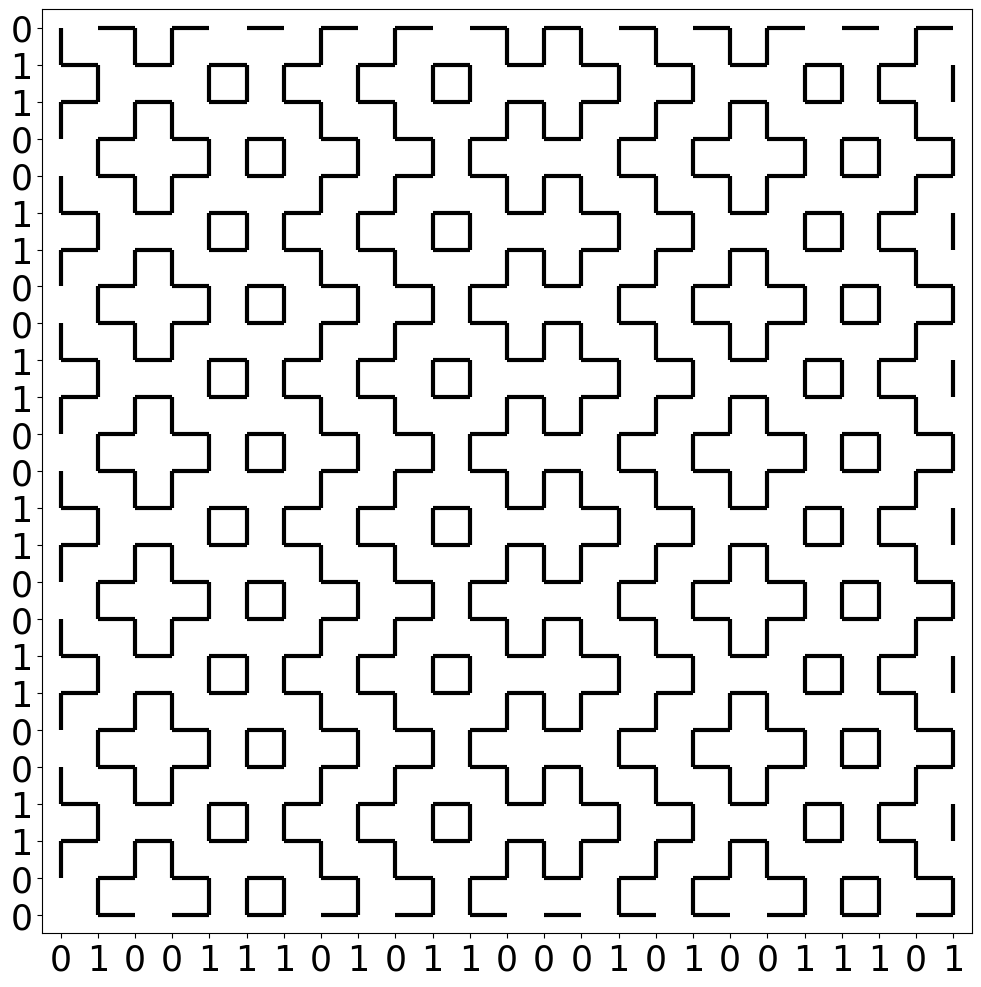}

\caption{Type $cm/pm$ where $x=(1011000101001110)^{\infty}$ and $y=(1100)^{\infty}$.}
\label{fig:cm/pm}
\end{figure}

\begin{figure}[h!tbp] 
\centering

\includegraphics[width=4.5cm]{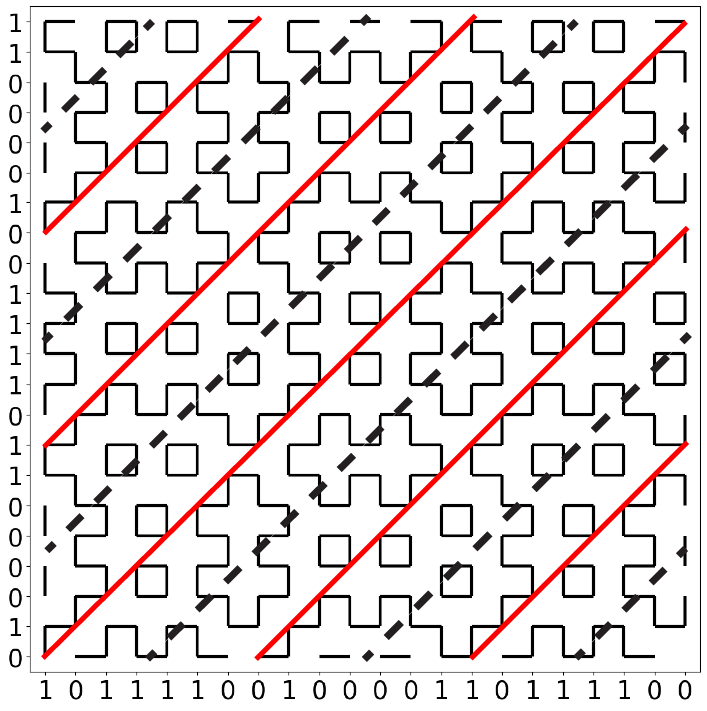}
\includegraphics[width=4.5cm]{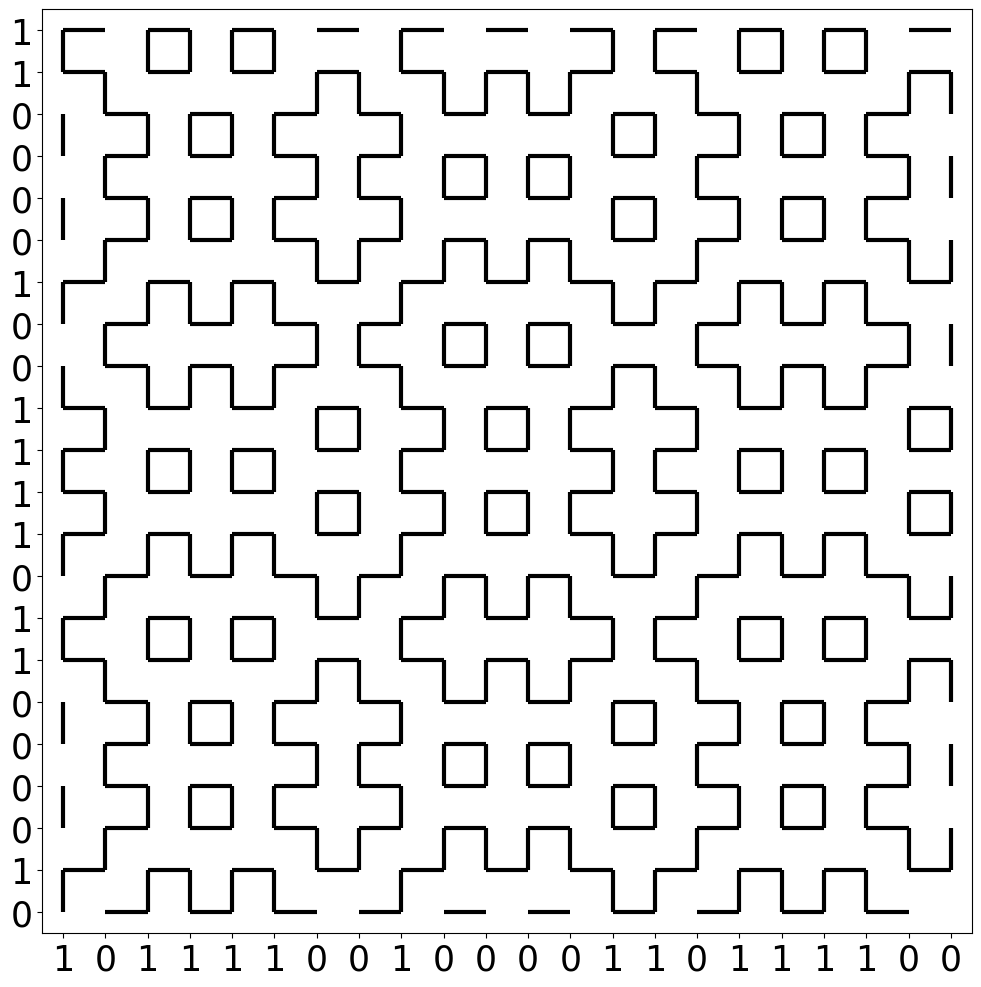}
\includegraphics[width=4.5cm]{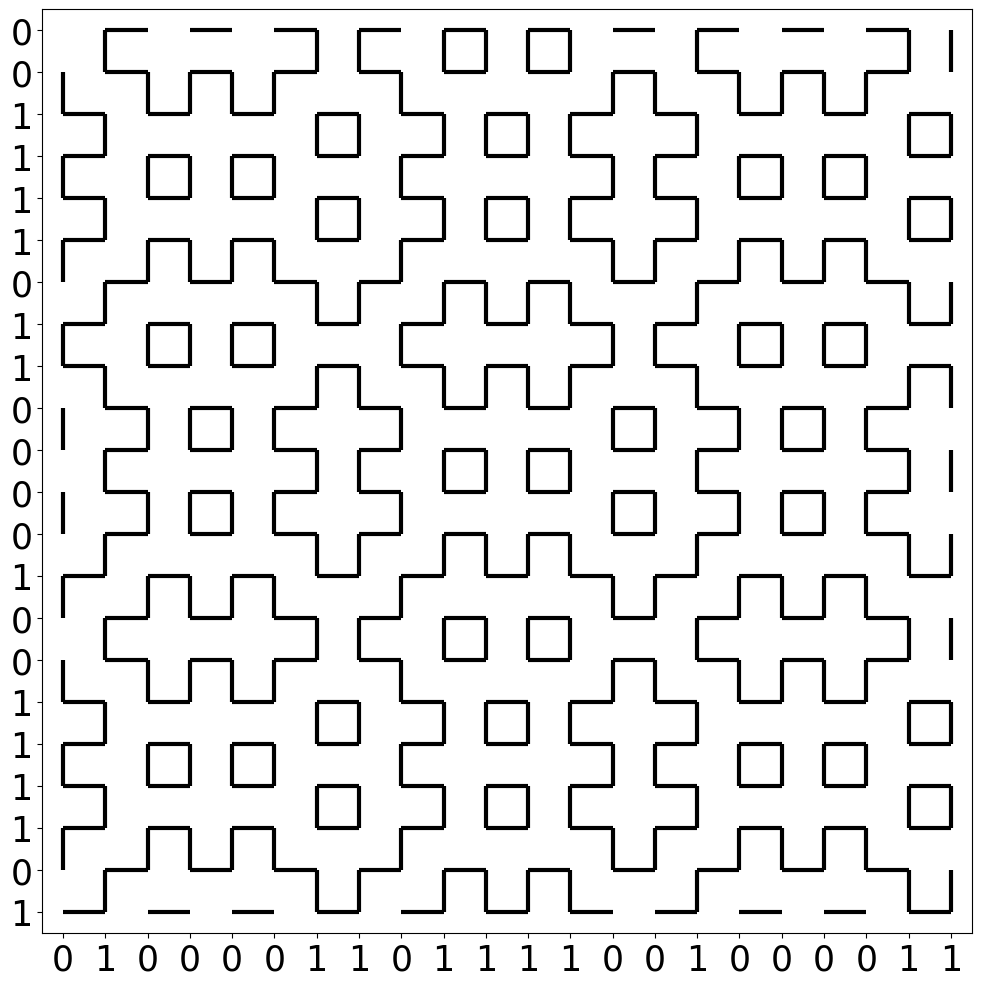}

\caption{Type $cm/pg$ where $x=(10111100100001)^{\infty}$ and $y=(01000011011110)^{\infty}$.}
\label{fig:cm/pg}
\end{figure}

\begin{figure}[h!tbp] 
\centering
\includegraphics[width=4.5cm]{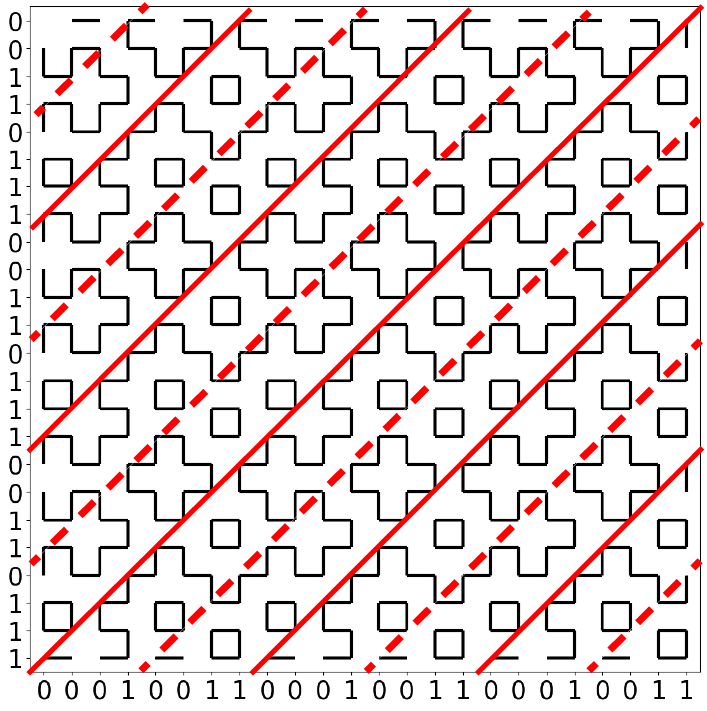}
\includegraphics[width=4.5cm]{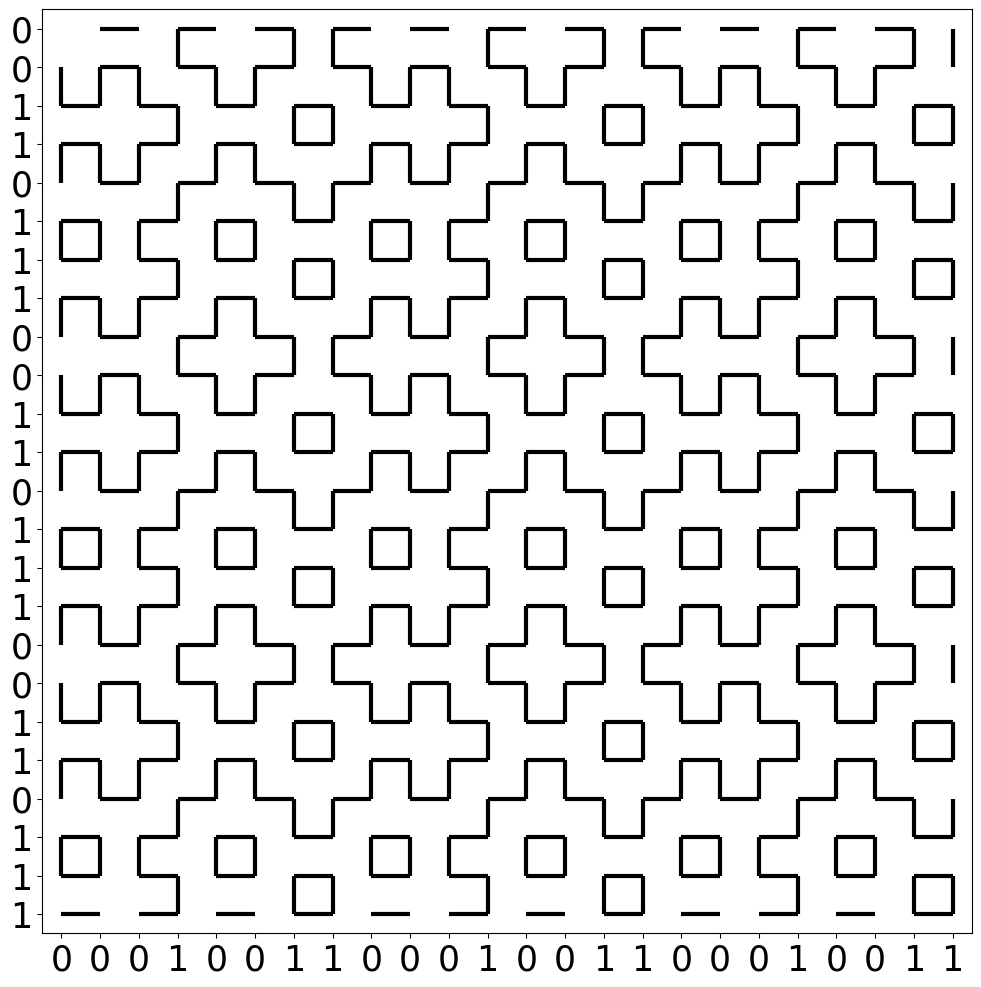}
\includegraphics[width=4.5cm]{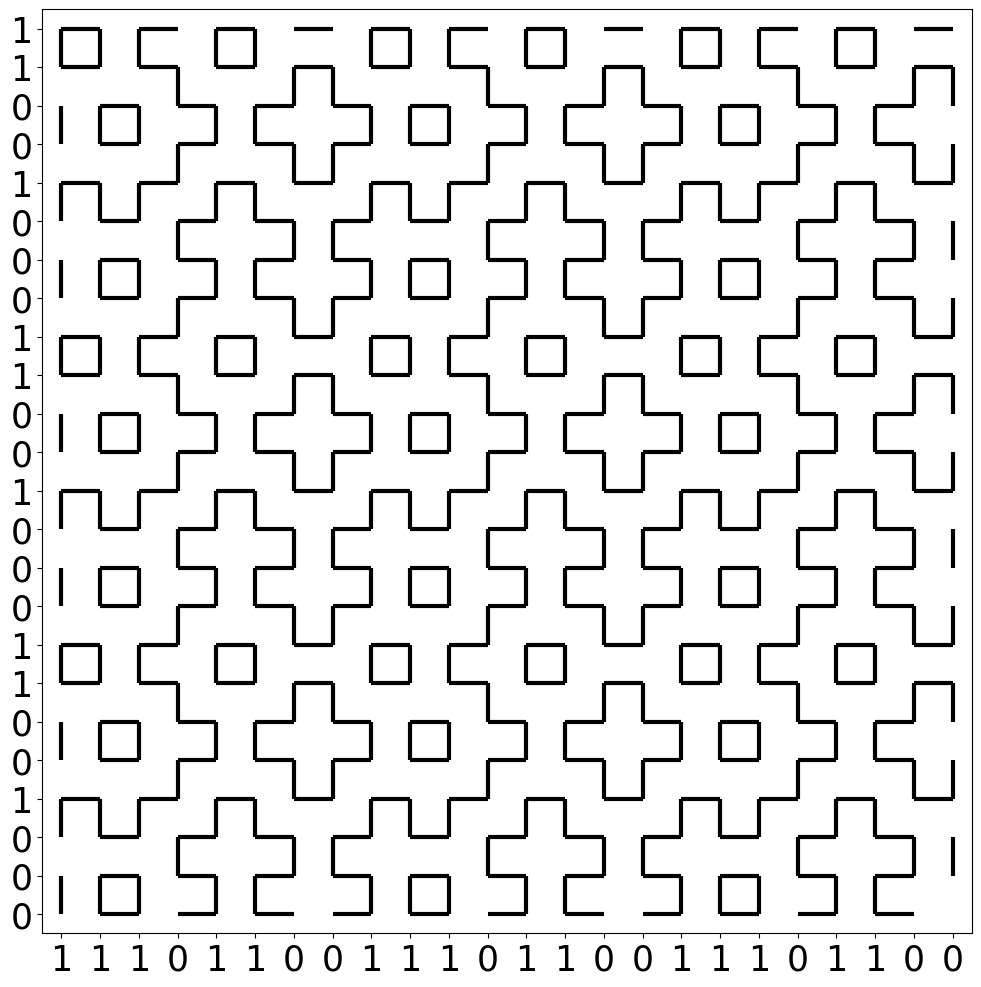}
\caption{Type $cm/p1$ where $x=(00010011)^{\infty}$ and $y=(11101100)^{\infty}$.}
\label{fig:cm/p1}
\end{figure}

\begin{figure}[h!tbp] 
\centering
\includegraphics[width=4.5cm]{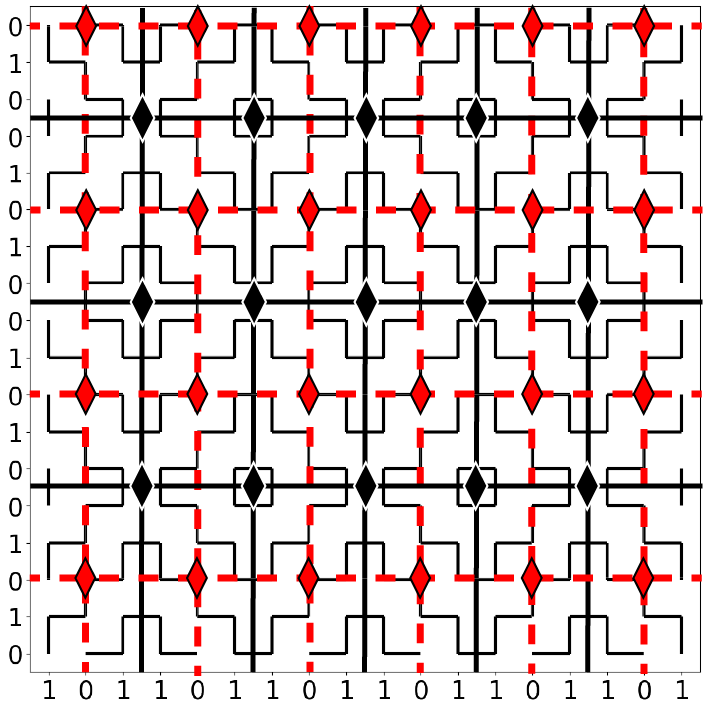}
\includegraphics[width=4.5cm]{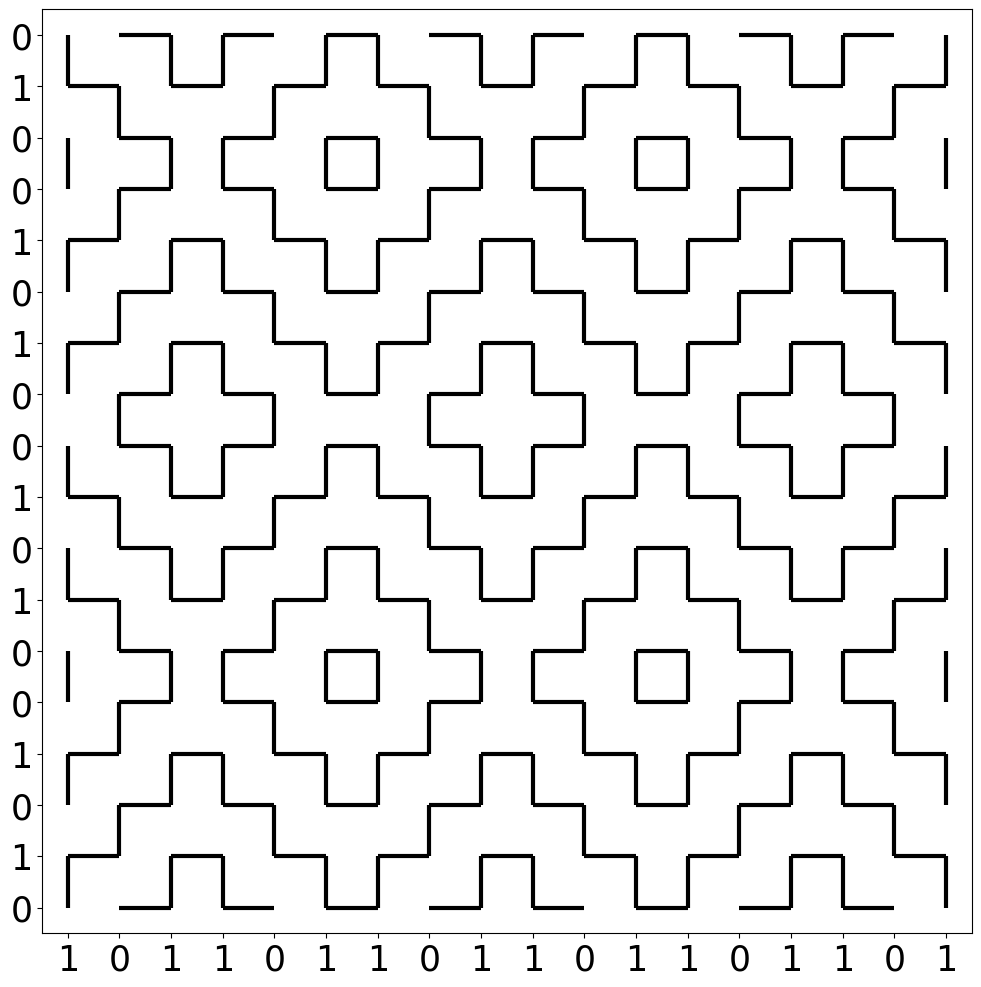}
\includegraphics[width=4.5cm]{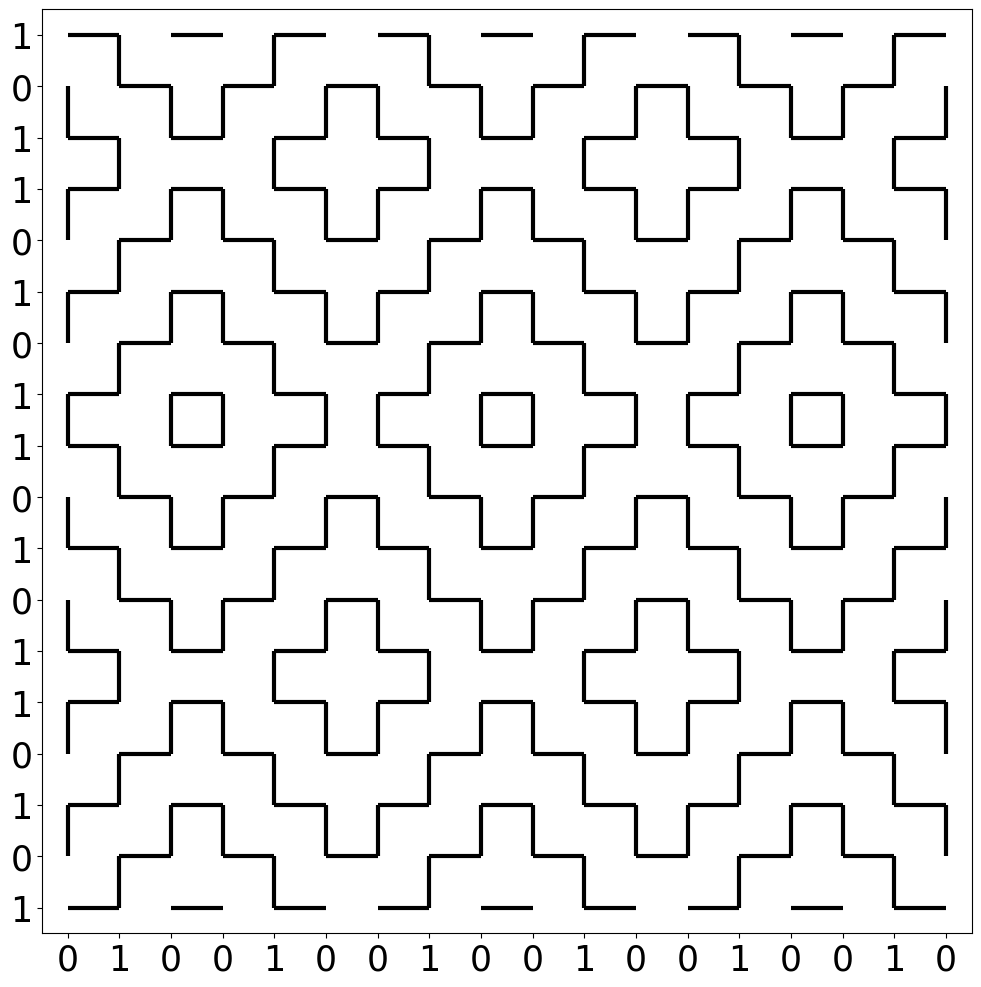}
\caption{Type $cmm/pmm$ where $x=(101101)^{\infty}$ and $y=(01010)^{\infty}$.}
\label{fig:cmm/pmm}
\end{figure}

\begin{figure}[h!tbp] 
\centering
\includegraphics[width=4.5cm]{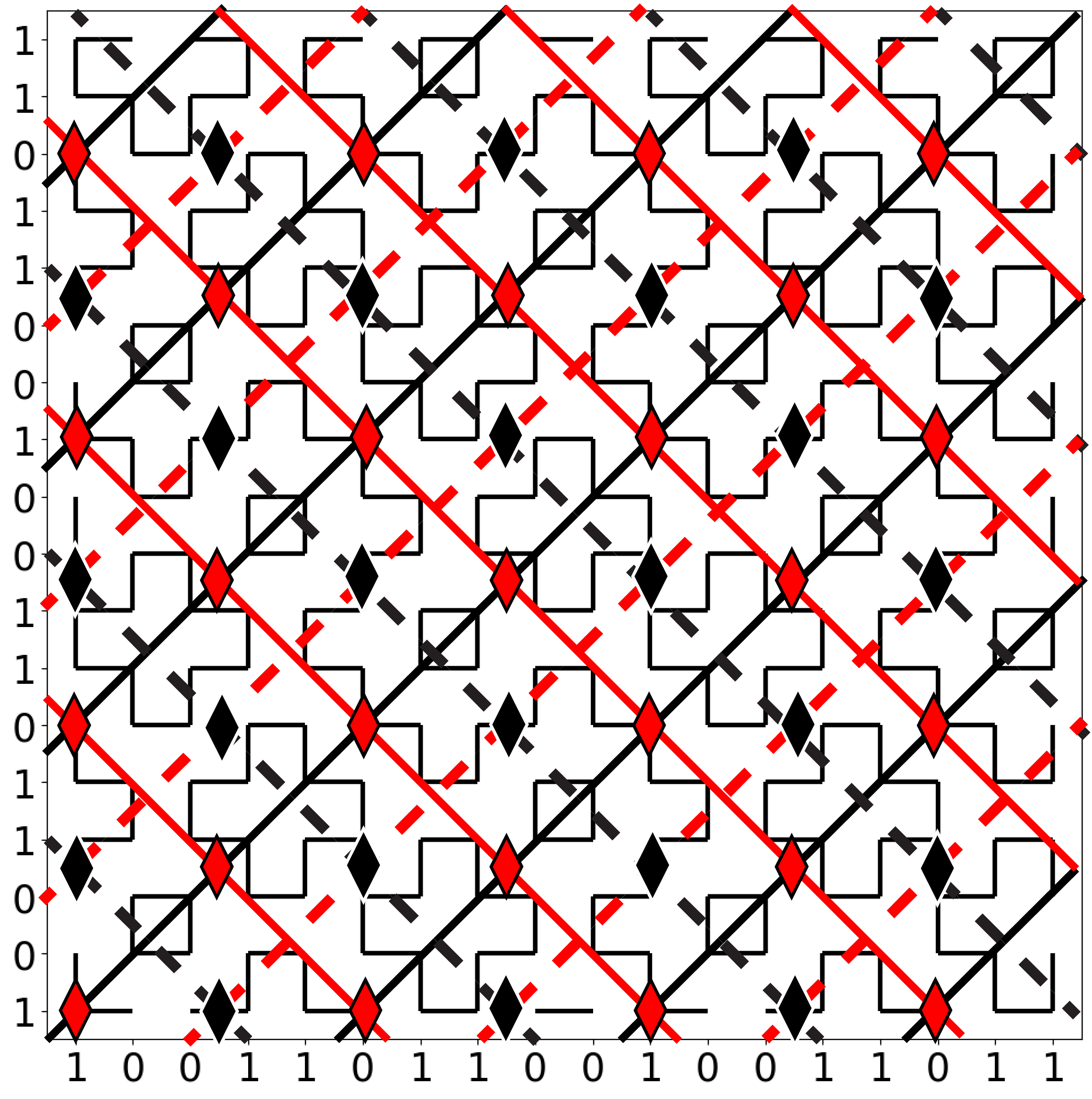}
\includegraphics[width=4.5cm]{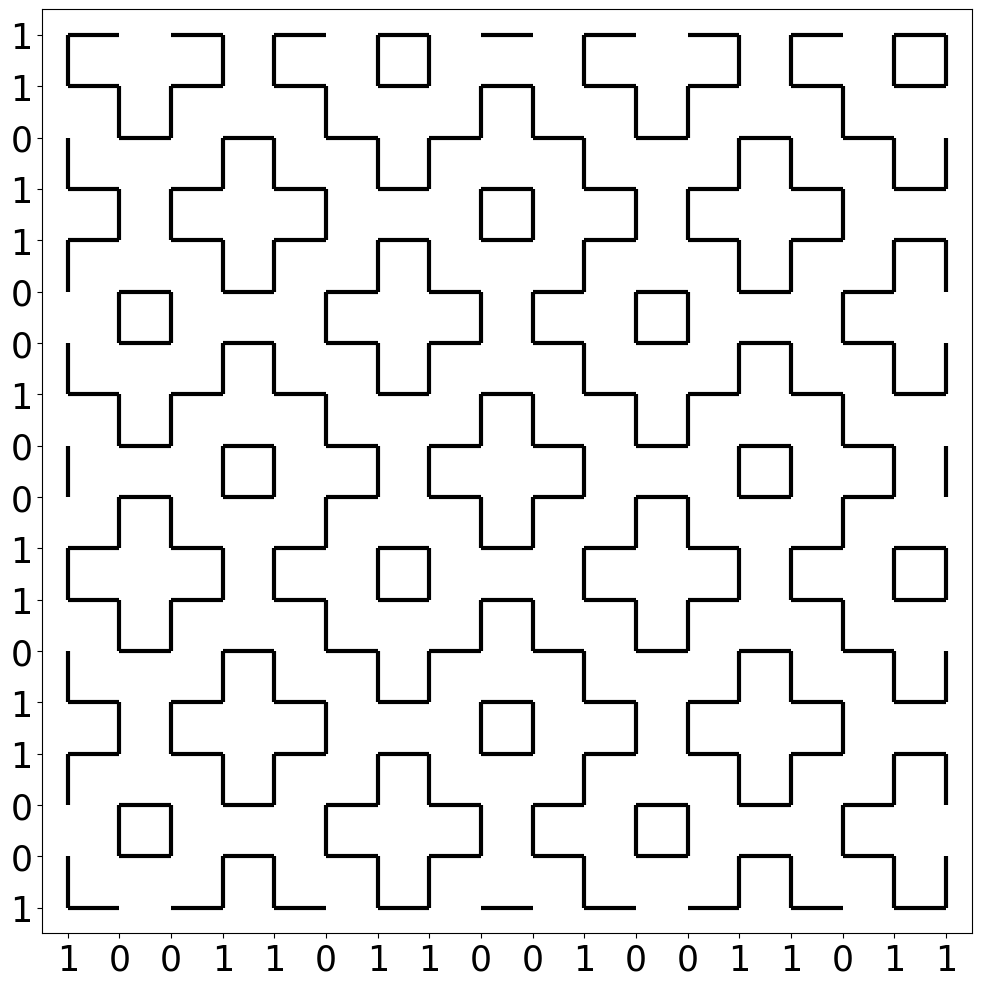}
\includegraphics[width=4.5cm]{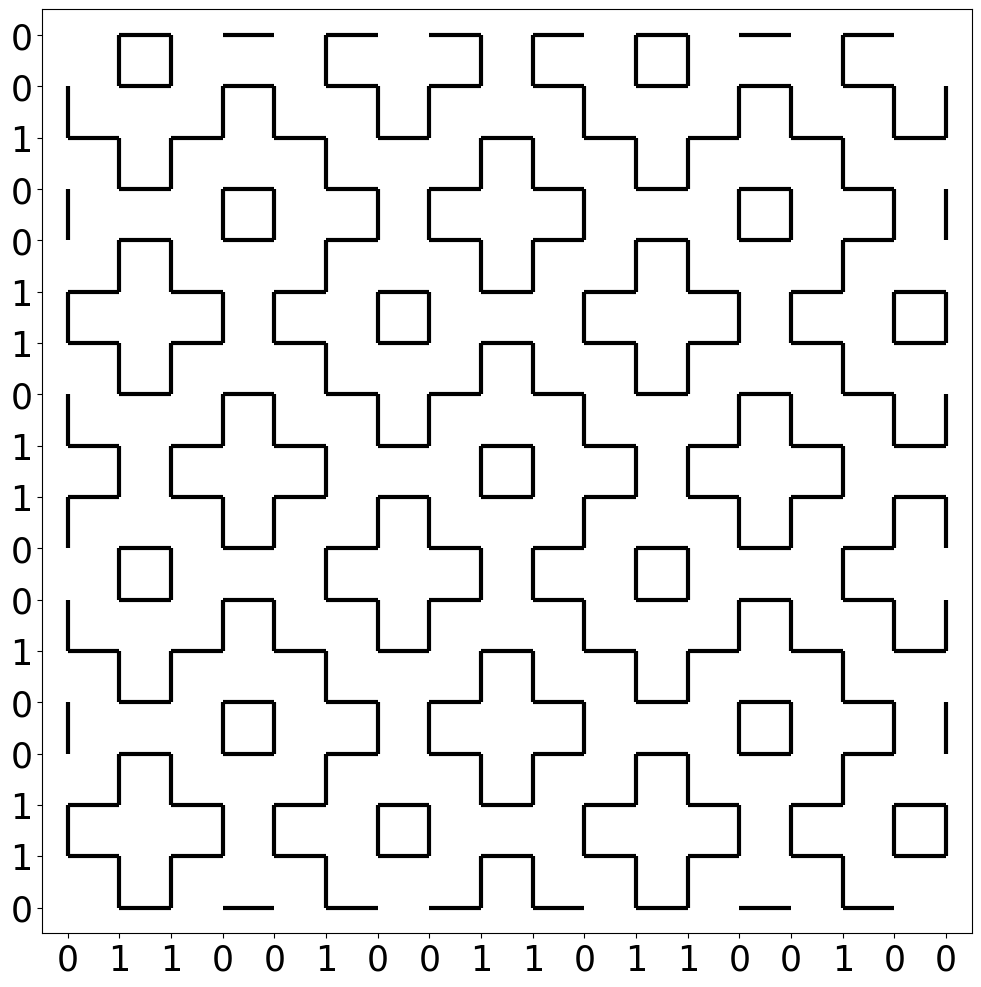}
\caption{Type $cmm/pmg$ where $x=(1001101)^{\infty}$ and $y=(1001101)^{\infty}$.}
\label{fig:cmm/pmg}
\end{figure}

\begin{figure}[h!tbp] 
\centering
\includegraphics[width=4.5cm]{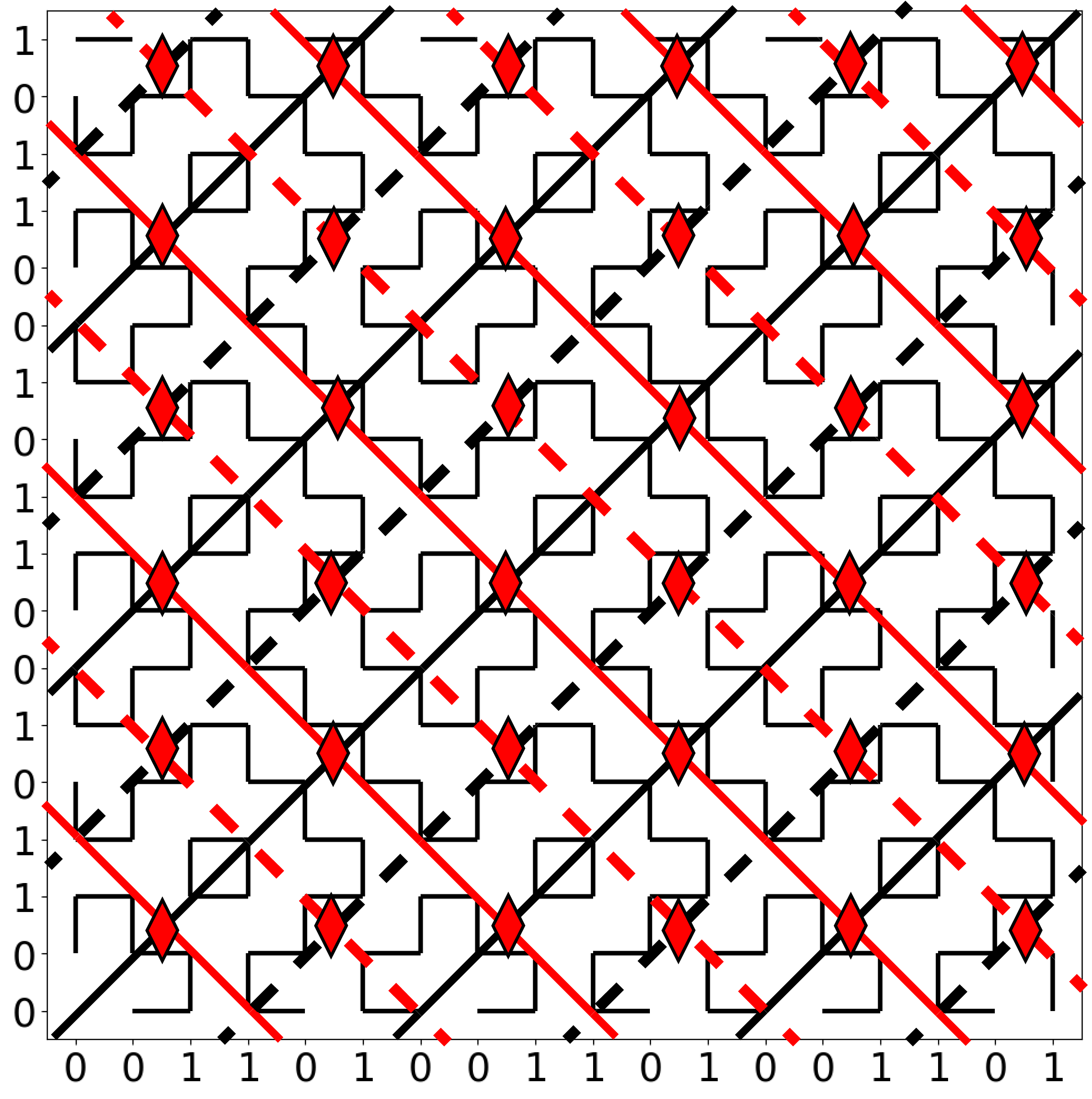}
\includegraphics[width=4.5cm]{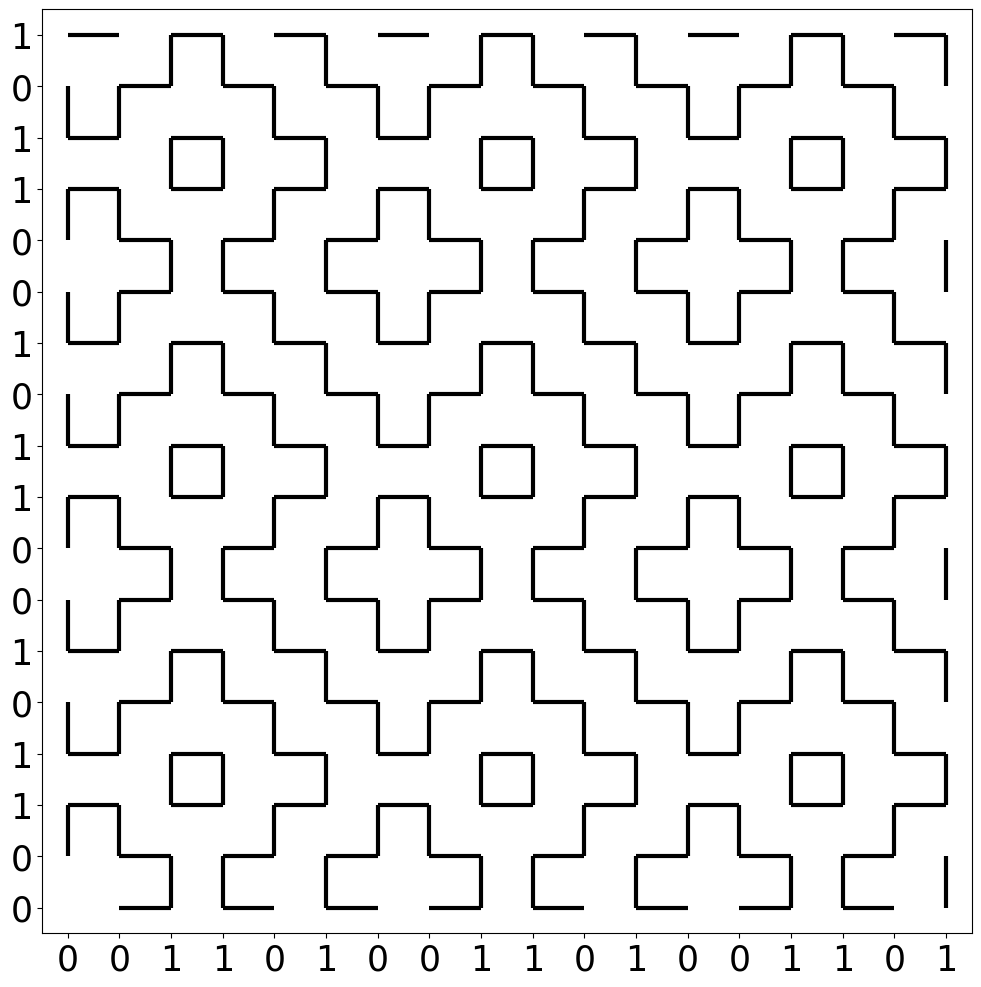}
\includegraphics[width=4.5cm]{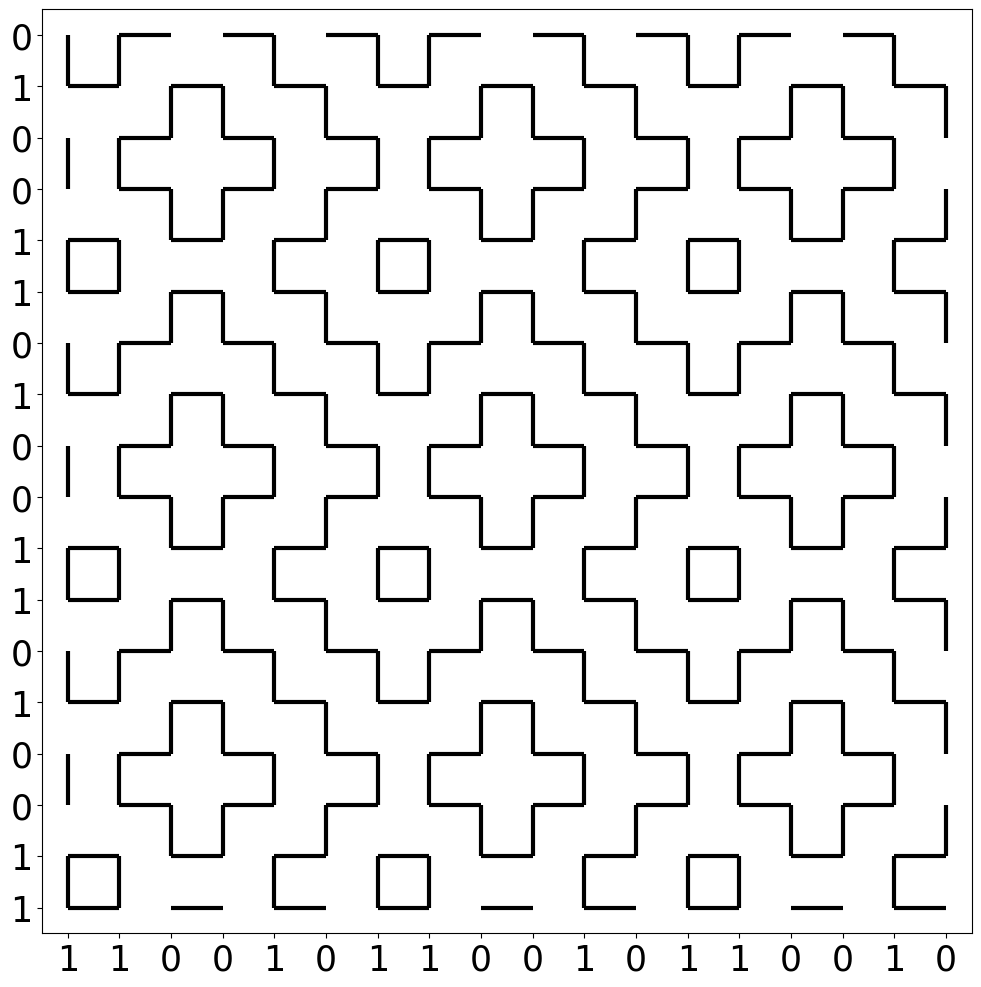}
\caption{Type $cmm/cm$ where $x=(001101)^{\infty}$ and $y=(001101)^{\infty}$.}
\label{fig:cmm/cm}
\end{figure}

\begin{figure}[h!tbp] 
\centering
\includegraphics[width=4.5cm]{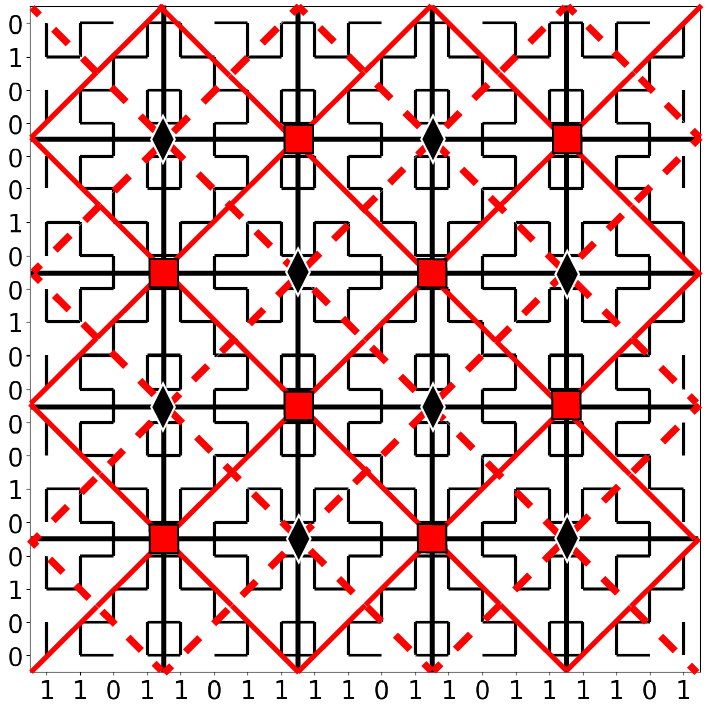}
\includegraphics[width=4.5cm]{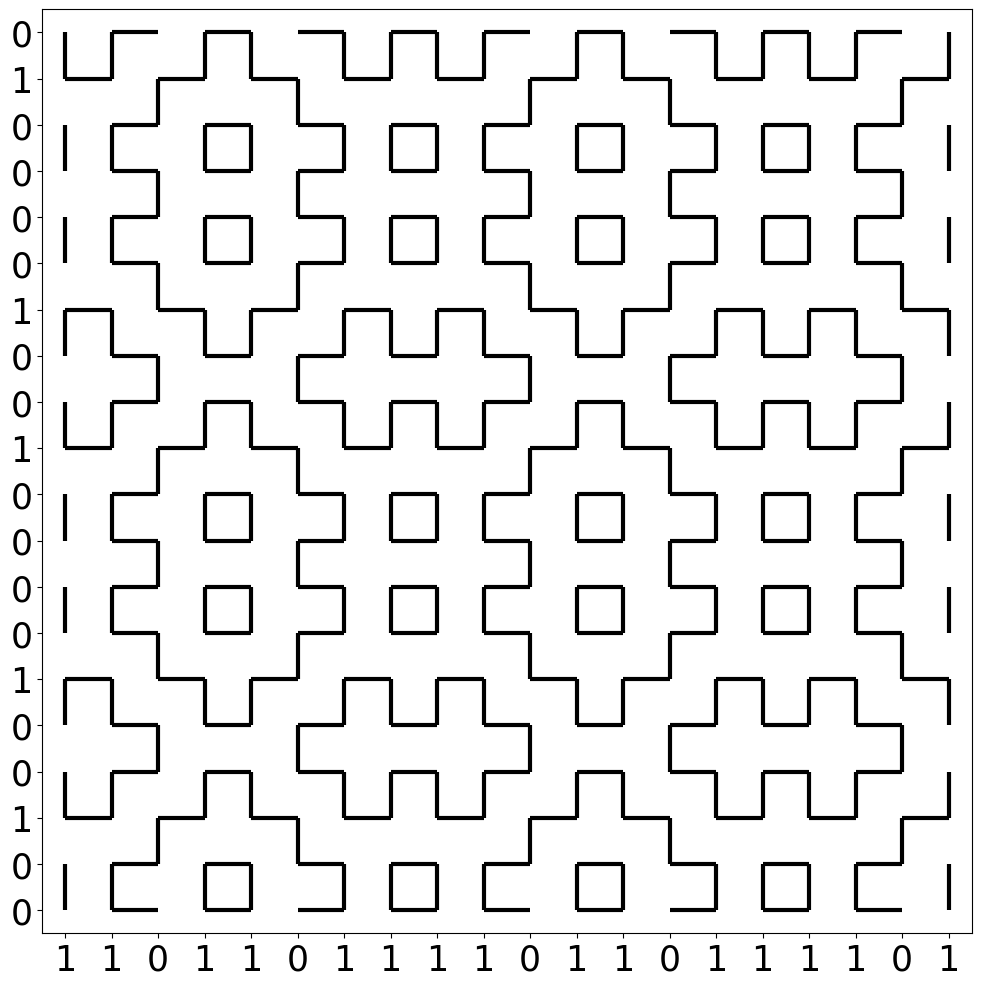}
\includegraphics[width=4.5cm]{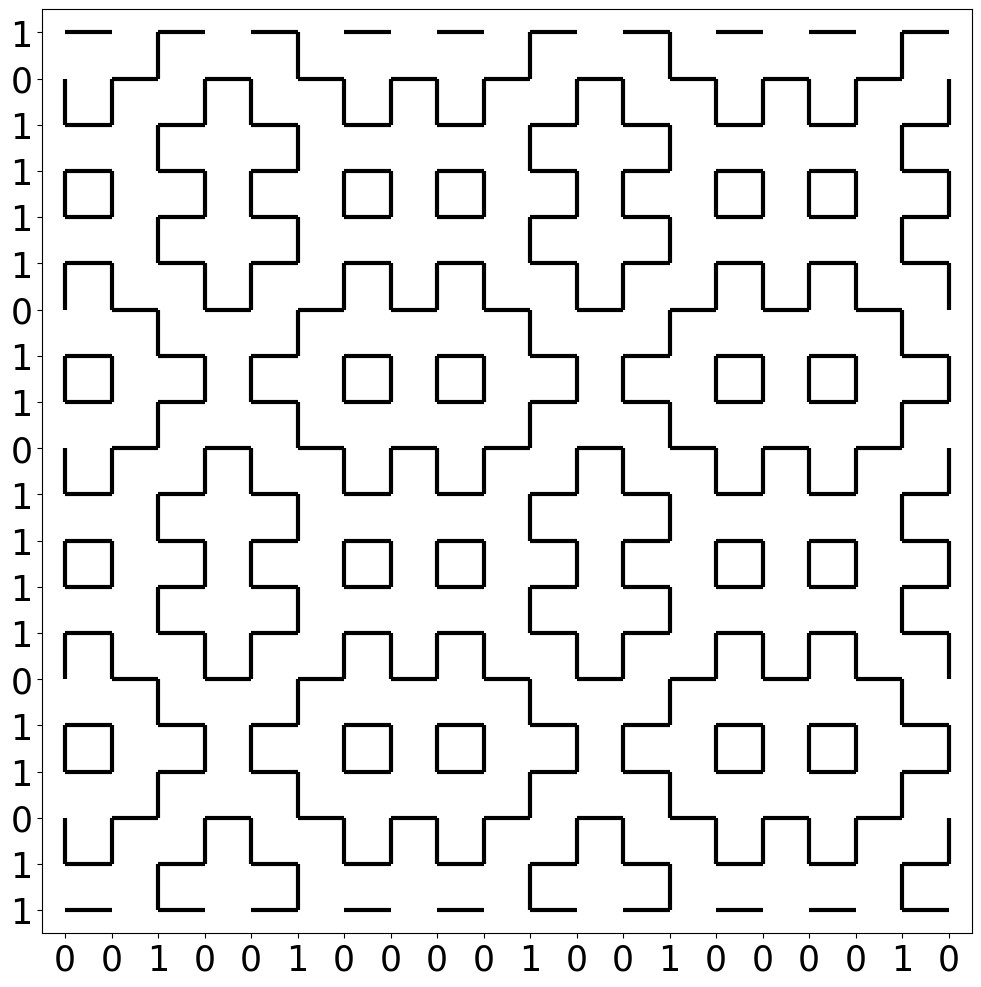}
\caption{Type $p4m/pmm$ where $x=(11011011)^{\infty}$ and $y=(00100100)^{\infty}$.}
\label{fig:p4m/pmm}
\end{figure}

\begin{figure}[h!tbp] 
\centering
\includegraphics[width=4.5cm]{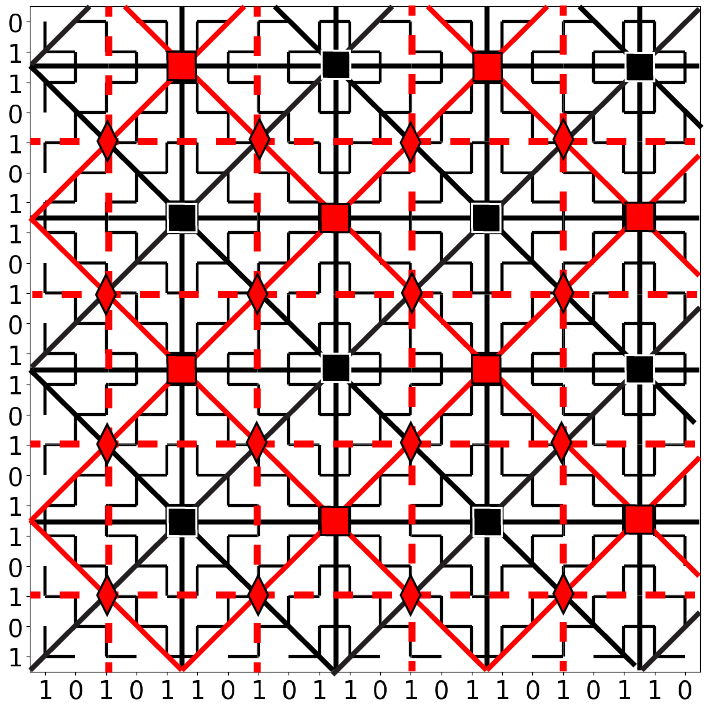}
\includegraphics[width=4.5cm]{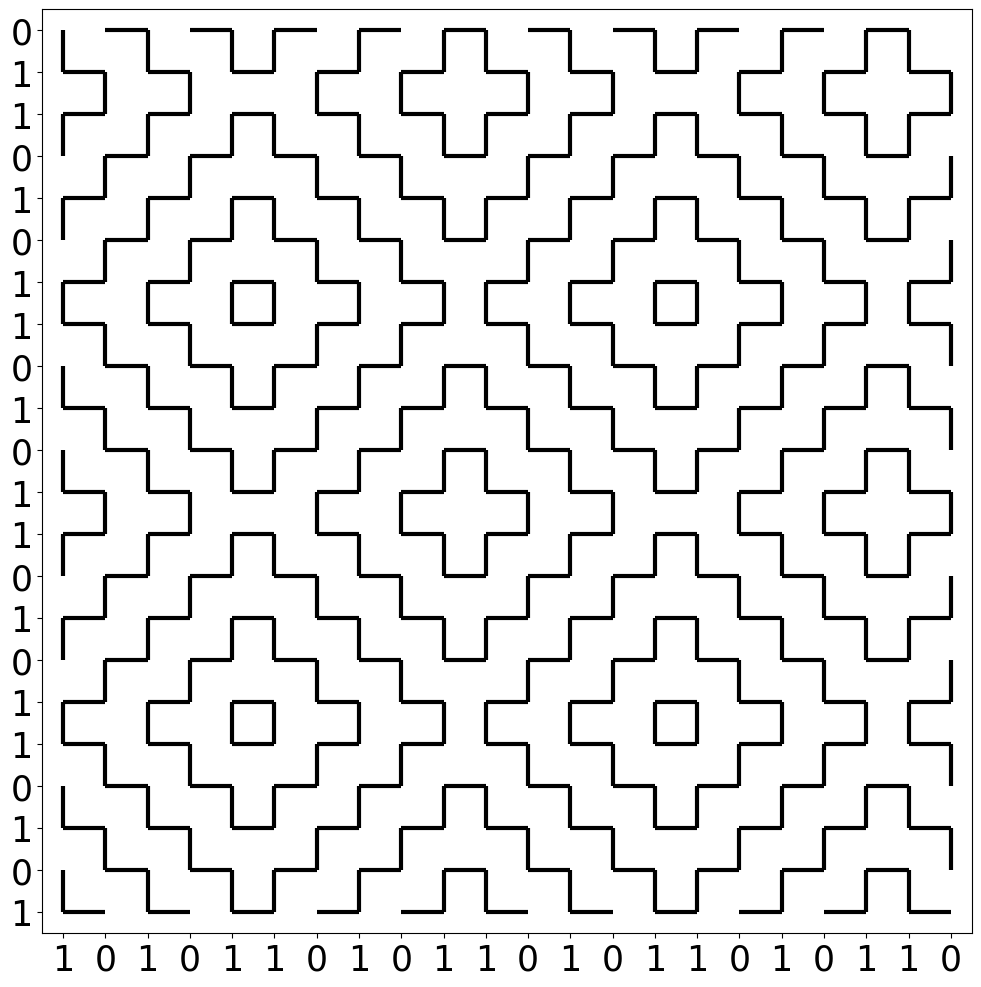}
\includegraphics[width=4.5cm]{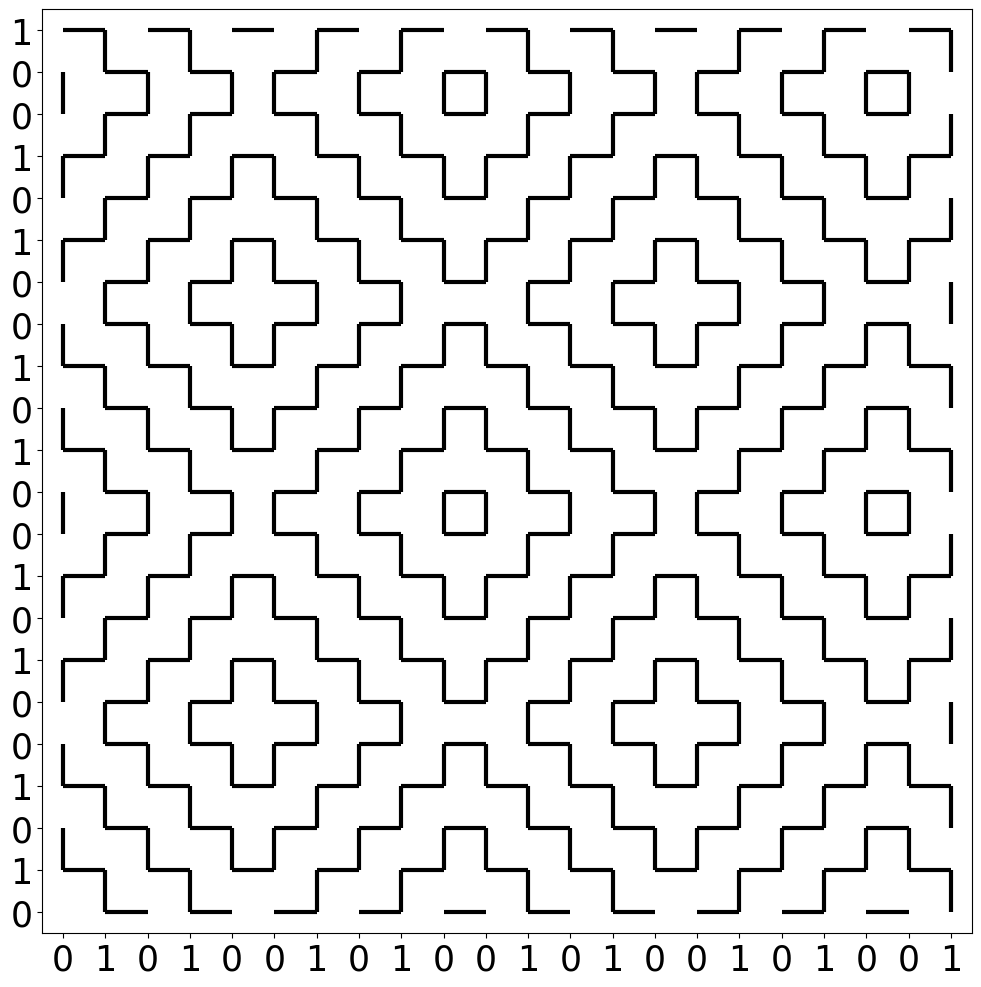}
\caption{Type $p4m/p4m$ where $x=(10101)^{\infty}$ and $y=(10101)^{\infty}$.}
\label{fig:p4m/p4m}
\end{figure}

Table~\ref{table:wallpapers} provides the list of all 46 two-colour wallpaper groups and whether they are compatible or incompatible with PGHPs. When the type is incompatible, the corresponding lemma proving that fact is referenced. For the types that are compatible, they have been illustrated in Figures~\ref{fig:p1/p1} through \ref{fig:p4m/p4m}. In each of these figures, the first image shows the front stitches with symmetry markings, the second shows just the front stitches, and the third shows that back stitches as they appear as if looking through transparent fabric.

\begin{table}[h!tbp]
\centering
\begin{tabular}{|c|c||c|c|} \hline
2-colour &  & 2-colour &  \\
Wallpaper Group & Compatible? & Wallpaper Group & Compatible?\\\hline \hline
$p1/p1$ & Yes & $cm/pg$ & Yes \\ \hline
$p2/p1$ & Yes & $cm/p1$ & Yes \\ \hline
$p2/p2$ & Yes & $cmm/pgg$ & No (Lemma~\ref{lem:one_color_results}) \\ \hline
$pm/p1$ & No (Lemma~\ref{lem:flip_reflection_lemma_results}) & $cmm/pmm$ & Yes \\ \hline
$pm/pm(m')$ & No (Lemma~\ref{lem:flip_reflection_lemma_results}) & $cmm/pmg$ & Yes \\ \hline
$pm/pg$ & No (Lemma~\ref{lem:flip_reflection_lemma_results}) & $cmm/cm$ & Yes \\ \hline
$pm/pm(m)$ & Yes & $cmm/p2$ & No (Lemma~\ref{lem:d'2}) \\ \hline
$pm/cm$ & Yes & $p4/p2$ & No (Lemma~\ref{lem:c4'}) \\ \hline
$pg/p1$ & Yes & $p4/p4$ & No (Lemma~\ref{lem:one_color_results}) \\ \hline
$pg/pg$ & Yes & $p4m/cmm$ & No (Lemma~\ref{lem:d'2}) \\ \hline
$pmm/pm$ & No (Lemma~\ref{lem:flip_reflection_lemma_results}) & $p4m/pmm$ & Yes \\ \hline
$pmm/p2$ & No (Lemma~\ref{lem:d'2}) & $p4m/p4$ & No (Lemma~\ref{lem:one_color_results}) \\ \hline
$pmm/pmg$ & No (Lemma~\ref{lem:d'2}) & $p4m/p4m$ & Yes\\ \hline
$pmm/pmm$ & No (Lemma~\ref{lem:pmmpmm}) & $p4m/p4g$ & No (Lemma~\ref{lem:one_color_results})\\ \hline
$pmm/cmm$ & No (Lemma~\ref{lem:d'2}) & $p4g/pgg$ & No (Lemma~\ref{lem:one_color_results})\\ \hline
$pmg/pm$ & Yes & $p4g/cmm$ & No (Lemma~\ref{lem:c4'})\\ \hline
$pmg/pg$ & No (Lemma~\ref{diagonal_glide_ref_flip}) & $p4g/p4$ & No (Lemma~\ref{lem:one_color_results})\\ \hline
$pmg/p2$ & No (Lemma~\ref{diagonal_glide_ref_flip}) & $p31m/p3$ & No (Lemma~\ref{lem:3_fold_wallpaper})  \\ \hline
$pmg/pgg$ & No (Lemma~\ref{lem:one_color_results}) & $p3m1/p3$ & No (Lemma~\ref{lem:3_fold_wallpaper})  \\ \hline
$pmg/pmg$ & Yes & $p6/p3$ & No (Lemma~\ref{lem:3_fold_wallpaper})  \\ \hline
$pgg/pg$ & No (Lemma \ref{lem:pgg/pg}) & $p6m/p3m1$ & No (Lemma~\ref{lem:3_fold_wallpaper})  \\ \hline
$pgg/p2$ & No (Lemma \ref{lem:pgg/pg}) & $p6m/p31m$ & No (Lemma~\ref{lem:3_fold_wallpaper}) \\ \hline
$cm/pm$ & Yes & $p6m/p6$ & No (Lemma~\ref{lem:3_fold_wallpaper}) \\ \hline
\end{tabular}
\caption{A summary of whether or not each of the 46 two-colour wallpaper groups are compatible with PGHPs. When the wallpaper group is incompatible with PGHPs, the lemma explaining that fact is referenced.}
\label{table:wallpapers}
\end{table}

\subsection{Examples of PGHPs}

This section contains example PGHPs for each of the compatible 2-colour wallpaper groups. All the computer-generated hitomezashi diagrams in this section (and indeed, the entire paper) were generated using the Python program written by Antara Sen.

\newpage

\section{Conclusion}
This paper is a full accounting of the way dual symmetries can be presented in both GHPs and PGHPs. While Table~\ref{rosettes} gives explicit instructions for how to generate the possible symmetries, we have left out descriptions on what binary strings should be chosen to produce particular symmetries in a PGHP. When constructing examples of different wallpaper symmetries, Lemmas~\ref{diagonal_glide2}, \ref{diagonal_glide3}, and \ref{diagonal_reflection_flip} provide clues in how to produce certain kinds of glide-reflective and reflective symmetry. Table~\ref{rosettes} can be leveraged as well to create certain symmetries in a patch of the PGHP---the difficulty in applying this information to PGHPs is that unintended symmetries often present in the periodic design that were not in the patch. However, through trial and error, it is possible to collate the information in this paper to improve your chances of stumbling upon a chosen set of wallpaper symmetries. The lessons of this paper should provide a useful starting point for anyone wanting to generate their own hitomezashi-inspired art.

Hitomezashi designs are not limited to the kinds described in this paper. There is ample opportunity to explore designs on different grids and with different stitch patterns. We hope the readers will use the tools developed here to fuel their own explorations into the beauty of hitomezashi.

\bibliographystyle{apacite}
\bibliography{ref.bib}

\end{document}